# The Seiberg-Witten equations and the Weinstein conjecture II: More closed integral curves for the Reeb vector field

Clifford Henry Taubes[†]

Department of Mathematics
Harvard University
Cambridge MA  02133

chtaubes@math.harvard.edu

Let M denote a compact, orientable 3-dimensional manifold and let a denote a contact 1-form on M; thus $a \wedge da$ is nowhere zero.  This article explains how the Seiberg-Witten Floer homology groups as defined for any given $Spin_{\mathbb{C}}$ structure on M give closed, integral curves of the vector field that generates the kernel of da.

[†]Supported in part by the National Science Foundation

## 1. Introduction

Let M denote a smooth, compact 3-manifold and let a denote a smooth contact 1-form on M. This is to say that the 3-form $a \wedge da$ is nowhere zero. Let v denote the vector field that generates the kernel of da, and pairs with a to give 1. The prequel to this article, [T1], asserts that v has closed, integral curves. To be more specific, [T1] asserted the following:

> *Let $K^{-1}$ denote the 2-plane bundle kernel(a) $\subset$ TM, and let $e \in H^2(M; \mathbb{Z})$ denote any class that differs by a torsion class from half the Euler class of K. There exists a finite, non-empty set of closed, integral curves of v and a corresponding set of positive weights such that the formal sum of weighted loops in M generates the Poincaré dual to e.*

(1.1)

Here, M is implicitly oriented by $a \wedge da$, and K is oriented so that the isomorphism $T^*M = \mathbb{R}v \oplus K^{-1}$ is orientation preserving. Note that the Euler class of K is, in all cases, divisible by two.

The purpose of this article is to describe other classes in $H_1(M; \mathbb{Z})$ that are represented by formal, positively weighted sums of closed integral curves of v. To set the stage for a precise statement of what is proved here, digress for a moment to introduce the set, $\mathcal{S}_M$, of $Spin_\mathbb{C}$ structures on M, this a principle homogeneous space for the group $H^2(M; \mathbb{Z})$. Each element in $\mathcal{S}$ has a canonically associated $\mathbb{Z}$-module, a version of Seiberg-Witten Floer homology that is described below. As explained momentarily, the contact form, a, supplies a canonical, $H^2(M; \mathbb{Z})$-equivariant isomormophism from $\mathcal{S}_M$ to $H^2(M; \mathbb{Z})$. With this isomorphism understood, each element in $H^2(M; \mathbb{Z})$ has an associated Seiberg-Witten Floer homology $\mathbb{Z}$-module.

**Theorem 1.1**: *Let $e \in H^2(M; \mathbb{Z})$. When $e \neq 0$, the Poincaré dual of e is represented by a positive integer weighted sum of closed, integral curves of the vector field v if the associated Seiberg-Witten Floer homology is non-trivial. The class $0 \in H_1(M; \mathbb{Z})$ is represented by a non-trivial, positive, integer weighted sum of closed integral curves of v provided the following: Either the Seiberg-Witten Floer homology for $0 \in H^2(M; \mathbb{Z})$ has two or more generators, or the Seiberg-Witten Floer homology for 0 vanishes when tensored with $\mathbb{Q}$.*

The case in Theorem 1.1 when e differs from half K's Euler class by a torsion element is the case that is considered in [T1] and (1.1). As a consequence, it is assumed in what follows, unless directed otherwise, that these elements differ by a non-torsion element.

There are various existence theorems for non-zero classes in the relevant version of Seiberg-Witten Floer homology. For example, Theorem 41.5.2 in Kronheimer and Mrowka's book [KM1] asserts an existence theorem for non-trivial Floer homology groups that can be used with Theorem 1.1 Their theorem refers to the dual Thurston norm for a certain subspace in $H^2(M; \mathbb{Z})$. This dual norm is defined using the Thurston norm on $H_2(M; \mathbb{Z})$. The latter is a semi-norm that assigns to a non-zero class h the infimum of the absolute values of the Euler characteristics of embedded surfaces in M with fundamental class equal to h in $H_2(M; \mathbb{Z})$. The dual of this semi-norm defines a norm on the subspace in $H^2(M; \mathbb{Z})$ that annihilates the span of the classes in $H_2(M; \mathbb{Z})$ with zero Thurston norm. Kronheimer and Mrowka's Theorem 41.5.2 and Corollary 40.1.2 in [KM1] imply the following:

> *If* M *is an irreducible 3-manifold, then the unit ball for the dual Thurston norm is the convex hull of the set of elements in* $H^2(M; \mathbb{Z})$ *with the following two properties: First, the element differs from half the first Chern class of* K *by a non-torsion element. Second, the associated Seiberg-Witten Floer homology has some non-torsion elements. Here, {0} is deemed the convex hull of the empty set.*

(1.2)

For another example, suppose that M fibers over the circle with fiber genus greater than one. It follows in this case from [MT] that there exist classes $e \in H^2(M; \mathbb{Z})$ that differ from half K's Euler class by a non-torsion element and have non-trivial Seiberg-Witten Floer homology. This is also the case when M is obtained via zero surgery on a knot in $S^3$ with non-constant Alexander polynomial.

Meanwhile, (1.2) also says something about the vanishing of the Seiberg-Witten Floer homology when tensored with $\mathbb{Q}$. For example, when M is an irreducible 3-manifold, (1.2) and Theorem 1.1 lead to

**Corollary 1.2**: *If* M *is irreducible, suppose that the Euler class of* K *has the following properties: It is not torsion and it has dual Thurston norm greater than 1. Then there is a non-empty set of closed, integral curves of* v *with positive weights that generate the trivial class in* $H_1(M; \mathbb{Z})$.

There are also theorems about the *Seiberg-Witten invariants* that can be applied to the case when M is reducible to deduce the existence of a set of closed integral curves of v that generate the trivial class in $H^2(M; \mathbb{Z})$. In this regard, note that the Seiberg-Witten invariant for a class in $H^2(M; \mathbb{Z})$ is the Euler characteristic of the corresponding Seiberg-Witten Floer homology as defined with respect to a certain $\mathbb{Z}/2\mathbb{Z}$ grading. The results in [MT] with Theorem 1.1 lead to

**Corollary 1.3**: *Suppose that* M *is the connect sum of manifolds* $M_1$ *and* $M_2$. *If the Euler class of* K *is torsion, then (1.1) applies and there is a non-empty set of closed integral curves of v with positive weights that generate the the trivial class in* $H_1(M; \mathbb{Z})$. *If the Euler class of* K *is not torsion, there is a non-empty set of closed, integral curves of* v *with positive weights that generate the trivial class in* $H_1(M; \mathbb{Z})$ *if neither* $M_1$ *nor* $M_2$ *has the* $\mathbb{Z}$*-homology of* $S^3$.

Note that the existence of a closed integral curve in the case considered by Corollary 1.2 can also be deduced using Hofer's celebrated results [H] and a theorem of Eliashberg [E]. Indeed, Eliashberg's result implies the following: If the Euler class of K has dual Thurston norm greater than 1, then the resulting contact 2-plane field is overtwisted. This being the case, Hofer's existence theorem from [H] can be applied. Hofer [H] also established a stronger version of what is said in Corollary 1.3: He proved that a Reeb vector field for a contact structure on a 3-manifold with non-zero $\pi_2$ has a contractible closed orbit.

The theorem of Hofer [H] noted above asserts that the Reeb vector field for a contact form whose kernel is overtwisted has contractible, closed orbits. As is explained in what follows, results from [MR] and [KM2] about Seiberg-Witten Floer homology lead to a proof of the following homological version of Hofer's theorem:

**Theorem 1.4**: *Let* M *denote a compact 3-manifold and let* a *denote a contact form on M whose kernel is overtwisted. Then there is a finite, non-empty set of closed integral curves of the Reeb vector field with an associated set of positive weights such that the weighted sum of curves is homologically trivial.*

Abbas, Cielebak and Hofer [ACH] conjecture that there exists in all cases a union of closed, integral curves of v with positive integer weights that define a homologically trivial cycle.

The proof of Theorem 1.1 uses much the same technology as that used in [T1] to prove what is stated in (1.1). Even so, there are two substantial issues that complicate the application. First, the Seiberg-Witten Floer homology groups are not $\mathbb{Z}$ graded when e differs from half K's Euler class by a non-torsion element. The second substantive issue concerns the Chern-Simons functional: It is not gauge invariant when e differs from half K's Euler class by a non-torsion element. As a consequence, the function whose critical points give the cycles for the Seiberg-Witten Floer homology is not gauge invariant. As is explained below, these two issues constitute the two sides of a single coin, and a single stroke circumvents both.

Theorem 1.4 is proved by coupling the techniques used to prove Theorem 1.1 with constructions and ideas from three sources, [KM1], [KM2], and [MR]. In particular, [MR] provides a vanishing theorem that plays a key role in the proof.

By way of a warning to the reader: This article is truly a sequel to [T1] as it borrows heavily from the latter. In particular, many arguments used here are straightforward modifications of those in [T1]. When such is the case, the reader is referred to the corresponding argument and equations in [T1].

The remainder of this section describes the relevant versions of the Seiberg-Witten equations, and also some of the most quoted results from [T1]. It also supplies a brief definition of the relevant version of the Seiberg-Witten Floer homology. Section 2 closely follows the strategy from [T1] so as to reduce the proof of Theorem 1.1 to a proposition that relates spectral flow to the Chern Simons functional. This key proposition, Proposition 1.9, refines some of the spectral flow results that are stated in Section 5 of [T2]. Section 3 proves this proposition. Section 4 reduces the proof of Theorem 1.4 to a proof of a uniqueness assertion about solutions to a version of the Seiberg-Witten equations on $\mathbb{R} \times M$. The final section proves this uniqueness assertion.

As is the case with [T1], what appears below owes much to previous work of Peter Kronheimer and Tom Mrowka. Moreover, Tom Mrowka suggested using [KM2] and [MR] for Theorem 1.4.

**a) The Seiberg-Witten equations**

Fix a metric on M so that the contact form, a, has norm 1 and is such that da = 2∗a. Here, ∗ denotes the metric's Hodge star operator. With this metric fixed, a $\text{Spin}_{\mathbb{C}}$ structure is no more and no less than an equivalence class of lifts of the oriented, orthonormal frame bundle of M to a principle U(2) bundle.

Such a lift, F, defines an associated, hermitian $\mathbb{C}^2$ bundle, $\mathbb{S} = F \times_{U(2)} \mathbb{C}^2$. The bundle $\mathbb{S}$ is a module for a Clifford multiplication action of T*M. This action is denoted by cl. To be precise, cl is defined so that $cl(\nu_1) cl(\nu_2) = - cl(*(\nu_1 \wedge \nu_2))$ when $\nu_1$ and $\nu_2$ are orthogonal 1-forms. Of course, cl(ν) is anti-hermitian and $cl(\nu)^2 = -|\nu|^2$. A section of $\mathbb{S}$ is called a spinor. Let det($\mathbb{S}$) denote the complex line bundle $F \times_{\text{det}} U(1)$ where det: $U(2) \to U(1)$ is the determinant homomorphism. The Seiberg-Witten equations constitute a system of equations for a pair consisting of a connection on det($\mathbb{S}$) and a section of $\mathbb{S}$.

To write the version of the Seiberg-Witten equations used here, first note that the bundle $\mathbb{S}$ decomposes as a direct summand of eigenbundles for the endomorphism cl(a). This decomposition has the form $\mathbb{S} = E \oplus EK^{-1}$. Here, E is a complex line bundle and K is viewed as a complex line bundle. The convention used for this decomposition of $\mathbb{S}$ takes cl(a) to act as +i on the left-most summand in $\mathbb{S}$ while acting as -i on the right-most summand. A given contact form defines the *canonical Spin$_\mathbb{C}$-structure*; this is the Spin$_\mathbb{C}$ structure whose spinor bundle $\mathbb{S} = \mathbb{S}_I = I_\mathbb{C} \oplus K^{-1}$, where $I_\mathbb{C}$ denotes a trivial line bundle. The assignment $\mathbb{S} \to c_1(E)$ supplies the 1-1 correspondence between $\mathcal{S}_M$ and $H^2(M; \mathbb{Z})$ that is used in (1.1) and Theorem 1.2. With a contact form chosen, the decomposition $\mathbb{S} = E \oplus EK^{-1}$ finds det($\mathbb{S}$) isomorphic to $E^2K^{-1}$.

As is explained momentarily, there is a Hermitian connection on $K^{-1}$ that is canonical up to the action of the gauge group, this the space of maps from M to U(1). With such a 'cannonical connection' chosen, one can view the Seiberg-Witten equations as equations for a pair (A, $\psi$) where A is a connection on E and $\psi$ is a section of $\mathbb{S}$. This is the view taken in what follows.

To be more explicit about the canonical connection, note that a connection on det($\mathbb{S}$) with the Levi-Civita connection on TM determines a Hermitian connection on $\mathbb{S}$. Let $\nabla: C^\infty(M; \mathbb{S}) \to C^\infty(M; \mathbb{S} \otimes T^*M)$ denote the associated covariant derivative. The composition, first $\nabla$ and then the endomorphism from $C^\infty(M; \mathbb{S} \otimes T^*M)$ to $C^\infty(M; \mathbb{S})$ induced by cl defines the Dirac operator, this a first order, hermitian differential operator from $C^\infty(M; \mathbb{S})$ to itself. With the canonical connection chosen, $D_A$ is used in what follows to denote the Dirac operator that is defined by a given connection, A, on E. Granted this notation, suppose that a unit length, trivializing section, $1_\mathbb{C}$, is chosen for the bundle $I_\mathbb{C}$. Let $A_I$ denote the connection on $I_\mathbb{C}$ that makes $1_\mathbb{C}$ covariantly constant. There is a unique connection on K with the property that $1_\mathbb{C}$ is annihilated by the $A = A_I$ version of $D_A$. The latter connection is the promised canonical connection.

The specification of the Seiberg-Witten equations requires the choice of a number, r, with $r \geq 1$. Also required is the choice of a smooth, 1-form on M with small $C^3$ norm; in particular the norm should be less than 1. This form is denoted in what follows by $\mu$. The form $\mu$ should be $L^2$-orthogonal to all harmonic 1-forms on M and it should be coclosed, thus $d*\mu = 0$. Granted these choices, the simplest of the relevant versions of the Seiberg-Witten equations read:

- $B_A = r(\psi^\dagger \tau \psi - ia) + i(*d\mu + \varpi_K)$,
- $D_A \psi = 0$.

(1.3)

Here, $B_A$ denotes the Hodge dual of A's curvature 2-form. Meanwhile, $\psi^\dagger \tau \psi$ is a short-hand way to write the image in $C^\infty(M; iT^*M)$ of the endomorphism $\psi \otimes \psi^\dagger$ via the adjoint of cl. Finally, $\varpi_K$ is the harmonic 1-form whose Hodge dual is -$\pi$ times the first Chern class of the bundle K.

In what follows, the form $\mu$ is chosen from a certain Banach space of smooth forms on M. The precise norm for this space is not relevant. Suffice it to say here that the space is dense in $C^\infty(M; T^*M)$. In particular, the Banach space norm bounds a suitable, positive multiple of any given $C^k$ norm. The Banach space used here is denoted by $\Omega$.

Note that the equations in (1.3) are gauge invariant. This means the following: Suppose that $(A, \psi)$ satifies (1.3) and that $u$ is a smooth map from M to $S^1 = U(1)$. Then the pair $(A - u^{-1}du, u\psi)$ also satisfies (1.3). The set of gauge equivalence classes of solutions to (1.3) is compact for any choice of $\mu$.

Note also that there are no $\psi = 0$ solutions to (1.3) in the cases under consideration, those where $c_1(E) - \frac{1}{2}c_1(K)$ is not a torsion class.

Introduce Conn(E) to denote the space of smooth, Hermitian connections on E. The equations in (1.3) are the critical points of a functional on $\text{Conn}(E) \times C^\infty(M; \mathbb{S})$. It proves convenient to write this functional, $\mathfrak{a}$, as

$$\mathfrak{a} = \tfrac{1}{2}(\mathfrak{cs} + 2\mathfrak{e}_\mu - r\mathrm{E}) + r\int_M \psi^\dagger D_A \psi \ .$$

(1.4)

where $\mathfrak{cs}$, $\mathfrak{e}_\mu$ and $\mathrm{E}$ are functionals on Conn(E) that are defined as follows: The definition of $\mathfrak{cs}$ requires the choice of a connection, $A_E$, on E. Here, $A_E$ is chosen so that its curvature 2-form is harmonic. Write $A = A_E + \hat{a}_A$, where $\hat{a}_A$ is a section of $iT^*M$. In the case when $E = 1_\mathbb{C}$, choose $A_E$ to be the trivializing connection $A_I$ that was introduced above. Use $\hat{c}_{1\mathbb{S}}$ to denote the harmonic 2-form on M that represents the first Chern class of the line bundle $\det(\mathbb{S})$. Then

- $\mathfrak{cs}(A) = -\int_M \hat{a}_A \wedge *d\hat{a}_A - 2i\pi \int_M \hat{a} \wedge \hat{c}_{1\mathbb{S}}$ ,
- $\mathfrak{e}_\mu = i \int_M \mu \wedge *B_A$ ,
- $\mathrm{E}(A) = i \int_M a \wedge *B_A$ .

(1.5)

The introduction referred to the Chern Simons functonal. This is $\mathfrak{cs}$. All but $\mathfrak{cs}$ in (1.4) are fully gauge invariant. The functional $\mathfrak{cs}$ is gauge invariant if and only if $c_1(\det(\mathbb{S})) = 2c_1(E) - c_1(K)$ is a torsion class. In the cases under consideration, this class is <u>not</u> torsion. In any event, $\mathfrak{cs}(A - u^{-1}du) = \mathfrak{cs}(A)$ if $u$ is a homotopically trivial map from M to $S^1$. In general, $\frac{1}{4\pi^2}\mathfrak{cs}$ differs by a multiple of the greatest integer divisor of $c_1(\det(\mathbb{S}))$ on gauge equivalent pairs.

**b) Properties of solutions to the Seiberg-Witten equations**

This subsection restates for ease of reference some salient facts from [T1] about solutions to (1.3). The first of these explains how solutions to (1.3) lead to closed integral curves of v.

**Theorem 1.5:** *Fix a complex line bundle* E *so as to define a Spin$_\mathbb{C}$-structure on* M *with spinor bundle* $\mathbb{S}$ *given by* $E \oplus EK^{-1}$. *Let* $\{r_n\}_{n=1,2,...}$ *denote an increasing, unbounded sequence of postive real numbers, and for each* n, *let* $\mu_n$ *denote a co-exact 1-form on M. For each* n, *let* $(A_n, \psi_n) \in \text{Conn}(E) \times C^\infty(M; \mathbb{S})$ *denote a solution to the* $r = r_n$ *version of (1.3) as defined using* $\mu_n$. *Suppose that there is an* n-*independent bound for the* $C^3$ *norm of* $\mu_n$. *In addition, suppose that there exists an index* n *independent upper bound for* $\mathcal{E}(A_n)$ *and a strictly positive,* n *independent lower bound for* $\sup_M (1 - |\psi_n|)$. *Then there exists a non-empty set of closed, integral curves of the Reeb vector field. Moreover, there exists a positive, integer weight assigned to each of these integral curves such that the corresponding formal, integer weighted sum of loops in* M *gives the class in* $H_1(M; \mathbb{Z})$ *that is Poincaré dual to the first Chern class of the bundle* E.

This is Theorem 2.1 in [T1] and its proof is in Section 6 of [T1].
       To say more about what is happening in Theorem 1.5, agree to write a spinor $\psi$ with respect to the decomposition $\mathbb{S} = E \oplus K^{-1}E$ as $\psi = (\alpha, \beta)$. If the sequences $\{r_n\}$, $\{\mu_n\}$ and $\{(A_n, \psi_n = (\alpha_n, \beta_n))\}$ satisfy the conditions stated by Theorem 1.5, then the corresponding sequence $\{\alpha_n^{-1}(0)\}_{n=1,2,...}$ of closed subsets of M has a subsequence that converges in a suitable sense to the desired set of closed integral curves.
       The next three lemmas summarizes various facts about the solutions to (1.3). They restate what is asserted in Lemmas 2.2, 2.3 and 2.4 in [T1]. The proofs of the latter are also in Section 6 of [T1]

**Lemma 1.6**: *Fix the line bundle* $E \to M$ *and hence the Spin$_\mathbb{C}$ structure. There is a constant* $\kappa \geq 1$ *with the following significance: Let* $\mu$ *denote a co-exact 1-form with* $C^3$ *norm bounded by 1. Fix* $r \geq 1$ *and suppose that* $(A, \psi = (\alpha, \beta))$ *is a solution to the version of (1.3) given by* r *and* $\mu$. *Then*
- $|\alpha| \leq 1 + \kappa\, r^{-1}$ ,
- $|\beta|^2 \leq \kappa\, r^{-1}(|1 - |\alpha|^2| + r^{-1})$ .

The next lemma concerns the derivatives of $\alpha$ and $\beta$. To state the lemma, suppose that A is a given connection on E. Introduce $\nabla$ to denote the associated covariant derivative. The covariant derivative on $K^{-1}E$ that is defined by A and the canonical connection is denoted in what follows by $\nabla'$. Powers of these covariant derivatives are defined using A and the Levi-Civita induced connection on the appropriate tensor power of T*M.

**Lemma 1.7:** *Fix the line bundle* $E \to M$ *and hence the Spin$_\mathbb{C}$ structure. For each integer* $q \geq 0$ *and constant* $c \geq 0$, *there is a constant* $\kappa \geq 1$ *with the following significance: Let* $\mu$ *denote a co-exact 1-form whose $C^{3+q}$ norm is bounded by* c. *Fix* $r \geq 1$ *and suppose that* $(A, \psi = (\alpha, \beta))$ *is a solution to the version of (1.3) defined by* r *and* $\mu$. *Then*

- $|\nabla^q \alpha| \leq \kappa (r^{q/2} + 1)$.
- $|\nabla'^q \beta| \leq \kappa (r^{(q-1)/2} + 1)$.

The final lemma says something about the connection A.

**Lemma 1.8**: *Fix the line bundle* E *and hence the Spin$_\mathbb{C}$ structure. Also, fix a connection* $A_E$ *on* E. *There exists a constant* $\kappa \geq 1$ *with the following significance: Let* $\mu$ *denote a co-exact 1-form with $C^3$ norm bounded by 1. Fix* $r \geq 1$ *and suppose that* $(A, \psi = (\alpha, \beta))$ *is a solution to the version of (1.3) defined by* r *and* $\mu$. *Then, there exists a smooth map* u: $M \to S^1$ *such that* $\hat{a} = A - u^{-1}du - A_E$ *obeys* $|\hat{a}| \leq \kappa (r^{2/3} E^{1/3} + 1)$.

Note that Lemma 1.6 implies that

$$-c_\mu \leq E(A) \leq r \, vol(M) + c_\mu ,$$

(1.6)

when $(A, \psi)$ satisfies (1.3). Here, $c_\mu$ depends only on the $C^3$ norm of $\mu$.

**c) Spectral Flow**

Each pair of $r \in [1, \infty)$ and element $(A, \psi) \in \text{Conn}(E) \times C^\infty(M; \mathbb{S})$ has an associated, first order, elliptic differential operator. This operator is denoted by $\mathfrak{L}$; it maps $C^\infty(M; iT^*M \oplus \mathbb{S} \oplus i\mathbb{R})$ to itself, and it is defined as follows: Let $\mathfrak{h} = (b, \eta, \phi) \in C^\infty(M; iT^*M \oplus \mathbb{S} \oplus i\mathbb{R})$. The respective $iT^*M$, $\mathbb{S}$ and $i\mathbb{R}$ summands of the $\mathfrak{L}\mathfrak{h}$ are:

- $*db - d\phi - 2^{-1/2} r^{1/2} (\psi^\dagger \tau \eta + \eta^\dagger \tau \psi)$,
- $D_A \eta + 2^{1/2} r^{1/2} (cl(b)\psi + \phi\psi)$,
- $*d*b - 2^{-1/2} r^{1/2} (\eta^\dagger \psi - \psi^\dagger \eta)$.

(1.7)

The domain of this operator extends to map $L^2(M; iT^*M \oplus \mathbb{S} \oplus i\mathbb{R})$ to itself as a self-adjoint operator with domain $L^2_1(M; iT^*M \oplus \mathbb{S} \oplus i\mathbb{R})$. Its spectrum is discrete, all eigenvalues have finite multiplicity, there are no accumulation points for the spectrum, and the spectrum is unbounded from above and from below. A pair $(A, \psi)$ is called a *non-degenerate* solution to (1.3) when the r and $(A, \psi)$ version of $\mathfrak{L}$ has trivial kernel.

Let $A_E$ denote the connection that was chosen above to define $\mathfrak{cs}$. Fix $r_E \geq 1$. If $\psi$ is a suitably generic section of $C^\infty(M; \mathbb{S})$ with $|\psi| \leq 1$, then the $r = r_E$ and $(A = A_E, \psi = \psi_E)$ version of $\mathfrak{L}$ will have trivial kernel. Agree to fix $\psi_E$ once and for all. In the case when E $= 1_\mathbb{C}$, it follows from Lemma 5.4 in [T1] that the pair $(A_I, \psi_I)$ is non-degenerate when r is greater than some constant $r_*$. This understood, take $r_E$ in to be $r_* + 1$ and $\psi_E = \psi_I$ when E is the trivial bundle.

Now, suppose that $r \geq 1$ and $(A, \psi) \in \text{Conn}(E) \times C^\infty(M; \mathbb{S})$ are such that the r and $(A, \psi)$ version of $\mathfrak{L}$ also has trivial kernel. Then one can associate a number, $f(A, \psi)$, to the pair $(A, \psi)$ as follows: Fix a path $t \to (r(t), \mathfrak{c}(t)) \in [1, r) \times \text{Conn}(E) \times C^\infty(M; \mathbb{S})$ that is parameterized by $t \in [0, 1]$, starts at $(1, (A_E, \psi_E))$ and ends at $(r, (A, \psi))$. Take $f(A, \psi)$ equal to the spectral flow for the corresponding path, $t \to \mathfrak{L}_t$ of versions of (1.7). For those unfamiliar with this notion, see Section 5a in [T1] or [T2]. The value of $f$ does not depend on the chosen path. Note in what follows that $f$ defines a locally constant function on a dense, open subset of $\text{Conn}(E) \times C^\infty(M; \mathbb{S})$. (It can be shown that $f$ is defined on the complement of a codimension 1 subvariety.)

When r is given and $(A, \psi)$ is a non-degenerate solution to (1.3), the value of $-f$ is called the *degree* of $(A, \psi)$. This notion of degree assigns an integer to each non-degenerate solution to (1.3). Note that the degree as just defined depends on the choice of $A_E$ and also on $\psi_E$, but the difference between the degrees of any two elements from $\text{Conn}(E) \times C^\infty(M; \mathbb{S})$ does not depend on these choices.

The function $f$ is not gauge invariant when $c_1(E) - \frac{1}{2} c_1(K)$ is not torsion; its values can differ on gauge equivalent configurations by a multiple of the greatest common divisor of $c_1(\det(\mathbb{S}))$. However,

$$\mathfrak{cs}^f = \mathfrak{cs} - 4\pi^2 f$$

(1.8)

is gauge invariant. This understood, the following proposition plays a central role in the proof of Theorems 1.1 and 1.5.

**Proposition 1.9**: *There exists a constant $\kappa \geq 1$ with the following significance: Fix $r \geq \kappa$ and $\mu \in \Omega$ with $C^3$ norm less than 1, and let $(A, \psi)$ denote a non-degenerate solution to the r and $\mu$ version of (1.3). Then $|\mathfrak{cs}^f(A, \psi)| \leq \kappa r^{2/3} (1 + |E|)^{4/3} (\ln r)^\kappa$.*

This proposition is new; its proof is given in Section 3.

Under certain circumstance, Proposition 1.9 is a refinement of a version of Proposition 5.1 from [T1] that holds in the case when $c_1(\det(\mathbb{S}))$ is not torsion. The latter also plays a role in what follows.

**Proposition 1.10**: *There exists $\kappa$ with the following significance: Fix $r \geq \kappa$ and $\mu \in \Omega$ with $C^3$ norm less than 1, and let $(A, \psi)$ denote a non-degenerate solution to the r and $\mu$ version of (1.3). Then $|\mathfrak{cs}^f(A, \psi)| \leq \kappa r^{31/16}$.*

*Proof of Proposition 1.10*:  Repeat verbatim the arguments given in Sections 5c and 5d of [T1] and then appeal to Proposition 5.5 in [T1].

**d) Seiberg-Witten Floer homology**

The definition used here of Seiberg-Witten Floer homology in the case when $c_1(E) - \frac{1}{2} c_1(K)$ is not torsion comes from [KM1]. What follows is a summary. To start, the definition requires Seiberg-Witten equations that are perturbed versions of those that are depicted in (1.3). Each such perturbation of (1.3) is the variational equations of a function on $\text{Conn}(E) \times C^\infty(M; \mathbb{S})$ that is obtained from what is depicted in (1.4) by adding one of an admissable class of gauge invariant functionals to what appears on (1.4)'s right hand side. The precise nature of this class of admissable additions is described in Chapter 11.6 of [KM1], see Definition 11.6.3. Suffice it to say here that these admissable functionals come from a separable Banach space of functions on $\text{Conn}(E) \times C^\infty(M; \mathbb{S})$, and that one need only consider additions with very small norm. This Banach space is denoted in what follows by $\mathcal{P}$, and its norm is denoted by $\|\cdot\|_\mathcal{P}$. Note that the function $\mathfrak{e}_\mu$ is in $\mathcal{P}$ when $\mu$ is in the Banach space $\Omega$. A function from $\mathcal{P}$ is called a *perturbation function*, and the terms that it supplies to (1.3) are called *perturbation terms*. Fix $r \geq 1$ and a perturbation function $\mathfrak{g}$ on $\text{Conn}(E) \times C^\infty(M; \mathbb{S})$. The resulting version of the Seiberg-Witten equations reads

- $B_A = r(\psi^\dagger \tau \psi - ia) + \mathfrak{T}(A, \psi) + i \varpi_K$,
- $D_A \psi = \mathfrak{S}(A, \psi)$.

(1.9)

where $\mathfrak{T}$ and $\mathfrak{S}$ are defined by requiring that their respective $L^2$ inner products with sections b of $iT^*M$ and $\eta$ of $\mathbb{S}$ are $\mathfrak{T} = -\frac{d}{dt} \mathfrak{g}(A + tb, \psi)|_{t=0}$ and $\mathfrak{S} = -\frac{d}{dt} \mathfrak{g}(A, \psi + t\eta)|_{t=0}$. For example, the equations that appear in (1.3) use $\mathfrak{g} = \mathfrak{e}_\mu$. A solution to (1.9) is deemed to be non-degenerate when a certain perturbed version of (1.7) has trivial kernel. This perturbed version of (1.7) is depicted in (3.1) of [T1]. The perturbation subtracts respective terms $\mathfrak{t}_{(A,\psi)}(b, \eta) = \frac{d}{dt} \mathfrak{T}(A+tb, \psi+t\eta)|_{t=0}$ and $\mathfrak{s}_{(A,\psi)} = \frac{d}{dt} \mathfrak{S}(A+tb, \psi+t\eta))|_{t=0}$ from what is written in the first and second lines of (1.7).

As in the cases that are considered in [T1], there is a dense, open set of forms $\mu \in \Omega$ such that with r fixed, then (1.3) has a finite set of gauge equivalence classes of solutions, with each non-degenerate in the sense that the operator in (1.7) has trivial kernel. With it understood that $c_1(E) - \frac{1}{2} c_1(K)$ is not torsion, the arguments from [KM1] that justified Lemma 3.1 and 3.2 in [T1] also prove

**Lemma 1.11**: *Given $r \geq 0$ there is a residual set of $\mu \in \Omega$ such that all solutions to the corresponding version of (1.3) are non-degenerate. There is an open dense set of $\mathfrak{g} \in \mathcal{P}$* *such that all solutions to the corresponding version of* (1.9) *are non-degenerate. For a given* r ≥ 1 *and* μ ∈ Ω *and* q ∈ 𝒫, *suppose that* c = (A, ψ) *is a non-degenerate solution to the* r *and* g = e_μ + q *version of (1.9). Then the following is true:*

- *There exist* ε > 0 *such that if* (A´, ψ´) *is a solution to (1.9) that is not gauge equivalent to* (A, ψ), *then the* $L^2_1$ *norm of* (A - A´, ψ - ψ´) *is greater than* ε.
- *There is a smooth map,* c(·), *from the ball of radius* ε *centered at the origin in* 𝒫 *to* Conn(E) × C^∞(M; 𝕊) *such that* c(0) = c *and such that* c(p) *is a non-degenerate solution to the version of (1.9) defined by* r *and the perturbation* g = e_μ + q + p.

With r chosen, fix μ ∈ Ω so that all solutions to the resulting version of (1.3) are non-degenerate. The cycles for the Seiberg-Witten Floer homology consist of the elements in the vector space over ℤ that is generated by the set of gauge equivalence classes of solutions to (1.3). Let p denote the greatest integer divisor of $c_1$(det(𝕊)). Any given generator for the Seiberg-Witten Floer homology has a ℤ/pℤ grading; and for the purpose of this article, the ℤ/pℤ grading for a given generator is defined as follows: The grading of the generator (A, ψ) is the image of -ƒ(A, ψ) in ℤ/pℤ. The grading is defined only modulo p because values of ƒ on gauge equivalent pairs in Conn(E) × C^∞(M; 𝕊) can differ by multiples of p.

The differential on the chain complex is obtained using instanton solutions to the Seiberg-Witten equations on ℝ × M. Given r ≥ 1 and a perturbation function g ∈ 𝒫, an instanton is a smooth map s → ∂(s) from ℝ to Conn(E) × C^∞(M; 𝕊) that satisfies the equation

- $\frac{\partial}{\partial s}$ A = -B_A + r(ψ^† τ^k ψ - ia) + 𝔗(A, ψ) + i϶_K,
- $\frac{\partial}{\partial s}$ ψ = -D_A ψ + 𝔖(A, ψ).

(1.10)

and is such that both $\lim_{s \to -\infty}$ ∂(s) and $\lim_{s \to \infty}$ ∂(s) are solutions to the r and g version of (1.9). The instantons determine the differential as follows: Let c and c´ denote a pair of solutions to (1.9) for a given r and perturbation g. Use ℳ(c, c´) to denote the space of instanton solutions s → ∂(s) = (A(s), ψ(s)) to (1.10) whose s → -∞ limit is equal to c and whose s → +∞ limit equal to û·c´ for some û ∈ C^∞(M; S^1). Given an integer ɩ, introduce ℳ_ɩ(c, c´) to denote the union of those components of ℳ(c, c´) with the following property: Define the family of operators $\{\mathfrak{L}_s\}_{s \in \mathbb{R}}$ with $\mathfrak{L}_s$ defined by the given value of r and the configuration ∂(s) = (A(s), ψ(s)) in (3.1). Then ∂ ∈ ℳ_ɩ(c, c´) if and only if the spectral flow for $\{\mathfrak{L}_s\}_{s \in \mathbb{R}}$ is equal to -ɩ. Theorems 15.1.1 and 16.1.3 in [KM1] provide the following analog to Lemma 3.6 in [T1]:

**Lemma 1.12**: *Given $r \geq 1$, there exists an open, dense set of $\mu \in \Omega$ with $C^3$ norm less than 1 for which the following is true*:
- *Each solution to the $r$ and $\mu$ version of* (1.3) *is non-degenerate.*
- *Given $\mu$ for which the preceding conclusions hold, there exists $\varepsilon > 0$ and a dense, open subset of the radius $\varepsilon$ ball in $\mathcal{P}$ such that if $\mathfrak{q}$ is in this set, then the $\mathfrak{g} = \mathfrak{e}_\mu + \mathfrak{q}$ versions of $\{\mathcal{M}_{\iota=1}(\cdot, \cdot)\}$ can be used to compute the differential on the Seiberg-Witten Floer complex. In particular, for such $\mathfrak{q}$,*
   a) $\mathcal{M}_\iota(\cdot, \cdot) = \emptyset$ *when $\iota < 0$.*
   b) $\mathcal{M}_{\iota=0}(\mathfrak{c}, \mathfrak{c}') = \emptyset$ *unless $\mathfrak{c} = \mathfrak{c}'$ in which case it contains only the constant map $s \to \mathfrak{d}(s) = \mathfrak{c}$.*
   c) $\mathcal{M}_{\iota=1}(\mathfrak{c}, \mathfrak{c}')$ *is a finite set of copies of $\mathbb{R}$.*

A pair $(\mu, \mathfrak{q}) \in \Omega \times \mathcal{P}$ as described by Lemma 1.12 is called r-admissable. When $(\mu, \mathfrak{q})$ is r-admissable, then the $(\mathfrak{c}, \mathfrak{c}')$ version of $\mathcal{M}_{\iota=1}(\mathfrak{c}, \mathfrak{c}')$ is 1-dimensional and each component has a well defined sign, this described in Chapter 22.1 of [KM1]. The sum of these signs defines an integer that is denoted here by $\sigma(\mathfrak{c}, \mathfrak{c}')$. Take $\sigma(\mathfrak{c}, \mathfrak{c}')$ to equal 0 when $\mathcal{M}_{i=1}(\mathfrak{c}, \mathfrak{c}') = \emptyset$. An integer so defined is associated to any generators, $\mathfrak{c}$ and $\mathfrak{c}'$, of the Seiberg-Witten Floer context with the property that degree($\mathfrak{c}$) = degree($\mathfrak{c}'$) + 1. This understood, the differential on the Seiberg-Witten Floer complex is given by the rule:

$$\partial \mathfrak{c} = \sum_{\mathfrak{c}'} \sigma(\mathfrak{c}, \mathfrak{c}') \, \mathfrak{c}'.$$

(1.11)

The differential increases the mod(p) degree of any given generator by 1.

An appeal to Proposition 22.1.4 and Corollary 23.1.6 in [KM1] gives the following analog to Proposition 3.7 in [T1]:

**Proposition 1.12**: *Fix $r \geq 1$ and an r-admissable pair $(\mu, \mathfrak{q})$ to define the generators and differential on the Seiberg-Witten Floer comple. Then $\delta^2 = 0$. Moreover, given two r admissable pairs, there exists an isomorphism between the corresponding versions of Seiberg-Witten Floer homology. In addition, the homology so defined for different values of $r \geq 1$ are isomorphic.*

This last proposition ends the introduction to Seiberg-Witten Floer homology.

**e) Conventions**

The reader is forwarned here about a notational convention that is used in the subsequent sections of this paper: Whenever $c_0$ appears, it represents a postive constant whose value depends only on a chosen metric, the 1-forms $a$ and $\mu$, the bundle $E$, and perhaps other fixed objects. In particular, the value of $c_0$ is independent of the value of $r$ and of any particular solution to (1.3), (1.7), (1.9) or (1.10) that is under consideration. The precise value for $c_0$ will typically change from appearance to appearance.

## 2. The proof of Theorem 1.1

The plan is to fix a suitable, small normed element $\mu \in \Omega$, and then find an increasing sequence $\{r_n\}_{n=1,2,...}$ and, for each n, a pair $(A_n, \psi_n) \in \text{Conn}(E) \times C^\infty(M; \mathbb{S})$ that solves the $r = r_n$ and $\mu$ version of (1.3); and then appeal to Theorem 1.5 to find the desired set of closed integral curves of v. For this purpose, it is necessary to find $(A_n, \psi_n)$ so that the corresponding sequence $\{E(A_n)\}_{n=1,2,...}$ has a finite upper bound, and so that the sequence $\{\sup_M (1 - |\psi_n|)\}_{n=1,2,...}$ has a positive lower bound.

As will be clear from what follows, the purpose of Proposition 1.12 is to supply the solutions to (1.3). For example, when the Seiberg-Witten invariant is non-zero, then any r and $\mu$ version of (1.3) has at least one gauge equivalence class of solution. As in [T1], the subtle issue is that of finding a sequence $\{r_n\}_{n=1,2,...}$, and the corresponding sequence $\{(A_n, \psi_n)\}_{n=1,2,...}$ such that $\{E(A_n)\}_{n=1,2,...}$ is bounded. The arguments that are used in [T1] are modified in what follows to find such a sequence. The existence of a positive lower bound for $\{(1 - |\psi_n|)\}_{n=1,2,...}$ for the case when E is non-trivial follows immediately from Lemma 1.6 by virtue of the fact that the first component of any section of $\mathbb{S} = E \oplus EK^{-1}$ must vanish at some point in M in this case. Of course, this argument can not be made in the case when E is trivial. The need for an alternate argument when $E = I_\mathbb{C}$ is the proximate cause for the separate assertion in Theorem 1.1 for $e = 0$ case.

What follows outlines how the arguments from [T1] procede in the situation at hand. Proposition 1.12 asserts that the Floer homologies that are defined for different values of r and different perturbation functions are isomorphic. This understood, it is first necessary to construct isomorphisms that lift in a suitable fashion to give homomorphisms between the relevant chain complexes. This is done so as to facilitate the assignment of a certain continuous function of r to any given Seiberg-Witten Floer homology class. In the case considered in [T1], this function was the value of $\mathfrak{a}$ on a certain generator in a certain representative cycle for the class. The cycle and generator are found using a min-max construction. In the present situation, $\mathfrak{a}$ is not gauge invariant and so lacks a definite value on any given generator of the Seiberg-Witten Floer complex. However,

$$\mathfrak{a}^f = \mathfrak{a} - 2\pi^2 f.$$

(2.1)

is defined on any non-degenerate solution to (1.3) and is gauge invariant. As is explained in what follows, the min-max construction from [T1] work just as well with $\mathfrak{a}^f$ in lieu of $\mathfrak{a}$ to assign a continuous function of r to each Seiberg-Witten Floer homology class. As it turns out, the r dependence of $\mathfrak{a}^f$ is such that the arguments from Section 4 of [T1] can be employed with but one significant change to find the sequence $\{r_n, (A_n, \psi_n)\}$ with $\{E(A_n)\}$ bounded. This change requires new input, and this is supplied by Proposition 1.9.

### a) Assigning certain functions of r to Seiberg-Witten Floer homology classes

The purpose of this subsection is to associate values of $\mathfrak{a}^f$, E and $\mathfrak{cs}^f$ to a given Seiberg-Witten Floer homology class. This association will provide functions of r that are analogous to the functions that are described in Definitions 4.1 and 4.4 of [T1].

The first step for this task is to define chain homotopies between the respective Seiberg-Witten Floer complexes that are defined at different values of r. The construction starts with an analog of Proposition 3.11 from [T1].

**Proposition 2.1**: *There is a residual subset in $\Omega$ with $C^3$ norm less than 1 and with the following properties: Let $\mu$ denote a form from this subset. There is a locally finite set $\{\rho_j\} \subset (1, \infty)$ with $\rho_1 < \rho_2 < \cdots$ such that if $r > 1$ and $r \notin \{\rho_j\}$, then*
1) *Each solution to the r and $\mu$ version of (1.3) is non-degenerate.*
2) *Define $\mathfrak{a}$ using the r and $\mathfrak{g} = \mathfrak{e}_\mu$ version of (1.4). If $\mathfrak{c}$ and $\mathfrak{c}'$ are solutions the r and $\mu$ version of (1.3) that are not gauge equivalent, then $\mathfrak{a}^f(\mathfrak{c}) \neq \mathfrak{a}^f(\mathfrak{c}')$.*

*Proof of Proposition 2.1*: The proof is identical save for some minor notational changes to the proof of Proposition 3.11 given in Section 7b of [T1].

This last proposition can be used to define a 'canonical' basis and for the Seiberg-Witten Floer complex when $r \notin \{\rho_j\}$. This is done with the help of the following lemma:

**Lemma 2.2**: *If all solutions to the r and $\mu$ version of (1.3) are non-degenerate, then there exists $\varepsilon > 0$ with the following significance:*
- *If $\mathfrak{q} \in \mathcal{P}$ has norm less than $\varepsilon$, then all solutions to the $\mathfrak{g} = \mathfrak{e}_\mu + \mathfrak{q}$ version of (1.4) are non-degenerate.*
- *If $\mathfrak{q} \in \mathcal{P}$ has norm less than $\varepsilon$, then there is a 1-1 correspondence between the set of solutions to the r and $\mu$ version of (1.3) and the set of solutions to the r and $\mathfrak{g} = \mathfrak{e}_\mu + \mathfrak{q}$ version of (1.9).*
- *In particular, if $\mathfrak{c}$ is a solution to the r and $\mu$ version of (1.3), then there exists a smooth map, $\mathfrak{c}(\cdot)$, from the radius $\varepsilon$ ball in $\mathcal{P}$ into $\mathrm{Conn}(E) \times C^\infty(M; \mathbb{S})$ such that $\mathfrak{c}(\mathfrak{q})$ solves the r and $\mathfrak{e}_\mu + \mathfrak{q}$ version of (1.9) and such that $\mathfrak{c}(0) = \mathfrak{c}$.*
- *If $\mathfrak{c}$ is a solution to the r and $\mu$ version of (1.3), then the assignment $\mathfrak{q} \to \mathfrak{a}^f(\mathfrak{c}(\mathfrak{q}))$ defines a smooth function on the radius $\varepsilon$ ball in $\mathcal{P}$.*
- *If $\mathfrak{a}^f(\mathfrak{c}) \neq \mathfrak{a}^f(\mathfrak{c}')$ when $\mathfrak{c}$ and $\mathfrak{c}'$ are distinct gauge equivalence classes of solutions to the r and $\mu$ version of (1.3), then $\varepsilon$ can be chosen so that the ranges of the functions $\mathfrak{a}^f(\mathfrak{c}(\cdot))$ and $\mathfrak{a}^f(\mathfrak{c}'(\cdot))$ are disjoint.*

*Proof of Lemma 2.2:* This follows from Proposition 2.1 with Lemma 1.11 given that $f$ is locally constant on its domain of definition in $\text{Conn}(E) \times C^\infty(M; \mathbb{S})$.

Suppose now that all solutions to the $r$ and $\mu$ version of (1.3) are non-degenerate and are such that $\mathfrak{a}^f$ distinguishes distinct gauge equivalence classes of solutions. For each degree $k \in \mathbb{Z}/p\mathbb{Z}$, label the degree $k$ gauge equivalence classes of solutions to (1.3) by the consecutive integers starting at 1 using the convention that $\mathfrak{a}^f(\mathfrak{c}_\nu) > \mathfrak{a}^f(\mathfrak{c}_{\nu+1})$. This labeling is deemed the *canonical labeling*. Now suppose that $\varepsilon > 0$ is given by Lemma 2.2, that $\mathfrak{q} \in \mathcal{P}$ has norm less than $\varepsilon$ and that $(\mu, \mathfrak{q})$ is $r$-admissable. If $\mathfrak{q}$ has these properties, then $(\mu, \mathfrak{q})$ is called *strongly r-admissable*. Assume that $\mathfrak{q}$ has these properties. For each degree $k$ and each labeling integer $\nu$, let $\mathfrak{c}_\nu(\cdot)$ denote the map that is provided by Lemma 2.2 with $\mathfrak{c}_\nu(0) = \mathfrak{c}_\nu$. Then the collection $\{\mathfrak{c}_\nu(1)\}$ defines a labeling of the generators for the degree $k$ generators of the Seiberg-Witten Floer complex as defined using the $r$ and $\mathfrak{g} = \mathfrak{e}_\mu + \mathfrak{q}$. The corresponding basis for this complex is deemed the *canonical basis*.

Granted this terminology, what follows is an analog of Proposition 3.12 in [T1].

**Proposition 2.3**: *Let $\mu$ be as in Proposition 2.1. Fix $\rho_j \in (r_k, \infty)$ from the set described in Proposition 2.1. There exists a possibly empty, but contiguous set $\mathbb{J}(i) \subset \mathbb{Z}$, and a corresponding sequence $\{t_m\}_{m \in \mathbb{J}(i)} \in (\rho_i, \rho_{i+1})$ with the following properties:*

- *The sequence is increasing, and it has no accumulation points in the open interval.*
- *For any given $m \in \mathbb{J}(i)$, there exists $\mathfrak{q}_m \in \mathcal{P}$ of small norm such that $(\mu, \mathfrak{q}_m)$ is strongly r-admissable for all $r \in [t_m, t_{m+1}]$.*
- *When the canonical basis is used for the $r \in [t_m, t_{m+1}]$ and $\mathfrak{g} = \mathfrak{e}_\mu + \mathfrak{q}_m$ versions of the Seiberg-Witten Floer complex, then the differential on this complex is independent of $r$ as $r$ varies in $[t_m, t_{m+1}]$. This differential is denoted by $\delta_m$.*
- *Let $m \in \mathbb{J}(i)$. In each $\mathbb{Z}/p\mathbb{Z}$ degree, there is an upper triangular, integer valued matrix, $\mathbb{A}$, with 1 on the diagonal such that $\delta_m = \mathbb{A}^{-1}\delta_{m-1}\mathbb{A}$. Here, both $\delta_m$ and $\delta_{m-1}$ are written with respect to the canonical basis.*

*Proof of Proposition 2.3*: But for two observations and some terminology, the proof is identical to the proof of Proposition 3.12 in [T1]. Here is the first observation: Let $\mathbb{V}$ denote the $\mathbb{Z}/p\mathbb{Z}$ graded vector space with basis in any given degree $k$ given by the canonical basis for the Seiberg-Witten Floer complex as defined using the gauge equivalence classes of solutions to the $r = t_{m-1}$ and $\mu$ version of (1.5). Fix a path $s \to \mathfrak{q}(s)$ as in Part 1 of the proof in Section 7c in [T1] of Proposition 3.12 in [T1]. Use this path to define the equations that appear in (7.12) of [T1]. Let $\mathfrak{c}_-$ denote a solution to the $r = t_{m-1}$ and $\mathfrak{g} = \mathfrak{e}_\mu + \mathfrak{q}_{m-1}$ version of (1.9) and let $\mathfrak{c}_+$ denote a solution to the $r = t_m$ and $\mathfrak{g} = \mathfrak{e}_\mu + \mathfrak{q}_m$ version. Assume that $\mathfrak{c}_-$ and $\mathfrak{c}_+$ have the same degree in $\mathbb{Z}/p\mathbb{Z}$. Let $\mathcal{M}_{q(\cdot)}(\mathfrak{c}_-, \mathfrak{c}_+)$ denote the space of solutions $s \to \mathfrak{d}(s)$ to (1.9) with $s \to -\infty$ limit equal to $\mathfrak{c}_-$ and with $s \to \infty$ limit equal to $\hat{u} \cdot \mathfrak{c}_+$ with $\hat{u}$ a smooth map from M to $S^1$. Let $s \to (A(s), \psi(s))$ denote an element in $\mathcal{M}_{q(\cdot)}(\mathfrak{c}_-, \mathfrak{c}_+)$. For each $s \in \mathbb{R}$, let $\mathfrak{L}_s$ denote the version of the operator in (3.1) of [T1] that uses the pair $(A(s), \psi(s))$ and the perturbation data coming from $\mathfrak{g} = \mathfrak{e}_\mu + \mathfrak{q}(s)$. Use $\mathcal{M}_{q(\cdot)0}(\mathfrak{c}_-, \mathfrak{c}_+) \subset \mathcal{M}_{q(\cdot)}(\mathfrak{c}_-, \mathfrak{c}_+)$ to denote the subset of $\mathcal{M}_{q(\cdot)}(\mathfrak{c}_-, \mathfrak{c}_+)$ for which the family $\{\mathfrak{L}_s\}_{s \in \mathbb{R}}$ has zero spectral flow. If $q(\cdot)$ is suitably generic, then $\mathcal{M}_{q(\cdot)0}(\mathfrak{c}_-, \mathfrak{c}_+)$ is a finite set of points, and each point comes with a sign. (See Chapter 25.2 in [KM1].) Let $\sigma(\mathfrak{c}_-, \mathfrak{c}_+)$ denote the sum of these signs, or else 0 if $\mathcal{M}_{q(\cdot)0}(\mathfrak{c}_-, \mathfrak{c}_+) = \emptyset$. Use these numbers to define a matrix $\mathbb{T}: \mathbb{V} \to \mathbb{V}$ as in (7.13) of [T1]. The results in Chapter 25.3 in [KM1] imply that $\mathbb{T}\delta_{m-1} = \delta_m \mathbb{T}$ and that $\mathbb{T}$ induces an isomorphism from the $\delta_{m-1}$ version of the Seiberg-Witten Floer homology to the $\delta_m$ version.

Granted the preceding, then the arguments given in Part 2 of [T1]'s Section 7c prove that $\mathbb{A}$ is upper triangular given the following: First, if $s \to \mathfrak{d}(s) \in \mathcal{M}_{q(\cdot)0}(\mathfrak{c}_-, \mathfrak{c}_+)$ and if the $s \to \infty$ limit of $\mathfrak{d}(s)$ is written as $\hat{u} \cdot \mathfrak{c}_+$ with $\hat{u} \in C^\infty(M; S^1)$, then $f(\hat{u} \cdot \mathfrak{c}_+) = f(\mathfrak{c}_-)$. Indeed, this condition defines what is required for membership in $\mathcal{M}_{q(\cdot)0}(\mathfrak{c}_-, \mathfrak{c}_+)$. Second, the total change of the function $s \to \mathfrak{a}(\mathfrak{d}(s))$ between $-\infty$ and $\infty$ must be negative unless $\mathfrak{c}_- = \mathfrak{c}_+$. This follows using the same argument as that given for the proof of Lemma 7.2 in [T1]. The point here, as in the proof of Lemma 7.2 in [T1], is that the function $\mathfrak{a}(\mathfrak{d}(\cdot))$ on $\mathbb{R}$ is non-increasing in the case when $s \to \mathfrak{d}(s)$ is a solution to (1.10), and it is constant if and only if $\mathfrak{d}(\cdot)$ is constant.

Proposition 2.3 allows for a canonical identification of the Seiberg-Witten Floer homology at different values of r in any given interval of $(1, \infty) - \{\rho_j\}$. This is done just as in Definition 3.13 of [T1] which treats the case when $c_1(\det(\mathbb{S}))$ is torsion. To say more, fix $\rho_i \in \{\rho_j\}$ and consider the interval $(\rho_i, \rho_{i+1})$. For $m \in \mathbb{J}(i)$, the Seiberg-Witten Floer homology is defined for $r \in [t_m, t_{m+1})$ by the differential $\delta_m$. The various $k \in \mathbb{Z}/p\mathbb{Z}$ versions of Proposition 2.3's matrix $\mathbb{A}$ are used to identify the $r \in [t_{m-1}, t_m)$ version of the complex with the version that is defined for $r \in [t_m, t_{m+1})$. This definition of the Seiberg-Witten Floer complex for the interval $(\rho_i, \rho_{i+1})$ is used implicitly in what follows when reference is made to the *Seiberg-Witten Floer complex* and the *Seiberg-Witten Floer homology* for the interval $(\rho_i, \rho_{i+1})$.

What follows is an analog of Definitions 4.1 and 4.4 of [T1].

**Definition 2.4**: *Fix $\mu$ as in Proposition 2.1, and let $\{\rho_j\} \subset (1, \infty)$ denote the resulting set. Fix $\rho_i \in \{\rho_j\}$ and let $\{t_n\}_{n \in \mathbb{J}(i)}$ be as described in Proposition 2.3. Fix $t_m \in \{t_n\}_{n \in \mathbb{J}(i)}$ and* *introduce the perturbation $\mathfrak{q}_m$ from Proposition 2.3. Given $\mathfrak{r} \in [t_m, t_{m+1})$, use $\mathfrak{r}$ and the perturbation $\mathfrak{g} = \mathfrak{e}_\mu + \mathfrak{q}_m$ to define the Seiberg-Witten Floer homology complex. Use the canonical labeling to identify the generators with the solutions to the $\mathfrak{r}$ and $\mu$ version of (1.3). Let $\theta$ denote a non-zero class for a given degree in the resulting Seiberg-Witten Floer homology.*

- *Suppose that $\mathfrak{n} = \sum_v z_v \mathfrak{c}_v$ is a cycle defined for the given value of $\mathfrak{r}$ and $\mu$ that represents the class $\theta$. Define $\hat{\mathfrak{a}}^f(\mathfrak{n}, \mathfrak{r})$ to be the maximum value of $\mathfrak{a}^f$ on the set of gauge equivalence classes of solutions to the $\mathfrak{r}$ and $\mu$ version of (1.3) that appear in the sum for $\mathfrak{n}$. Then define $\mathfrak{a}^f_\theta(\mathfrak{r})$ to be the minimum over all such $\mathfrak{n}$ of the values of $\hat{\mathfrak{a}}^f(\mathfrak{n}, \mathfrak{r})$.*

- *Suppose that $\mathfrak{n} = \sum_v z_v \mathfrak{c}_v$ is a cycle defined for the given value of $\mathfrak{r}$ and $\mu$ that represents the class $\theta$ and is such that $\hat{\mathfrak{a}}^f(\mathfrak{n}, \mathfrak{r}) = \mathfrak{a}^f_\theta(\mathfrak{r})$. Let $\hat{E}(\mathfrak{r}, \mathfrak{n})$ denote the infimum of the values of $E$ on the configurations $\mathfrak{c} \in \{\mathfrak{c}_v\}$ that appear in the sum for $\mathfrak{n}$ and have $\mathfrak{a}^f(\mathfrak{c}) = \mathfrak{a}^f_\theta(\mathfrak{r})$. Then, define*
   a) *$\hat{E}(\mathfrak{r})$ to be the infimum of the set $\{\hat{E}(\mathfrak{n}, \mathfrak{r})\}$ over all such $\mathfrak{n}$,*
   b) *$\mathfrak{v}(\mathfrak{r}) = 2\mathfrak{a}^f_\theta(\mathfrak{r}) + \mathfrak{r}\hat{E}(\mathfrak{r})$,*
   c) *$\mathfrak{f}(\mathfrak{r}) = -2\mathfrak{r}^{-1}\mathfrak{a}^f_\theta(\mathfrak{r}) = \hat{E}(\mathfrak{r}) - \mathfrak{r}^{-1}\mathfrak{v}(\mathfrak{r})$.*

The next proposition is the analog of Proposition 4.2 in [T1].

**Proposition 2.5**: *Fix $\mu$ as in Proposition 1.3, and let $\{\rho_j\} \subset (1, \infty)$ denote the resulting set. Then the various $\mathfrak{r} \in (1, \infty) - \{\rho_j\}$ versions of the Seiberg-Witten Floer homology can be identified so that the following is true: Let $\theta$ denote a Seiberg-Witten Floer homology class in a given degree. Then the function $\mathfrak{a}_\theta(\cdot) : (1, \infty) - \{\rho_j\}$ defines a continuous, piecewise differentiable function that extends to $(1, \infty)$ as a continuous piecewise differentiable function. Moreover, there exists a countable set in $(1, \infty)$ and a smooth map, $\mathfrak{c}(\cdot) = (A(\cdot), \psi(\cdot))$, from the complement of this set into $\mathrm{Conn}(E) \times C^\infty(M; \mathbb{S})$ such that for each $\mathfrak{r}$ in the domain of $\mathfrak{c}(\cdot)$, the configuration $\mathfrak{c}(\mathfrak{r})$ is a non-degenerate solution to the $\mathfrak{r}$ and $\mu$ version of (1.3) with the same degree as $\theta$ and with $\mathfrak{a}^f_\theta(\mathfrak{r}) = \mathfrak{a}^f(\mathfrak{c}(\cdot))$, $\hat{E} = E(A(\mathfrak{r}))$ and $\mathfrak{v} = \mathfrak{c}\mathfrak{s}^f(\mathfrak{c}(\mathfrak{r})) + 2\mathfrak{e}_\mu(\mathfrak{c}(\mathfrak{r}))$. In addition, $\mathfrak{c}(\mathfrak{r})$ defines a generator of the Seiberg-Witten Floer complex whose coefficient is non-zero in some cycle that represents the class $\theta$.*

### b) Proof of Proposition 2.5

There are four steps to the proof. The first three steps construct a suitable identification between the respective Seiberg-Witten Floer homologies as defined for different components of $(1, \infty) - \{\rho_j\}$. This is done with constructions that are very much like those given in Sections 3d and 3e of [T1]. Once these identifications are defined, the last step constitutes what is an essentially verbatim repeat of the arguments that are given in Section 4a of [T1] to prove the latter's Proposition 4.2.

   Step 1: Let $\rho_i \in \{\rho_j\}$. The constructions in Sections 3d and 3e of [T1] have their analogs here that provide a chain equivalence between the respective Seiberg-Witten Floer complexes as defined for the contiguous intervals $(\rho_{i-1}, \rho_i)$ and $(\rho_i, \rho_{i+1})$. As is explained below, almost verbatim analogs exist for each of five properties in Section 3d of [T1] including (3.6) from [T1] and Lemmas 3.14-3.17 from [T1]. This step provides some background stage setting and notation for the statements of these analogs.

   To start, suppose that $r \in (\rho_{i-1}, \rho_{i+1})$ and that some small element $\mathfrak{q} \in \mathcal{P}$ has been specified. Let $\mathfrak{c}$ denote a non-degenerate solution to the $r$ and $\mathfrak{g} = \mathfrak{e}_\mu + \mathfrak{q}$ version of (1.9). This is to say that the operator that is obtained from (1.7) by adding respective terms $\mathfrak{t}_{(A,\psi)}(b, \eta) = \frac{d}{dt}\mathfrak{T}(A+tb, \psi+t\eta)|_{t=0}$ and $\mathfrak{s}_{(A,\psi)} = \frac{d}{dt}\mathfrak{S}(A+tb, \psi+t\eta))|_{t=0}$ from what is written in the first and second lines of (1.7) has trivial kernel. Let $\mathfrak{L}_\mathfrak{q}$ denote this last operator. Define $f_\mathfrak{q}(\mathfrak{c})$ to be the spectral flow for a suitably generic path of self-adjoint operators on $L^2(M; iT^*M \oplus \mathbb{S} \oplus i\mathbb{R})$ that starts at the $(A_E, \psi_E)$ version of (1.7) and ends with $\mathfrak{L}_\mathfrak{q}$. Deem $-(f_\mathfrak{q}(\mathfrak{c}) + \frac{1}{4\pi^2}\mathfrak{cs}_E)$ to be the *degree* of $\mathfrak{c}$.

   If all solutions to the $r$ and $\mathfrak{g} = \mathfrak{e}_\mu + \mathfrak{q}$ version of (1.9) are non-degenerate, and if $\mathfrak{q}$ is suitably generic with regards to the behavior of the solutions space of the corresponding version of (1.10), then $\mathbb{Z}/p\mathbb{Z}$ graded Seiberg-Witten Floer homology can be defined as the homology of a square zero differential on the set of gauge equivalence classes of solutions to the $r$ and $\mathfrak{g} = \mathfrak{e}_\mu + \mathfrak{q}$ version of (1.9). Here, the grading of a solution $\mathfrak{c}$ is defined to be the $\mathbb{Z}/p\mathbb{Z}$ reduction of the degree of $\mathfrak{c}$ as just defined. Of course, if $r \in (\rho_{i-1}, \rho_i)$ of r if $r \in (\rho_i, \rho_{i+1})$ and if $\mathfrak{q}$ has sufficiently small norm, then the maps in Lemma 2.2 provide a canonical isomorphism from the p-periodic, $\mathbb{Z}$ graded complex as just defined to that defined by $\Lambda_-$ or $\Lambda_+$. This isomorphism identifies the respective differentials up to the action of some upper triangular matrix.

   Step 2: Fix $\varepsilon > 0$, fix $r_- \in (\rho_{i-1}, \rho_i)$, and fix $r_+ \in (\rho_i, \rho_{i+1})$. Let $m \in \mathbb{J}(i-1)$ be such that $r_- \in [t_m, t_{m+1})$, and set $\mathfrak{q}_- = \mathfrak{q}_m$. Let $m' \in \mathbb{J}(i)$ be such that $r_+ \in [t_{m'}, t_{m'+1})$ and set $\mathfrak{q}_+ = \mathfrak{q}_{m'}$. Given that both $r_-$ and $r_+$ sufficiently close to $\rho_i$, it can be assumed that both $\mathfrak{q}_-$ and $\mathfrak{q}_+$ are in the radius $\varepsilon$ ball about 0 in $\mathcal{P}$. The task at hand is to choose a path $r \to \mathfrak{q}(r)$ in this ball with certain desired properties. The path is parameterized by $r \in [r_-, r_+]$, it obeys $\mathfrak{q}(r_-) = \mathfrak{q}_-$ and $\mathfrak{q}(r_+) = \mathfrak{q}_+$. As explained below, if $|\rho_i - r_-|$ and $|\rho_i - r_+|$ are sufficiently small, the path can be chosen to have the five properties listed next.

*Property 1*: Let $\mathfrak{g}(r) = \mathfrak{e}_\mu + \mathfrak{q}(r)$, and let $\mathfrak{a}_{\mathfrak{g}(r)}$ denote the functional on Conn(E) × $C^\infty(M; S^1)$ that is obtained from (1.4) by replacing $\mathfrak{e}_\mu$ by $\mathfrak{g}(r)$. Fix $r \in (t_m, t_{m'})$ and any solution to the $r$ and $\mathfrak{g}(r)$ version of (1.9). The value of $\mathfrak{a}_{\mathfrak{g}(r)}$ on the solution is within $\epsilon^2$ of the value of the original version of $\mathfrak{a}$ on some solution to the $r = \rho_i$ and $\mu$ version of (1.3). Moreover, there is a finite, increasing subset, $\{y_n\} \subset (r_-, r_+)$, such that all solutions to the r and $\mathfrak{g}(r)$ version of (1.9) are non-degenerate when $r \notin \{y_n\}$. In addition, the values of the functional $\mathfrak{a}_{\mathfrak{g}(r)} - 2\pi^2 f_{\mathfrak{g}(r)} - \frac{1}{2}\mathfrak{cs}_E$ distinguish the various gauge equivalence classes of solutions to the r and $\mathfrak{g}(r)$ version of (1.9) when $r \notin \{y_n\}$.

*Property 2*: Let $I \subset [r_-, r_+] - \{y_n\}$ denote a component. There exists a consecutively labeled, increasing set, $\{w_n\}_{n \in \mathbb{K}(I)}$, in the interior of I that is finite or countable, but with no accumulation points in I. For each $m \in \mathbb{K}(I)$, there exists a perturbation $\mathfrak{p}_m \in \mathcal{P}$ of very small norm such that $(\mu, \mathfrak{q}(r) + \mathfrak{p}_m)$ is (k, r)-admissable at each $r \in [w_m, w_{m+1}]$. In addition, $\mathfrak{p}_m$ is such that the gauge equivalence classes of solutions to the r and $\mathfrak{g}(r, m) = \mathfrak{e}_\mu + \mathfrak{q}(r) + \mathfrak{p}_m$ version of (1.9) are in 1-1 correspondence with those of the r and $\mathfrak{g}(r)$ version of (1.9) with the same degree for all $r \in [t_m, t_{m+1}]$. This equivalence is given by the analog of the maps in Lemma 2.2 (see Lemma 3.2 in [T1]). The equivalence is such that the value of $\mathfrak{a}_{\mathfrak{g}(r,m)} - 2\pi^2 f_{\mathfrak{g}(r,m)} - \frac{1}{2}\mathfrak{cs}_E$ on a given gauge equivalence class of r and $\mathfrak{g}(r, m)$ solutions to (1.9) is very much closer to the value of the function $\mathfrak{a}_{\mathfrak{g}(r)} - 2\pi^2 f_{\mathfrak{g}(r)} - \frac{1}{2}\mathfrak{cs}_E$ on its partner gauge equivalence class of r and $\mathfrak{g}(r)$ solutions to (1.9) then it is to the value of $\mathfrak{a}_{\mathfrak{g}(r)} - 2\pi^2 f_{\mathfrak{g}(r)} - \frac{1}{2}\mathfrak{cs}_E$ on any other gauge equivalence class of solution to the r and $\mathfrak{g}(r)$ version of (1.9). In particular, the ordering of the r and $\mathfrak{g}(r)$ equivalence classes of solutions given by the values $\mathfrak{a}_{\mathfrak{g}(r)} - 2\pi^2 f_{\mathfrak{g}(r)} - \frac{1}{2}\mathfrak{cs}_E$ is the same as that defined by $\mathfrak{a}_{\mathfrak{g}(r,m)} - 2\pi^2 f_{\mathfrak{g}(r,m)} - \frac{1}{2}\mathfrak{cs}_E$ via this equivalence.

As in the cases studied by [T1], what is asserted by Properties 1 and 2 have the following consequence: Fix $I \subset [r_-, r_+] - \{y_n\}$ and $m \in \mathbb{K}(I)$. Then the Seiberg-Witten Floer homology can be defined for $r \in [w_m, w_{m+1}]$ using the solutions to the r and $\mathfrak{g}(r, m)$ versions of (1.9) and (1.10). Note in this regard that the vector space of cycles in a given degree can be identified using Property 2 with a fixed vector space, this the vector space that is generated by the gauge equivalence classes of solutions to the r and $\mathfrak{g}(r)$ version of (1.9), and with these generaters labeled by their ordering using $\mathfrak{a}_{\mathfrak{g}(r)} - 2\pi^2 f_{\mathfrak{g}(r)} - \frac{1}{2}\mathfrak{cs}_E$. Here, the convention is to label the basis of cycles with the larger numbered ones having smaller values of $\mathfrak{a}_{\mathfrak{g}(r)}$. This fixed, r-independent basis is called the I-canonical basis.

*Property 3*: Fix an interval $I \subset [r_-, r_+] - \{y_n\}$ and some $w_m \in \mathbb{K}(I)$. As r varies in $[w_m, w_{m+1}]$, the differentials as written for the I-canonical basis of the p-periodic, $\mathbb{Z}$-graded, Seiberg-Witten Floer complex are independent of r. Moreover, there exists an upper triangular, degree preserving matrix, $\mathbb{A} = \mathbb{A}(m)$ with 1's on the diagonal such that the differential, $\delta_{m-1}$ defined on $[w_{m-1}, w_m]$ and the differential $\delta_m$ defined on $[w_m, w_{m+1}]$ are related, when written using the $\mathbb{I}$-canonical basis, by the rule $\delta_m = \mathbb{A}^{-1} \delta_m \mathbb{A}$.

The next property addresses behavior of the solutions to the r and $\mathfrak{g}(r)$ version of (2.4) at any given $y \in \{y_n\}$. In what follows, $I_-$ denotes the component of $(r_-, r_+) - \{y_n\}$ whose closure adds y as its upper endpoint, and $I_+$ denotes the the component whose closure adds y as its lower endpoint.

*Property 4*: One and only one of the following two assertions holds:

- *All solutions to the $r = y$ and $\mathfrak{g}(y)$ version of (1.9) are non-degenerate, and there is precisely one pair of distinct gauge equivalence classes of solutions to the $r = y$ and $\mathfrak{g}(y)$ version of (1.9) that are not distinguished by the values of $\mathfrak{a}_{\mathfrak{g}(y)} - 2\pi^2 f_\mathfrak{g} - \frac{1}{2} \mathfrak{cs}_E$. In addition, there exist $y_- < y$ and $y_+ > y$ such that if $y_0 \in [y_-, y_+]$ and if $\mathfrak{c}$ is a solution to the $r = y_0$ and $\mathfrak{g}(y_0)$ version of (1.9), then there is a smooth map, $\mathfrak{c}(\cdot)$, from $(y_-, y_+)$ to $\mathrm{Conn}(E) \times C^\infty(\mathbb{S})$ such that $\mathfrak{c}(y_-) = \mathfrak{c}$ and $\mathfrak{c}(r)$ solves the r and $\mathfrak{g}(r)$ version of (1.3) for each $r \in [y_-, y_+]$.*
- *The function $\mathfrak{a}_{\mathfrak{g}(y} - 2\pi^2 f_{\mathfrak{g}(y)} - \frac{1}{2} \mathfrak{cs}_E$ distinguishes the gauge equivalence classes of solution to the $r = y$ and $\mathfrak{g}(y)$ version of (1.9). Meanwhile all but one gauge-equivalence class of solution to the $r = y$ and $\mathfrak{g}(y)$ version of (1.9) has non-degenerate solutions. In addition,*
    1) *The operator $\mathfrak{L}_{\mathfrak{g}(y)}$ for any solution in the one anomalous gauge equivalence class has kernel dimension 1.*
    2) *The number of gauge equivalence classes of solutions to the r and $\mathfrak{g}(r)$ version of (1.9) change by two as r crosses y, and the number of gauge equivalence classes of solutions to the $r = y$ and $\mathfrak{g}(y)$ version of (1.9) differs by 1 from the number on either side of y.*
    3) *Let $I \in \{I_-, I_+\}$ denote the component with the greater number of equivalence classes. Then there are respective representatives, $\mathfrak{c}(r)$ and $\mathfrak{c}'(r)$, of distinct equivalence classes of solutions to the r and $\mathfrak{g}(r)$ version of (1.9) that vary smoothly with $r \in I$ and converge in $\mathrm{Conn}(E) \times C^\infty(\mathbb{S})$ as $r \to y$ to the one anomalous $r = y$ equivalence class. Also, the $\mathbb{Z}/p\mathbb{Z}$ degree of $\mathfrak{c}$ is one greater than that of $\mathfrak{c}'$.*
    4) *Let $\mathfrak{n}$ denote a solution to the $r = y$ and $\mathfrak{g}(y)$ version of (1.9) that is not gauge equivalent to the one anomalous gauge equivalence class.. Then there is a smooth map $\mathfrak{n}(\cdot)\colon I_- \cup \{y\} \cup I_+ \to \mathrm{Conn}(E) \times C^\infty(\mathbb{S})$ such that $\mathfrak{n}(y) = \mathfrak{n}$, and such that $\mathfrak{n}(r)$ is a solution to the r and $\mathfrak{g}(r)$ version of (1.9) for all $r \in I_- \cup \{y\} \cup I_+$.*

The final property describes how the generators of the Seiberg-Witten Floer homology change as r crosses a given $y \in \{y_n\}$. To this end, define the respective $I_-$ and $I_+$ versions of the cSWF complex and homology in degrees greater than k using the points $y_-$ and $y_+$. This is to say that $y_-$ is in some $I_-$ version of $[w_m, w_{m+1}]$, and use the corresponding $r = y_-$ and $\mathfrak{g}(r, m)$ to define the Seiberg-Witten Floer homology in degrees greater than k using these points. There are three parts to property 5 that address three cases that are consistent with what is described in Property 4..

*Property 5a*: Assume here that the top bullet in Property 4 is relevant for y. Use the maps $\mathfrak{c}(\cdot)$ to extend the $I_+$-canonical basis as defined at $y_+$ to give a new basis for the Seiberg-Witten Floer complex at $y_-$. Let $\mathfrak{c}$ and $\mathfrak{c}'$ denote the two generators that are not distinguished by $\mathfrak{a}_{\mathfrak{g}(y)} - 2\pi^2 f_{\mathfrak{g}(y)} - \frac{1}{2} \mathfrak{cs}_E$. If $\mathfrak{c}$ and $\mathfrak{c}'$ have different degrees, then this new basis at $y_-$ is the same as the $I_-$-canonical basis. If $\mathfrak{c}$ and $\mathfrak{c}'$ have the same degree, make the convention that $\mathfrak{c}(y_+) = \mathfrak{c}_n$ and $\mathfrak{c}'(y_+) = \mathfrak{c}_{n+1}$ where $\mathfrak{c}_n$ and $\mathfrak{c}_{n+1}$ are $I_+$-canonical basis elements at $y_+$. With respect to the $I_-$-canonical basis at $y_-$, either $\mathfrak{c}(y_-) = \mathfrak{c}_n$ and $\mathfrak{c}'(y_-) = \mathfrak{c}_{n+1}$, or else $\mathfrak{c}(y_-) = \mathfrak{c}_{n+1}$ and $\mathfrak{c}'(y_-) = \mathfrak{c}_n$. If the labelings do not change, then the respective $I_-$ and $I_+$ canonical basis for Seiberg-Witten Floer complexes as defined at $y_-$ and $y_+$ agree. If these canonical basis agree, either for this reason, or because $\mathfrak{c}$ and $\mathfrak{c}'$ have distinct degrees, then the differential, $\delta_-$, at $y_-$ is related to the differential, $\delta_+$, defined at $y_+$ as follows: $\delta_+ = \mathbb{A}^{-1} \delta_- \mathbb{A}$, where $\mathbb{A}$ is a degree preserving, upper triangular matrix with 1's on the diagonal.

Suppose now that $\mathfrak{c}$ and $\mathfrak{c}'$ have the same degree and the labelings change as r crosses y. Let d denote the degree of $\mathfrak{c}$ and $\mathfrak{c}'$. In this case, the differentials are again related by $\delta_+ = \mathbb{A}^{-1} \delta_- \mathbb{A}$, where $\mathbb{A}$ is a degree preserving matrix of the following sort: In degrees not equal to d, the matrix $\mathbb{A}$ is upper triangular with 1's on the diagonal. In degree d, the matrix $\mathbb{A}$ is such that $\mathbb{A}_{n,n} = \mathbb{A}_{n+1,n+1} = 0$, $\mathbb{A}_{n,n+1} = \mathbb{A}_{n+1,n} = 1$, $\mathbb{A}_{v,v} = 1$ if $v \neq n$ or $n+1$, and $\mathbb{A}_{v,v'} = 0$ *if* $v > v'$ *and* $(v, v') \neq (n+1, n)$.

Properties 5b and 5c assume that the second bullet of Property 4 describes the situation at y. In what follows, $\mathfrak{c}$ and $\mathfrak{c}'$ denote respective representatives of the two equivalence classes that do not extend across y; and let d+1 and d denote their respective $\mathbb{Z}/p\mathbb{Z}$ degrees. The maps that are supplied by the the fourth item of the second bullet of Property 4 are used to identify the remaining generators for the $I_-$-canonical basis at $y_-$ with the generators of the $I_+$-canonical basis at $y_+$. This identifies the full $I_-$-canonical basis at $y_-$ with the full $I_+$-canonical basis at $y_+$ in degrees different from d+1 and d, and does so as the identity map.

*Property 5b*: Assume here that $\mathfrak{c}$ and $\mathfrak{c}'$ are defined for $r \in I_-$. In degree d+1, the canonical basis at $y_+$ is obtained from that at $y_-$ by deleting the generator $\mathfrak{c}$; and in degree d, the change is deletion of the generator $\mathfrak{c}'$. Note that this identification preserves the ordering given by the value of $\mathfrak{a}_{g(r)} - 2\pi^2 f_{g(r)} - \frac{1}{2}\mathfrak{cs}_E$. Let $\mathbb{V}_+$ denote the vector space of cycles as defined for $y_+$. With the preceding identifications understood, the vector space of cycles for $y_-$ is then $\mathbb{Z}\mathfrak{c} \oplus \mathbb{Z}\mathfrak{c}' \oplus \mathbb{V}_+$. Let $\delta_+$ denote the Seiberg-Witten Floer differential on $\mathbb{V}_+$ and let $\delta_-$ denote that on $\mathbb{Z}\mathfrak{c} \oplus \mathbb{Z}\mathfrak{c}' \oplus \mathbb{V}_+$. There is a degree preserving homomorphism, $\mathbb{T}: \mathbb{Z}\mathfrak{c} \oplus \mathbb{Z}\mathfrak{c}' \oplus \mathbb{V}_+ \to \mathbb{V}_+$ with the following properties:

- $\mathbb{T}\delta_- = \delta_+\mathbb{T}$.
- $\mathbb{T}$ *induces an isomorphism on homology*
- $\mathbb{T}$ *maps $\mathbb{V}_+$ to itself as an upper triangular matrix with 1's on the diagonal.*
- *The value of* $\mathfrak{a}_{g(y)} - 2\pi^2 f_{g(y)} - \frac{1}{2}\mathfrak{cs}_E$ *on any generator that appears in $\mathbb{T}\mathfrak{c}$ is less than its value on $\mathfrak{c}$.*
- *The value of* $\mathfrak{a}_{g(y)} - 2\pi^2 f_{g(y)} - \frac{1}{2}\mathfrak{cs}_E$ *on any generator that appears in $\mathbb{T}\mathfrak{c}'$ is less than its value on $\mathfrak{c}'$.*

*Property 5c*: Assume here that $\mathfrak{c}$ and $\mathfrak{c}'$ are defined for $r \in I_+$. Let $\mathbb{V}_-$ denote the vector space of cycles as defined at $y_-$. With the aforementioned identifications understood, the vector space of cycles at $y_+$ is $\mathbb{Z}\mathfrak{c} \oplus \mathbb{Z}\mathfrak{c}' \oplus \mathbb{V}_-$. Let $\delta_-$ denote the Seiberg-Witten Floer differential on $\mathbb{V}_-$ and let $\delta_+$ denote the differential on $\mathbb{Z}\mathfrak{c} \oplus \mathbb{Z}\mathfrak{c}' \oplus \mathbb{V}_-$. There is a degree preserving homomorphism $\mathbb{T}: \mathbb{V}_- \to \mathbb{Z}\mathfrak{c} \oplus \mathbb{Z}\mathfrak{c}' \oplus \mathbb{V}_-$ with the following properties:

- $\mathbb{T}\delta_- = \delta_+\mathbb{T}$.
- $\mathbb{T}$ *induces an isomorphism on homology*
- $\mathbb{T}$ *is upper triangular with ones on the diagonal in degrees different from d+1 and d,*
- *If $\mathfrak{u}$ has degree d+1, then $\mathbb{T}\mathfrak{u} = \mathbb{A}\mathfrak{u} + \kappa_\mathfrak{u}\mathfrak{c}$ where $\mathbb{A}: \mathbb{V}_+ \to \mathbb{V}_+$ is an upper triangular matrix with 1's on the diagonal. Here, $\kappa_\mathfrak{u} = 0$ for a generator $\mathfrak{u}$ if the value of $\mathfrak{a}_{g(y)} - 2\pi^2 f_{g(y)} - \frac{1}{2}\mathfrak{cs}_E$ is less than its value on $\mathfrak{c}$.*
- *If $\mathfrak{v}$ has degree d, then $\mathbb{T}\mathfrak{v} = \mathbb{A}\mathfrak{v} + \kappa_\mathfrak{v}\mathfrak{c}'$ where $\mathbb{A}: \mathbb{V}_+ \to \mathbb{V}_+$ is an upper triangular matrix with 1's on the diagonal. Here, $\kappa_\mathfrak{v} = 0$ for a generator $\mathfrak{v}$ if the value of $\mathfrak{a}_{g(y)} - 2\pi^2 f_{g(y)} - \frac{1}{2}\mathfrak{cs}_E$ on $\mathfrak{v}$ is less than its value on $\mathfrak{c}'$.*

But for two extra remarks, the arguments in Section 7d of [T1] can be used to prove that Properties 1-5 can be satisfied. These remarks concerns the matrices $\mathbb{A}$ and $\mathbb{T}$ that appear in the $\mathfrak{a}^f$ versions of what is written above.

Here is the first remark: The arguments in Section 7d of [T1] use the solutions to (7.23) of [T1] to produce a degree preserving homomorphism that relates the respective vector spaces over $\mathbb{Z}$ of cycles for the Seiberg-Witten Floer complex as defined by a given value of r and distinct pairs, $\mathfrak{g}_- = \mathfrak{e}_\mu + \mathfrak{q}_-$ and $\mathfrak{g}_+ = \mathfrak{e}_\mu + \mathfrak{q}_+$ from $\mathcal{P}$. This matrix, either $\mathbb{A}$ or $\mathbb{T}$ as the case may be, is denoted here as $\mathbb{T}$. If $\mathfrak{c}_-$ and $\mathfrak{c}_+$ are respective solutions to the corresponding r and $\mathfrak{g}_\pm$ versions of (1.9), then the component along the basis element $\mathfrak{c}_+$ of $\mathbb{T}\mathfrak{c}_-$ is obtained by counting with signs a certain set of instanton solutions to a version of (7.23) in [T1]; this version is defined by a suitably generic map s $\to \mathfrak{q}(s)$ from $\mathbb{R}$ into $\mathcal{P}$ with the property that $\mathfrak{q}(s) = \mathfrak{q}_-$ for $s < -1$ and $\mathfrak{q}(s) = \mathfrak{q}_+$ for $s > 1$. As in Section 7d of [T1], the instantons that contribute to the count have two salient features. To state them, write the instanton as $s \to \mathfrak{d}(s) = (A(s), \psi(s))$. Also, for each $s \in \mathbb{R}$, let $\mathfrak{L}_{\mathfrak{q}(s)}$ denote the $\mathfrak{q} = \mathfrak{q}(s)$ and $(A, \psi) = (A(s), \psi(s))$ version of the operator $\mathfrak{L}_\mathfrak{q}$ as described in Step 1. The first requirement for $\mathfrak{d}(\cdot)$ demands that $\lim_{s \to -\infty} \mathfrak{d}(s) = \mathfrak{c}_-$ and that $\lim_{s \to \infty} \mathfrak{d}(s) = \hat{u} \cdot \mathfrak{c}_+$ where $\hat{u}$ is a smooth map from M to $S^1$. The second requirement demands that the family $\{\mathfrak{L}_{\mathfrak{q}(s)}\}_{s \in \mathbb{R}}$ have zero spectral flow. This understood, it then follows that the $\mathfrak{c}_+$ component of $\mathbb{T}\mathfrak{c}_-$ is non-zero only if the $\mathfrak{q} = \mathfrak{q}_-$ version of $f_\mathfrak{q}(\mathfrak{c}_-)$ is equal to the $\mathfrak{q} = \mathfrak{q}_+$ version of $f_\mathfrak{q}(\hat{u} \cdot \mathfrak{c}_+)$. Meanwhile, with $\mathfrak{g}(s)$ used to denote $\mathfrak{e}_\mu + \mathfrak{q}(s)$, the behavior of the function $\mathfrak{a}_{\mathfrak{g}(\cdot)}(\mathfrak{d}(\cdot))$ on $\mathbb{R}$ is just as it was in the cases that are considered in Section 7d. These last observations are needed to justify various assertions that $\mathbb{T}$ or parts of $\mathbb{T}$ are upper triangular with 1's on the diagonal.

To begin the second remark, recall that when $c_1(\det(\mathbb{S}))$ is torsion, the differential on the Seiberg-Witten Floer complex has the following property: If $\mathfrak{c}_-$ and $\mathfrak{c}_+$ are solutions to some r and $\mathfrak{g} = \mathfrak{e}_\mu + \mathfrak{q}$ version of (1.9) (this is (2.4) in [T1]), and if, when viewed as basis vectors in the cSWF complex, the solution $\mathfrak{c}_+$ appears with a non-zero weight in $\delta\mathfrak{c}_-$, then $\mathfrak{a}(\mathfrak{c}_-) > \mathfrak{a}(\mathfrak{c}_+)$. This is because differential for the Seiberg-Witten Floer complex was defined by counting a certain set of non-constant solutions to (1.10) (this is (2.11) in [T1]) and the equations in (1.10) are gradient flow equations for the function $\mathfrak{a}_\mathfrak{g}$ given in (2.9) of [T1].

For the case at hand, the coefficient for $\mathfrak{c}_+$ in $\delta\mathfrak{c}_-$ is again obtained as an algebraic count of a certain set of solutions to (1.10). Indeed, recall that this coefficient is the count with ±1 weights of the components of a set, $\mathcal{M}_{t=1}(\mathfrak{c}_-, \mathfrak{c}_+)$, that consists of the instantons solutions to the r and $\mathfrak{g}$ version (1.10) that have the following two properties: Write the instanton as a map $s \to \mathfrak{d}(s) = (A(s), \psi(s))$. First, $\lim_{s \to -\infty} \mathfrak{d}(s) = \mathfrak{c}_-$ and $\lim_{s \to \infty} \mathfrak{d}(s) = \hat{u} \cdot \mathfrak{c}_+$ where $\hat{u}$ is a smooth map from M to $S^1$. To state the second property, define, for each $s \in \mathbb{R}$, the operator $\mathfrak{L}_{\mathfrak{q}(s)}$ to be the r and $\mathfrak{q}$ version of Step 1's operator $\mathfrak{L}_\mathfrak{q}$ as defined using the pair $(A, \psi) = (A(s), \psi(s))$. Then $\mathfrak{d}$ is in $\mathcal{M}_1(\mathfrak{c}_-, \mathfrak{c}_+)$ if and only if the spectral flow for the family $\{\mathfrak{L}_{q(s)}\}_{s \in \mathbb{R}}$ is equal to 1. This then means that $f_q(\mathfrak{c}_-) = f_q(\hat{u} \cdot \mathfrak{c}_+) + 1$. Meanwhile, it is still the case that the equations in (1.10) for a map $s \to \mathfrak{d}(s) = (A(s), \psi(s))$ are the gradient flow equations for the function $\mathfrak{a}_\mathfrak{g}$, this now defined defined as in (1.4) with $\mathfrak{g} = \mathfrak{e}_\mu + \mathfrak{q}$ replacing $\mathfrak{e}_\mu$. This being the case, the function $\mathfrak{a}_\mathfrak{g}(\mathfrak{d}(\cdot))$ on $\mathbb{R}$ is non-increasing and constant if and only if $\mathfrak{d}(\cdot)$ is constant. Granted all of this, it follows that the coefficient that multiplies $\mathfrak{c}_+$ in the basis expansion of $\delta \mathfrak{c}_-$ is non-zero only if the value on $\mathfrak{c}_-$ of the function $\mathfrak{a}_\mathfrak{g} - 2\pi^2 f_\mathfrak{g} - \frac{1}{2} \mathfrak{cs}_E$ is at least $2\pi^2$ more than the value of this function on $\mathfrak{c}_+$.

<u>Step 3</u>: Suppose here that $\varepsilon$, $r_-$, $r_+$ and a path $q(\cdot)$ have been chosen so as to satisfy the five properties stated in Step 2. Let $\{y_n\}$ be as described in Property 1, and suppose that $y \in \{y_n\}$. Assume in what follows that Property 5b is the relevant part of Property 5 for y. With the notation from Step 2 understood, what follows is the analog here of what is stated by Lemma 3.16 in [T1].

**Lemma 2.6**: *Let $\mathfrak{u} \in \mathbb{V}_+$ denote the class such that $\mathbb{T}\mathfrak{u} = \mathbb{T}\mathfrak{c}$ and let $\mathfrak{v} \in \mathbb{V}_+$ denote the class such that $\mathbb{T}\mathfrak{v} = \mathbb{T}\mathfrak{c}'$. If $p > 2$, there exists $A \in \{\pm 1\}$ such that $\delta_-(\mathfrak{c} - \mathfrak{u}) = A(\mathfrak{c}' - \mathfrak{v})$. As a consequence, there exists $\mathfrak{n} \in \mathbb{V}_+$ of degree d such that $\delta_-\mathfrak{c} = A\mathfrak{c}' + \mathfrak{n}$ and such that the value of $\mathfrak{a}_{\mathfrak{g}(y)} - 2\pi^2 f_{\mathfrak{g}(y)} - \frac{1}{2} \mathfrak{cs}_E$ on the generators that appear in $\mathfrak{n}$ with non-zero coefficient is no greater than its value on $\mathfrak{c}'$. If $p = 2$, either the preceding conclusion holds as stated, or it holds after switching the roles of $\mathfrak{c}$ and $\mathfrak{c}'$.*

*Proof of Lemma 2.6*: In what follows, recall that p denotes the greatest divisor of $c_1(\det(\mathbb{S}))$. Let $\mathfrak{v} \in \mathbb{V}_+$ denote the class with degree d such that $\mathbb{T}\mathfrak{v} = \mathbb{T}\mathfrak{c}'$. Let $\mathfrak{u}$ denote the class of degree d+1 such that $\mathbb{T}\mathfrak{u} = \mathbb{T}\mathfrak{c}$. If $p > 2$, then $\delta_-(\mathfrak{c}' - \mathfrak{v}) = 0$, and if $p = 2$, there exists $B \in \mathbb{Z}$ such that $\delta_-(\mathfrak{c}' - \mathfrak{v}) = B(\mathfrak{c} - \mathfrak{u})$. Suppose first that $B = 0$. Granted this, then no matter the value of p, the class $\mathfrak{c}' - \mathfrak{v}$ must equal $\delta_-(\mathfrak{w} + K\mathfrak{c})$ for some $K \in \mathbb{Z}$ and $\mathfrak{w} \in \mathbb{V}_+$. Indeed, were this not the case, then $\mathbb{T}$ could not induce an isomorphism on homology. Since $\delta_+\mathbb{T}(\mathfrak{w} + K\mathfrak{c}) = 0$, it follows that $\mathbb{T}(\mathfrak{w} + K\mathfrak{c}) = \delta_+\mathbb{T}(\mathfrak{o})$ for some class $\mathfrak{o} \in \mathbb{V}_+$ of degree $(d + 2) \bmod(p)$. Again, this is necessary if $\mathbb{T}$ induces an isomorphism on homology. As $\delta_-\mathfrak{o} = \mathfrak{w} + K\mathfrak{c} + A(\mathfrak{u} - \mathfrak{c})$ with $A \in \mathbb{Z}$, this implies that $A\delta_-(\mathfrak{c} - \mathfrak{u}) = \mathfrak{c}' - \mathfrak{v}$. As $\mathfrak{c}'$ is a generator, so $A = \pm 1$ and therefore $\delta_-\mathfrak{c} = A\mathfrak{c}' + \mathfrak{n}$ where $\mathfrak{n} = \mathfrak{u} - A\mathfrak{v}$. The claim about the relative values of the function $\mathfrak{a}_{\mathfrak{g}(y)} - 2\pi^2 f_{\mathfrak{g}(y)} - \frac{1}{2} \mathfrak{cs}_E$ follows from the last two points in Property 5b given that the respective values $\mathfrak{a}_{\mathfrak{g}(\cdot)}$ and $-f_{\mathfrak{g}(\cdot)}$ on the $s \to -\infty$ limit of an instanton are no greater than their values on the $s \to \infty$ limit of an instanton.

Now suppose that $B \neq 0$. In this case, $\delta_-(\mathfrak{c} - \mathfrak{u}) = 0$ and so the argument just given can be repeated verbatim after switching the roles of $\mathfrak{c}$ and $\mathfrak{c}'$.

Step 4: Given what is said in the previous steps, the discussion in Part 4a of [T1] and its proof of Proposition 4.2 in [T1] can be repeated with only minor cosmetic changes to identify the Seiberg-Witten Floer homology for different values of r in $[1, \infty)-\{\rho_j\}$ so as to obtain a continuous, piece-wise differentiable function, $r \to \mathfrak{a}^f_\theta(r)$ that is defined on the whole of $[1, \infty)$. The existence of the solution $\mathfrak{c}(r)$ with the stated properties follows from the min-max definition of $\mathfrak{a}^f_\theta$.

### c) Final arguments for Theorem 1.1

Proposition 2.5 is used in what follows to identify the respective Seiberg-Witten Floer homology groups at distinct value of $r \in [1, \infty)-\{\rho_j\}$. With these identifications from Proposition 2.5 understood, say that a class $\theta$ in the Seiberg-Witten Floer homology in a given degree is a <u>divergence</u> class when the following is true: Given $E > 0$, there exists $\rho_E \in (1, \infty)$ such that Definition 2.4's function $\hat{E}(r)$ is greater than E when $r > \rho_E$. Granted this definition, what follows is a crucial part of the story.

**Proposition 2.7**: *Fix $\mu$ as in Proposition 2.1. Suppose that $\theta$ is a divergence class. The class $\theta$ determines a constant, $c > 0$, with the following significance: Fix $r´ > 1$ and there exists $r > r´$ and a solution $(A, \psi)$ to the version of (1.3) determined by r and $\mu$ that has the same degree as $\theta$ and is such that $\mathfrak{cs}^f(A, \psi) > c r^2$ and $E(A) > c r$.*

*Proof of Proposition 2.7:* But for two points, the proof of this proposition differs in no essential aspects from the proofs of Proposition 4.6 and Corollary 4.7 in [T1]. Here is the first point: Suppose that r is in the domain of some map $\mathfrak{c}(\cdot)$ from Lemma 2.2. Then $f(\mathfrak{c}(\cdot))$ is constant on some neighborhood of r, and so $\frac{d}{dr}\mathfrak{a}^f = -\frac{1}{2}\hat{E}$. There is much more to the second point. This concerns the bound given in (4.9) of [T1] for $\mathfrak{cs}(A)$. In the case at hand, this bound must be replaced by a bound for $\mathfrak{cs}^f$. Proposition 1.9 supplies a useable bound.

Proposition 2.7 and Proposition 1.10 together imply that there are no divergence classes in the Seiberg-Witten Floer homology. To elaborate, suppose to the contrary that $\theta$ is a divergence class. Take $(A, \psi)$ as in Proposition 2.7 for some very large value of r. Proposition 1.10 asserts that $|\mathfrak{cs}^f| \leq c_0 r^{31/16}$ when r is large. Meanwhile, Proposition 2.7 asserts that $|\mathfrak{cs}^f| \geq c_0 r^2$; so these assertions are incompatible.

As a consequence, given a non-zero Seiberg-Witten Floer homology class, there is a form $\mu \in \Omega$ with $C^3$ norm less than 1 and a sequence $\{(r_n, (A_n, \psi_n))\}_{n=1,2,...} \subset (1, \infty) \times (\text{Conn}(E) \times C^\infty(M; \mathbb{S}))$ with the following properties: First, $\{r_n\}$ is increasing and unbounded from above. Second, $(A_n, \psi_n)$ solves the $r = r_n$ and $\mu$ version of (1.3). Third, the sequence $\{E(A_n)\} \subset \mathbb{R}$ is bounded.

If $c_1(E) \neq 0$, then $\{\sup_M (1 - |\psi_n|)\}$ also has a positive lower bound. To see why, note first that $1 - |\psi_n| \geq 1 - |\alpha_n| - c_0 r_n^{-1/2}$ by virtue of what is said in Lemma 1.6. As $c_1(E)$ is not zero, and as $\alpha_n$ is a section of E, there exist points in M where $|\alpha_n|$ vanishes. Thus $\sup_M (1 - |\psi_n|) \geq 1 - c_0 r_n^{-1/2}$. Granted these bounds, the statement of Theorem 1.1 in the case $c_1(E) \neq 0$ follows directly from the statement of Theorem 1.5.

Consider next the case where $c_1(E) = 0$. The argument used when $c_1(E) \neq 0$ to find a positive lower bound for $\{\sup_M(1 - |\psi_n|)\}$ won't suffice because E is the trivial bundle. In fact, there is a solution for all $r \geq 1$ with bounded energy to a special version of (1.3) for which the corresponding values of $\sup_M(1 - |\psi|)$ limit to zero. To explain, reintroduce the section, $1_\mathbb{C}$, of E with norm equal to 1. Use this section to trivialize E and let $A_I$ denote the corresponding product connection. The pair $(A_I, \psi_I = (1_\mathbb{C}, 0))$ has $\mathcal{E}(A) = 0$; and for any $r \geq 1$, this pair solves the equations

- $B_A = r(\psi^\dagger \tau \psi - ia)$.
- $D_A \psi = 0$.

(2.2)

As the next proposition shows, this solution to (2.2) has a counter-part that solves (1.3).

**Proposition 2.8**: *There exists $r_I > 1$ and $\delta \in (0, \frac{1}{2})$ with the following significance: Fix $\mu \in \Omega$ with $C^3$ norm less than 1 and fix $r \geq r_I$. Then there exists a unique gauge equivalence class of solutions to* (1.3) *with the norm of the spinor component nowhere less than $1 - \delta$. Moreover, $\mathcal{E}$ has an $r$ independent upper bound on this equivalence class.*

This proposition is proved momentarily.

Granted this proposition, suppose that the Seiberg-Witten Floer homology has two or more generators. As at least one of them won't be the gauge equivalence class supplied by Proposition 2.8, there is a sequence $\{(r_n, (A_n, \psi_n))\}$ that satisfies the assumptions of Theorem 1.5. Thus, the assertion of Theorem 1.1 follows in this case.

Suppose instead that the Seiberg-Witten Floer homology vanishes. Fix $\mu$ as in Proposition 2.1 and let $\{\rho_j\} \subset (r_I, \infty)$ denote the set supplied by this same proposition. Let $\rho_i \in \{\rho_j\}$. With Proposition 2.3 in mind, let $t_m \in \{t_n\}_{n \in \mathbb{J}(i)}$ and let $r \in [t_m, t_{m+1})$. Use the r and $\mathfrak{g} = \mathfrak{e}_\mu + \mathfrak{q}_m$ versions of (1.9) and (1.10) to defined the Seiberg-Witten Floer complex, but use the canonical basis as labeled by the gauge equivalence classes of solutions to the r and $\mu$ version of (1.3). Proposition 2.8 supplies a particular generator for each such version of the Seiberg-Witten Floer complex. Use $\mathfrak{c}_I(r)$ to denote this fiducial generator. Note that the degree of $\mathfrak{c}_I(r)$ in $\mathbb{Z}/p\mathbb{Z}$ is independent of r when r is sufficiently large. Indeed, this follows from Lemma 5.4 in [T1]. This understood, normalize the degree so that $\mathfrak{c}_I(\cdot)$ has degree zero in $\mathbb{Z}/p\mathbb{Z}$. Note that $\mathfrak{a}^f(\mathfrak{c}_I(\cdot))$ is a differentiable function on $(r_I, \infty)$.

For each $r \in (r_I, \infty) - \{\rho_j\}$ use $\mathcal{B}(r)$ to denote the set of cycles $\mathfrak{n}$ in degree $-1$ that can be written as $\mathfrak{n} = \delta(\mathfrak{c}_1 + \mathfrak{w})$ where $\mathfrak{w}$ can be written without the generator $\mathfrak{c}_1$. For each such cycle $\mathfrak{n}$, use $\mathfrak{a}^f(\mathfrak{n}, r)$ to denote the maximum of $\mathfrak{a}^f$ on the generators that appear in $\mathfrak{n}$ with non-zero weight. Now set $\mathfrak{a}^f_I(r)$ to denote the infimum of the set $\{\mathfrak{a}^f(\mathfrak{n}, r)\}_{\mathfrak{n} \in \mathcal{B}(r)}$.

With $\mathfrak{a}^f_I(r)$ now defined, replace $\theta$ in the statement of the second bullet of Definition 2.4 with $\mathcal{B}(r)$ and likewise replace $\mathfrak{a}^f_\theta(r)$ with $\mathfrak{a}^f_I(r)$. Use these replacements in Definition 2.4's second bullet to define $\hat{E}(r)$, $\mathfrak{v}(r)$ and $\mathfrak{f}(r)$.

**Proposition 2.9** *Given that the Seiberg-Witten Floer homology is trivial, there exists $r_{I*} \geq r_I$ such that $\mathcal{B}(r) \neq \emptyset$ when $r > r_{I*}$.*

This proposition is proved at the end of Section 5.

**Proposition 2.10**: *The function $\mathfrak{a}^f_I$ extends to $(r_{I*}, \infty)$ as a continuous and piecewise differentiable function. In fact, there is a countable set in $(r_{I*}, \infty)$ and a smooth map, $\mathfrak{c}(\cdot) = (A(\cdot), \psi(\cdot))$, from the complement of this set into $\mathrm{Conn}(E) \times C^\infty(M; \mathbb{S})$ such that for each $r$ in the domain of $\mathfrak{c}(\cdot)$, the configuration $\mathfrak{c}(r)$ is a non-degenerate solution to the $r$ and $\mu$ version of (1.3) with $\mathfrak{a}^f_I(r) = \mathfrak{a}^f(\mathfrak{c}(\cdot))$, $\hat{E} = E(A(r))$ and $\mathfrak{v} = \mathfrak{cs}^f(\mathfrak{c}(r))$. In addition, the generator that is defined by $\mathfrak{c}(r)$ appears with non-zero coefficient in a cycle from $\mathcal{B}(r)$.*

*Proof of Proposition 2.10*: Granted that $\mathcal{B}(r) \neq \emptyset$ when $r$ is large, the proof uses the chain maps from the first three steps of Proposition 2.5 to prove that $\mathfrak{a}^f_I$ has the desired extension. With these chain maps available, the argument is, but for cosmetics, the same as that given in Section 4a of [T1] to prove the latter's Proposition 4.2.

With Proposition 2.10 in hand, use Proposition 1.9 and the argument in Section 4 of [T1] for Proposition 4.6 and Corollary 4.7 of [T1] with $\mathfrak{a}_\theta$ replaced by $\mathfrak{a}^f_I$ to prove the following analog of Proposition 2.7: There exists a sequence $\{r_n, (A_n, \psi_n)\}$ with $\{r_n\}$ increasing and unbounded, with $(A_n, \psi_n)$ satisfying the $r = r_n$ and $\mu$ version of (1.3), with $(A_n, \psi_n)$ in $\mathcal{B}(r_n)$, and with one of the following two properties:

- $\{E(A_n)\}$ is *bounded*.
- $E(A_n) > c\, r_n$ and $\mathfrak{cs}^f(A_n, \psi_n) > c\, r_n^2$ where $c > 0$ is independent of $n$.

(2.3)

The second option in (2.3) is ruled out by Proposition 1.10. Meanwhile, Lemma 5.4 in [T1] guarantees that $\{\sup_M (1 - |\psi_n|)\}$ must be bounded away from zero by an n-independent, positive number when $n$ is large because the $\mathbb{Z}/p\mathbb{Z}$ degree of $(A_n, \psi_n)$ differs by 1 from that of $(A_I, \psi_I)$. Granted all of this, an appeal to Theorerm 1.5 finds the desired set of closed, integral curves of the vector field $v$.

### d) Proof of Proposition 2.8

Consider looking for a solution that has the form $(A, \psi) = (A_I, \psi_I) + ((2r)^{1/2} b, \eta)$ with $(b, \eta) \in C^\infty(M; iT^*M \oplus \mathbb{S})$. The pair $(A, \psi)$ will solve the $r$ and $\mu$ version of (1.3) if $\mathfrak{h} = (b, \eta, \phi) \in C^\infty(iT^*M \oplus \mathbb{S} \oplus i\mathbb{R})$ solves the system of equations

- $*db - d\phi - 2^{-1/2} r^{1/2} (\psi_I^\dagger \tau \eta + \eta^\dagger \tau \psi_I) - r^{1/2} \eta^\dagger \tau \eta = \frac{i}{2} r^{-1/2} (*d\mu + \varpi_K)$,
- $D_{A_I} \eta + 2^{1/2} r^{1/2} (\text{cl}(b)\psi_I + \phi\psi_I) + 2r^{1/2} (\text{cl}(b)\eta + \phi\eta) = 0$,
- $*d*b - 2^{-1/2} r^{1/2} (\eta^\dagger \psi_I - \psi_I^\dagger \eta) = 0$.

(2.4)

Let $\mathfrak{L}_0$ denote that $(A_I, \psi_I)$ the version of the operator that appears in (1.7). Then (2.4) has the schematic form $\mathfrak{L}_0 \mathfrak{h} + r^{1/2} \mathfrak{h} * \mathfrak{h} = r^{-1/2} v$. According to (5.25) in [T1], $\mathfrak{L}_0$ is invertible when $r$ is large, so a solution, $\mathfrak{h}$, to (2.4) can be viewed as a fixed point of the map

$$\mathfrak{h} \to T(\mathfrak{h}) = \mathfrak{L}_0^{-1}(r^{-1/2} v - r^{1/2} \mathfrak{h} * \mathfrak{h}).$$

(2.5)

To see that the map T has a fixed, point, introduce the Hilbert space $\mathbb{H}$ that is obtained by completing $C^\infty(M; iT^*M \oplus \mathbb{S} \oplus i\mathbb{R})$ using the norm $\|\cdot\|_\mathbb{H}$ whose square has value on $\mathfrak{h}$ given by

$$\|\mathfrak{h}\|_\mathbb{H}^2 = \int_M |\nabla_I \mathfrak{h}|^2 + \tfrac{1}{4} r \int_M |\mathfrak{h}|^2 ;$$

(2.6)

Here, $\nabla_I$ is defined so that $\nabla_I(b, \eta, \phi) = (\nabla b, \nabla_{A_I} \eta, d\phi)$. Equation (5.23) in [T1] guarantees that

$$\|\mathfrak{L}_0 \mathfrak{h}\|_2^2 \geq \|\mathfrak{h}\|_\mathbb{H}^2$$

(2.7)

when $r$ is large. Because $|\nabla_I \mathfrak{h}| \geq |d|\mathfrak{h}||$, a Sobolev inequality guarantees that $\|\mathfrak{h}\|_\mathbb{H}$ dominates the $L^p$ norms of $\mathfrak{h}$ for $p \leq 6$. In particular, there exists $c_0 > 0$ that is independent of $\mathfrak{h}$ and $r$ and such that $\|\nabla \mathfrak{h}\|_\mathbb{H} \geq c_0 \|\mathfrak{h}\|_6$. This inequality implies that

$$\|\mathfrak{h}\|_4^4 \leq c_0 \|\nabla \mathfrak{h}\|_\mathbb{H}^3 \|\mathfrak{h}\|_2 \leq c_0 r^{-1/2} \|\mathfrak{h}\|_\mathbb{H}^4.$$

(2.8)

The inequality in (2.8) implies that T extends to give a smooth map from $\mathbb{H}$ to itself.

Granted this fact, fix $R > 2\|v\|_2$ and let $\mathbb{B}_R \subset \mathbb{H}$ denote the ball of radius $r^{-1/2} R$. As is explained next, T maps $B_R$ to itself as a contraction mapping when $r$ is sufficiently large. Indeed, it follows from (2.7) and (2.8) that

$$\|T(\mathfrak{h})\|_{\mathbb{H}} \leq \tfrac{1}{2} r^{-1/2} R + r^{1/2} c_0 \|\mathfrak{h}\|_4^2 \leq \tfrac{1}{2} r^{-1/2} R + r^{1/2} c_0 R^2 r^{-5/4}) \leq \tfrac{1}{2} r^{-1/2} R (1 + 2c_0 R r^{-1/4})$$
(2.9)

when $\mathfrak{h} \in \mathbb{H}$. A similar calculation proves that $\|T(\mathfrak{h}) - T(\mathfrak{h}')\|_{\mathbb{H}} \leq c_0 r^{-1/4} R \|\mathfrak{h} - \mathfrak{h}'\|_{\mathbb{H}}$ when $\mathfrak{h}$ and $\mathfrak{h}'$ are in $\mathbb{H}$ and r is large. These bounds with the contraction mapping theorem imply that any large r version of T has a unique fixed point in $\mathbb{B}_R$. Standard elliptic regularity arguments can be employed to prove that the fixed point is smooth; thus the fixed point is a section of $C^\infty(iT^*M \oplus \mathbb{S} \oplus i\mathbb{R})$ that obeys (2.4).

Consider next the norm of $\psi = \psi_I + \eta$. To this end, write the middle equation in (2.4) as $D_{A_I} \eta = r^{1/2} ([\mathfrak{h}] + [\mathfrak{h} \otimes \mathfrak{h}])$ where $[\cdot]$ is shorthand in each case for some linear endomorphism with an r-independent pointwise norm. Use the Green's function for the operator $D_{A_I}$ to conclude that

$$|\eta| \leq c_0 \|\mathfrak{h}\|_2 + c_0 r^{1/2} \sup_{x \in M} \int_M \frac{1}{\text{dist}(x,\cdot)^2} (|\mathfrak{h}| + |\mathfrak{h}|^2) \ .$$
(2.10)

Now, if u is an $L_1^2$ function on M, then the function $\text{dist}(x,\cdot)^{-1} u$ is square integrable, and its norm is bounded by a constant multiple of the $L_1^2$ norm of u. Granted that such is the case, (2.10) implies the following: There exists an r independent constant $\rho_0$ such that for any $\rho \in (0, \rho_0)$ one has

$$|\eta| \leq c_0 \|\mathfrak{h}\|_2 + c_0 r^{1/2} (\rho^{-1/2} \|\mathfrak{h}\|_2 + \rho^{1/2} \|\mathfrak{h}\|_{\mathbb{H}} + \|\mathfrak{h}\|_{\mathbb{H}}^2) \ .$$
(2.11)

Take $\rho = r^{-1/4}$ in (2.11) and use the fact that $\|\mathfrak{h}\|_{\mathbb{H}} \leq c_0 r^{-1/2}$ and $\|\mathfrak{h}\|_2 \leq c_0 r$ to conclude that $|\eta| \leq c_0 r^{-1/4}$. Thus, $|\psi| \geq 1 - c_0 r^{-1/4}$.

Turn now to the uniqueness assertion that is made by Proposition 2.8. Suppose that $\delta \in (0, \tfrac{1}{2})$ has been given and that $(A, \psi = (\alpha, \beta))$ is a solution to the r and $\mu$ version of (1.3) with the property that $1 - |\psi| < \delta$ at all points in M. This being the case, it follows from the bounds in Lemma 1.6 that $|\alpha| \geq 1 - \delta - c_0 r^{-1/2}$ at all points in M. Given $\varepsilon > 0$, one can now argue as at the end of Section 8 in [T1] that there exists $\delta_0$ such that if $\delta < \delta_0$, then the bounds in (5.24) of [T1] hold. The same scaling arguments used in Section 8 of [T1] prove that

$$|\nabla(\nabla \alpha)| \leq \varepsilon r \quad \text{and} \quad |\nabla'(\nabla' \beta)| \leq \varepsilon r^{1/2} \ .$$
(2.12)

can be assumed as well.

Because $\alpha$ is nowhere zero when $\delta$ is small, there exists $u \in C^\infty(M; S^1)$ such that $e^{iu} \alpha = |\alpha| 1_\mathbb{C}$. Change $(A, \psi)$ to this new gauge, but use $(A, \psi)$ again to denote the resulting pair of connection and spinor. Because $\psi = (\alpha, \beta)$ with $\alpha = |\alpha| 1_{\mathbb{C}}$, the connection A has the form $A = A_I + \hat{a}$ with

$$\hat{a} = \tfrac{1}{2} (\alpha^{-1} \nabla \alpha - \bar{\alpha}^{-1} \nabla \bar{\alpha}) .$$

(2.13)

What with (5.24) in [T1] and (2.12), this implies that

$$r^{-1/2}|\hat{a}| + r^{-1}|\nabla \hat{a}| \leq c_0 \varepsilon .$$

(2.14)

It is now necessary to change to yet another gauge so that the result can be written as $(A_I + 2^{1/2} r^{1/2} b, \psi_I + \eta)$ where b and $\eta$ obey the equation that appears in the third bullet of (2.4). This is done with a gauge transformation that has the form $e^x$ with $x \in C^\infty(M; i\mathbb{R})$. Given any such x, the resulting $(b, \eta)$ has the form

$$b = 2^{-1/2} r^{-1/2} (\hat{a} - dx) \quad and \quad \eta = (e^{ix}|\alpha| - 1)\psi_I + e^{ix}(0, \beta) .$$

(2.15)

The equation in the third bullet in (2.4) is obeyed if x obeys the equation

$$d^\dagger dx + 2r|\alpha|\sin(x) = d^\dagger \hat{a} .$$

(2.16)

As explained momentarily, if $\varepsilon$ is small, and then r is sufficiently large, this equation for x can also be solved using a fixed point strategy and the solution so obtained obeys

$$|x| + r^{-1/2}|dx| \leq c_0 \varepsilon.$$

(2.17)

To see how this works, write $\Delta_r = d^\dagger d + 2r|\alpha|$. This operator is invertible and so a solution to (2.16) obeys $x = \mathcal{T}(x)$, where

$$\mathcal{T}(x) = \Delta_r^{-1}[d^\dagger \hat{a} + 2r |\alpha|(x - \sin(x))] .$$

(2.18)

In order to prove the existence of a fixed point of the right sort, introduce the Banach space $C^0(M; i\mathbb{R})$. Given $\rho \in (0, 1)$, let $\mathcal{B}_\rho$ denote the ball in this space where the norm is bounded by $\rho$. As is argued next, there exist r-independent constants $\rho$ and $\varepsilon_*$ such that if $\varepsilon < \varepsilon_*$ then $\mathcal{T}$ maps $\mathcal{B}_\rho$ to itself as a contraction mapping. This argument uses the following lemma.

**Lemma 2.11**: *There exists $\kappa > 0$ such that if $\varepsilon < \kappa$ and $r \geq 1$, then the following is true: If $g \in C^\infty(M; i\mathbb{R})$ and $y \in C^\infty(M; i\mathbb{R})$ obeys $\Delta_r y = g$, then $r|y| + r^{1/2}|dy| \leq \kappa \sup_M |g|$.*

This lemma is proved momentarily.

Given the bound from the lemma for the supremum norm of the solution to $\Delta_r y = g$, it follows from (2.14) that

$$|\mathcal{T}(x)| \le c_0(\varepsilon + (\sup_M |x|)^2) .$$
(2.19)

Thus, $\mathcal{T}$ maps $\mathcal{B}_\rho$ to itself if $\varepsilon < \frac{1}{2} c_0^{-1} \rho$ and $\rho < \frac{1}{2} c_0^{-1}$. A similar calculation finds that $\mathcal{T}$ is a contraction if $\rho < \rho_0$ with $\rho_0$ independent of r. Granted these bounds, then $\mathcal{T}$ has a unique fixed point in $\mathcal{B}_\rho$ and the fixed point, x, obeys $|x| \le c_0 \varepsilon$. Elliptic regularity then implies that x obeys (2.16). The bound on $|dx|$ asserted by (2.17) follows from (2.14) given what is said about $|dy|$ in Lemma 2.11.

To continue with the proof of the uniqueness assertion in Proposition 2.8, now use $(A, \psi)$ to denote the gauge transformed pair with $A = A_I + 2^{1/2} r^{1/2} b$ and $\psi = \psi_I + \eta$ where the pair $(b, \eta)$ obeys the equation given in the third bullet of (2.4). Write $\mathfrak{b} = (b, \eta, 0)$ and note that by virtue of (2.14) and (2.16), this section of the Hilbert space $\mathbb{H}$ obeys $|\mathfrak{b}| \le c_0 \varepsilon$. It is also the case that $\mathfrak{b}$ is a fixed point of T. Thus, if $\mathfrak{b}$ is in the ball $\mathbb{B}_R$, then it must be the contraction mapping solution. To prove that such is the case, recall that the fixed point equation is equivalent to the equation $\mathfrak{L}_0 \mathfrak{b} = r^{-1/2} \mathfrak{v} + r^{1/2} \mathfrak{b} * \mathfrak{b}$. This, (2.7) and the bound of $c_0 \varepsilon$ on $\sup_M |\mathfrak{b}|$ imply that

$$\|\mathfrak{b}\|_{\mathbb{H}} \le r^{-1/2} \|\mathfrak{v}\|_2 + c_0 r^{1/2} \varepsilon \|\mathfrak{b}\|_2 \le r^{-1/2} \|\mathfrak{v}\|_2 + c_0 \varepsilon \|\mathfrak{b}\|_{\mathbb{H}} .$$
(2.20)

Hence if $\varepsilon < \frac{1}{2} c_0^{-1}$, then $\mathfrak{b}$ is in the radius $r^{-1/2} R$ ball in $\mathbb{H}$ if R is greater than a suitable r-independent lower bound and r is large. Thus, the given solution is gauge equivalent to the contraction mapping solution.

*Proof of Lemma 2.11*: Suppose that $g \in C^0(M; i\mathbb{R})$ is given and that y solves the equation $d^\dagger dy + 2r|\alpha| y = g$. If $|\alpha| \ge \frac{1}{2}$, then

$$d^\dagger d |y| + r|y| \le |g| ,$$
(2.21)

and so the maximum principle implies that $|y| \le r^{-1} \sup_M |g|$. To obtain the bound for dy, differentiate the equation $\Delta_r y = g$ to obtain

$$\nabla^\dagger \nabla dy + \text{Ric}(dy) + 2r|\alpha| dy = dg - 4r d|\alpha| y ,$$
(2.22)

where Ric(·) is defined using the metric's Ricci tensor. When $p \in M$, use $G_p(\cdot)$ to denote the Green's function for the operator $\nabla^\dagger \nabla + \text{Ric}(\cdot) + 2r$ on $C^0(M; iT^*M)$ with pole at p. A maximum principle argument just like that imployed for (2.21) establishes the following: When $w \in C^0(M; iT^*M)$, then $|G_p(w)| \le c_0 r^{-1} \sup_M |w|$. This then implies that

$$\sup_M |dy| \le c_0 (\sup_M |1 - |\alpha||) \sup_M |dy| + c_0 \sup_M |d\alpha| \sup_M |y| + \sup_{p \in M} |G_p(dg)| .$$
(2.23)

Thus,

$$(1 - c_0 \varepsilon) \sup_M |dy| \le c_0 \varepsilon r^{-1/2} |g| + \sup_{p \in M} |G_p(dg)| .$$
(2.24)

To bound the last term in (2.24) in terms of the supremum norm of g, it is necessary to integrate by parts and then obtain a bound for the norms of the derivatives of the section $x \to G_p(x)$. This can be done readily using standard parametrix techniques and results in the bound

$$|\nabla G_p(x)| \le c_0 \frac{1}{\text{dist}(p,x)^2} e^{-\sqrt{r} \, \text{dist}(p,x)} .$$
(2.25)

Thus, an integration by parts shows that $G_p(dg)$ obeys

$$|G_p(dg)| \le c_0 \int_M \frac{1}{\text{dist}(p,\cdot)^2} e^{-\sqrt{r} \, \text{dist}(p,\cdot)} |g| \le c_0 r^{-1/2} \sup_M |g| .$$
(2.26)

and so $|G_p(dg)| \le c_0 r^{-1/2} \sup_M |g|$. This bound with (2.24) completes the proof of the lemma.

### 3. Proof of Proposition 1.9

Before starting, it is worth noting that what is asserted in Proposition 1.9 holds whether or not the 1-form a that appears in (1.3) is a contact form. It is only required that a be smooth and have norm 1. To elaborate, note that the proof that follows assumes that $(A, \psi)$ obeys (1.3) and that the conclusions of Lemmas 1.6-1.7 hold. These three lemmas are proved in Sections 6a-6c of [T1], and their proofs make no use of the assumption that $da = 2*a$ or that $a \wedge da > 0$.

Assume in what follows that $r \ge 2$ and $E \ge 2$. Let $(A_E, \psi_E)$ denote the pair in $\text{Conn}(E) \times C^\infty(M; \mathbb{S})$ that is used for the definition of $f$. Note that the curvature of $A_E$ is a harmonic form whose $L^2$ norm, and thus $C^{k \ge 0}$ norms are bounded apriori. An appeal to Lemma 1.8 finds a smooth map $u: M \to S^1$ such that $A - u^{-1}du = A_E + \hat{a}$ where $|\hat{a}| \le c_0 r^{2/3} (1 + |E|)^{1/3}$. Use A now to denote $A_E + \hat{a}$ and $\psi$ to denote the corresponding gauge transformation of what was originally called $\psi$.

The bound just given for $\hat{a}$ implies that

$$|\mathfrak{cs}(A)| \leq c_0 r^{2/3} (1 + |E|)^{4/3} ,$$

(3.1)

as can be seen by repeating the argument for (4.2) in [T1]. Granted (3.1), the assertion made by Proposition 1.9 follows with a suitable bound on $|f(A, \psi)|$. Note in this regard that the bound given by Proposition 1.10 does not serve for the purposes at hand.

a) **The bound for $f$**

A bound for the spectral flow is obtained in two steps. The first step bounds the absolute value of the spectral flow for the 1-parameter family of operators $\{\mathfrak{L}_s\}_{s \in [0,1]}$ where $\mathfrak{L}_s$ is given the version of (1.7) that uses $(A, \psi)$ but has $r$ replaced by $s^2 r$. The second step bounds the spectral flow for the family of Dirac operators $\{D_{A^s}\}_{s \in [0,1]}$ where $A^s = A + s\hat{a}$.

As explained next, both steps employ the strategy that is outlined in Section 5b of [T1], and in [T2]. To begin, suppose $\{\mathcal{L}_s = \mathcal{L} + \mathfrak{q}_s\}_{s \in [0,1]}$ is as described in Section 5a in [T1]. Take the diffeomorphism $\Phi$ from Section 5b of [T1] to be the identity map from $\mathbb{R}$ to itself, and fix $T > 0$. Let $\mathfrak{n}_s$ denote the number of linearly independent eigenvectors of $\mathcal{L}_s$ whose eigenvalue has absolute value no greater than $T$, and let $\mathfrak{n} = \sup\{\mathfrak{n}_s\}_{s \in [0,1]}$. Introduce the function $\wp(s)$ as in (5.6) of [T1]. Then the spectral flow for the family $\{\mathcal{L}_s\}_{s \in [0,1]}$ has absolute value no greater than

$$\tfrac{1}{2T} \mathfrak{n} \sup \{\|\tfrac{d}{ds} \mathfrak{q}_s\|_{op}\}_{s \in [0,1]} ,$$

(3.2)

where the norm $\|\cdot\|_{op}$ here denotes the operator norm.

In the case of $\{\mathfrak{L}_s\}$, the supremum in (3.2) is bounded by $c_0 r^{1/2}$. For the family of Dirac operators, the analogous norm is bounded by $c_0 |\hat{a}|$, and thus by $c_0 r^{2/3}(1 + |E|)^{1/3}$. In both cases, $T$ will be taken to be $(1 + |E|)^{1/2}$ with $E = E(A)$. This understood, the spectral flow in the case of $\{\mathfrak{L}_s\}$ is bounded in absolute value by

$$c_0 r^{1/2} (1 + |E|)^{-1/2} \mathfrak{n} .$$

(3.3)

and in the second case, by

$$c_0 r^{2/3} (1 + |E|)^{-1/6} \mathfrak{n} .$$

(3.4)

Granted the preceding, the proof Proposition 1.9 requires only a suitable bound for $\mathfrak{n}$ in the two cases. Note in this regard that the bound given by Proposition 5.2 of [T1] does not suffice. The following proposition supplies the desired bound.

**Proposition 3.1**: *There exists a constant, $\kappa > 1$, with the following significance: Suppose that $r \geq \kappa$, that $\mu \in \Omega$ has $C^3$ norm bounded by 1, and that $(A, \psi)$ is solution to the $r$ and* µ *version of (1.3). For* $s \in [0, 1]$, *define* $\mathcal{L}_s$ *to be either the version of (1.7) that is defined using* (A, ψ) *but with* r *replaced by* $s^2 r$, *or the Dirac operator* $D_{A_E + s(A-A_E)}$. *Set* $E = E(A)$ *and let* $\mathfrak{n}_s$ *denote the number of linearly independent eigenvectors of* $\mathcal{L}_s$ *whose eigenvalue has absolute value less than* $(1 + |E|)^{1/2}$. *Then* $\mathfrak{n} \leq \kappa(1 + |E|)^{3/2} (\ln r)^\kappa$.

Proposition 1.9 follows from Proposition 3.1 with (3.1), (3.3) and (3.4).

The strategy that is used to bound $\mathfrak{n}$ is described next and the remaining subsections supply the details. To begin the story on $\mathfrak{n}$, let $\mathbb{H}$ denote the domain Hilbert space of the family $\{\mathcal{L}_s\}$. Thus $\mathbb{H} = L^2_1(iT^*M \oplus \mathbb{S} \oplus i\mathbb{R})$ when $\{\mathcal{L}_s = \mathfrak{L}_s\}_{s \in [0,1]}$ where $\mathfrak{L}_s$ is the version of (1.7) that has r replaced by $s^2 r$. Meanwhile, $\mathbb{H} = L^2_1(\mathbb{S})$ when $\{\mathcal{L}_s = D_{A_E + s(A-A_E)}\}_{s \in [0,1]}$. In each case, the bound on $\mathfrak{n}$ is obtained by exhibiting a set, $\Theta$, of points in M with the following properties:

- *The set* $\Theta$ *has at most* $\kappa(1 + |E|)^{3/2} (\ln r)^\kappa$ *points.*
- *Let* $\vartheta_E$ *denote the span of the set of eigenvectors of* $\mathcal{L}_s$ *whose eigenvalue has absolute value* $1 + |E|$ *or less. If* $\mathfrak{j} \in \vartheta_E$ *and vanishes on* $\Theta$, *then* $\mathfrak{j}$ *is identically zero.*

(3.5)

The proof that such a set $\Theta$ exists exploits certain pointwise and $L^2$ bounds on the covariant derivatives of the elements in $\vartheta_E$. The desired $L^2$ bound uses little more than the Bochner-Weitzenboch formula for $\mathcal{L}_s^2$. The pointwise bounds require more work since the maximum allowed size of the covariant derivative of an element in $\vartheta_E$ at any given point in M is mostly determined by the size of $|B_A|$ in a surrounding ball. In particular, a region in M where $|B_A|$ is relatively large must contain more of $\Theta$'s points than a corresponding region where $|B_A|$ is small. This unavoidable complication is accomodated by decomposing M into cylindrical regions where $r(1 - |\alpha|^2)$ is significantly greater than $1 + |E|$, and the complementary part where this function is roughly $1 + |E|$ or less.

Such a large/small curvature decomposition is facilitated by the introduction of the notion of an adapted coordinate chart map from $C \times [-\delta, \delta]$ into M. Here, $C \subset \mathbb{C}$ is the disk of radius δ centered at the origin. Given a positive number, R, that is less than $\tfrac{1}{2}\delta$, use $\Delta_R \subset C$ to denote the disk with center at the origin and radius R. In what follows, R is taken to be the maximum of $100 r^{-1/2}$ and $\upsilon(1 + |E|)^{-1/2}$ with $\upsilon \in (0, 1)$ specified below. The following lemma supplies a cover of M by the images of adapted coordinate chart maps. Its proof is straightforward, and so omitted.

**Lemma 3.2**: *There exists a constant,* $\kappa > 1$ *with the following significance: Given* $R \in (0, \tfrac{1}{2}\delta)$, *there is a set,* Φ, *of adapted coordinate chart maps from* $C \times [-\delta, \delta]$ *into* M *with the following properties*:
- *There are at most* $\kappa R^{-3}$ *elements in* Φ.

- $\cup_{\phi \in \Phi} \phi(\Delta_R \times [-\frac{1}{4} R, \frac{1}{4} R])$ *covers* M.
- *No point in* M *is contained in more than* $\kappa$ *elements of the set* $\{\phi(\Delta_R \times [-R, R])\}_{\phi \in \Phi}$.

The set $\Theta$ has a decomposition as $\oplus_{\phi \in \Phi} \Theta_\phi$, where the points in $\Theta_\phi$ reside in the $\phi$-image of $\Delta_R \times [-R, R]$. A bound on the number of points in $\Theta_\phi$ is determined by

$$E_\phi = r \int_{\phi(\Delta_R \times [-R,R])} |1 - |\alpha|^2| \ .$$

(3.6)

In particular, the number of points in $\Theta_\phi$ is no greater than $\kappa(1 + R^{-1} E_\phi)(\ln r)^\kappa$ where $\kappa > 1$ is independent of r, (A, $\psi$), $\phi$, and R. Since $\Phi$ has at most $\kappa R^{-3}$ elements, it then follows that $\Theta$ has fewer than $\kappa R^{-3} (\ln r)^\kappa (1 + R^2 \sum_{\phi \in \Phi} E_\phi)$ points. What with the second bullet in Lemma 3.2 and Lemma 1.6, there is a constant, $c_0$, such that

$$\sum_{\phi \in \Theta} E_\phi \leq c_0 (1 + |E|) \ .$$

(3.7)

Thus, the set $\Theta$ has fewer than $\kappa R^{-3}(1 + R^2|E|))(\ln r)^\kappa$ points. Since $R^{-1} \leq \upsilon^{-1}(1 + |E|)^{1/2}$, this gives the asserted bound on $\mathfrak{n}$ given that $\upsilon$ is chosen to be independent of r, E and (A, $\psi$).

**b) Integral bounds for norms of elements in $\vartheta_E$.**

The desired bounds exploit the formulas for $\mathcal{L}_s$ and $\mathcal{L}_s^2$ in a number of ways. Such formula are first used to derive an apriori bound on the $L^2$ norm of the covariant derivate along v of any given $j \in \vartheta_E$. The following lemma states what is needed.

**Lemma 3.3**: *There exists* $\kappa > 1$ *with the following significance: Fix* $r \geq 1$, *a form* $\mu \in \Omega$ *with* $C^3$ *norm bounded by* 1, *and a solution,* (A, $\psi$), *to the* r *and* $\mu$ *version of (1.3). Set* E $= E(A)$. *Fix* $s \in [0, 1]$ *and let* $\mathcal{L}_s$ *denote either the version of the operator in (1.7) with* r *replaced by* $s^2 r$, *or the Dirac operator* $D_{A_E + s(A-A_E)}$. *Let* j *denote a linear combination of those eigenvectors of* $\mathcal{L}_s$ *whose eigenvalue has absolute value less than* $E^{1/2}$. *Then*

$$\|\nabla_v j\|_2 \leq \kappa (1 + |E|)^{1/2} \|j\|_2 \ ;$$

*here* $\nabla_v$ *denotes the covariant derivative in the direction* v *as defined using the connection* A *in the case where* $\mathcal{L}_s$ *is given by the* $r \to s^2 r$ *version of (1.7), and as defined using the connection* $A_E + s(A - A_E)$ *in the case where* $\mathcal{L}_s = D_{A_E + s(A-A_E)}$.

*Proof of Lemma 3.3*: The simpler case is that where $\mathcal{L}_s$ is the Dirac operator and so this case is considered first. To this end, suppose first that j is as described by the lemma. No generality is lost by assuming that $\|j\|_2 = 1$. Use $q$ to denote $\mathcal{L}_s^2 j$. The Bochner-Weitzenboch formula for the Dirac operator reads

$$\nabla^\dagger \nabla j + s\,\mathrm{cl}(B_A) j + \tfrac{1}{4} \mathfrak{R} j = q \tag{3.8}$$

Here $\mathfrak{R}$ denotes an endomorphism of $\mathbb{S}$ that is independent of $r$ and $(A, \psi)$; in particular, it is uniformly bounded. Meanwhile, $\nabla$ denotes the covariant derivative on $\mathbb{S}$ that is defined by the connection $A_E + s(A - A_E)$ and the canonical connection on $K$.

Write $j = (j_0, j_1)$ and $q = (q_0, q_1)$ using the decomposition of $\mathbb{S}$ as $E \oplus EK^{-1}$. Take the $L^2$ inner product of (3.8) with $j_1$ and, what with Lemma 1.6, the resulting equation implies the inequality

$$\|\nabla' j_1\|_2^2 + sr\langle(1 - |\alpha|^2), |j_1|^2\rangle_2 \leq c_0(\langle j_1, q_1\rangle_2 + 1) \tag{3.9}$$

Here, $\langle\,,\,\rangle_2$ denotes the $L^2$ inner product. Since $\|q_1\|_2 \leq \|q\|_2 \leq (1 + |E|)$, this then finds $\|\nabla' j_1\|_2 \leq c_0 (1 + |E|)^{1/2}$. This last inequality implies the desired bound for the $L^2$ norm of the directional derivative of $j_1$. To obtain the desired bound for $\nabla_v j_0$, introduce $q = (q_0, q_1)$ to denote $\mathcal{L}_s j$. The $E$ component of the equation $\mathcal{L}_s j = q$ equates $\nabla_v j_0$ with $q_0$, linear combinations of the directional derivatives of $j_1$ along the kernel of $a$, and multiples of $j_1$ by functions that depend only on the metric. This understood, it follows that $\|\nabla_v j_0\|_2 \leq c_0(\|q_0\|_2 + (1 + |E|)^{1/2}) \leq c_0 (1 + |E|)^{1/2}$ since the $L^2$ norm of $q_0$ is no greater than that of $q$, and the latter is no greater than $(1 + |E|)^{1/2}$.

Consider next the case when $\mathcal{L}_s$ is the operator $r \to s^2 r$ version of (1.7). As in (1.7), write the components of $j$ as $(b, \eta, \phi)$; and assume again that $\|j\|_2 = 1$. Use $q$ to denote $\mathcal{L}_s^2 j$. The corresponding three components of $q$ are

- $\nabla^\dagger \nabla b + \mathrm{Ric}(b) + 2s^2 r |\psi|^2 b + \sqrt{2}\,sr^{1/2}(\eta^\dagger \nabla \psi - (\nabla \psi)^\dagger \eta)$
- $\nabla_A^\dagger \nabla_A \eta + \mathrm{cl}(B_A)\eta + \tfrac{1}{4} R\eta - s^2 r[(\eta^\dagger \psi - \psi^\dagger \eta)\psi + (\psi^\dagger \tau^k \eta + \eta^\dagger \tau^k \psi)\tau^k \psi] - 2\sqrt{2}\,sr^{1/2} b \cdot \nabla \psi$,
- $d^\dagger d\phi + 2s^2 r |\psi|^2 \phi$.

$$\tag{3.10}$$

The last equation in (3.10) implies, that

$$\|d\phi\|_2^2 + s^2 r \|\phi\psi\|_2^2 \leq \|q\|_2 \|\phi\|_2, \tag{3.11}$$

and so $\|d\phi\|_2 \leq (1 + |E|)^{1/2}$. To continue, introduce $b_v$ to denote the contraction of $b$ with $v$. In addition, write $\eta = (\eta_0, \eta_1)$ as per the splitting of $\mathbb{S}$. Take the $L^2$ inner product of the top equation in (3.10) with $b_v a$, and take the $L^2$ inner product of the middle equation in (3.10) with $\eta_1$. With the aid of Lemma 1.6 and some integration by parts, the resulting expressions imply that

$$\|db_v\|_2^2 + \|\nabla'\eta_1\|_2^2 + s^2 r(\|b_v\psi\|_2^2 + \|\eta_1|\psi|\|_2^2) \le c_0(1 + \|j\|_2\|q\|_2) \,.$$

(3.12)

Since $\|q\|_2 \le (1 + |E|)$, this gives the desired bounds on the directional derivatives of $b_v$ and $\eta_1$.

The desired bounds on the directional derivative of $\eta_0$ and of $b - b_v a$ along v are obtained by noting that the $L^2$ norm of $\mathcal{L}_s j$ is bounded by $(1 + |E|)^{1/2}$, and noting the following: First, the component along the kernel of a of the top line in the $r \to s^2 r$ version of (1.3) contains $\nabla_v(b - b_v a)$ and no other derivatives of $b - b_v a$. Moreover, the remaining terms consist of derivatives of $b_v$ and $\phi$, and terms with norm bounded by $sr^{1/2}(|\beta||\eta_0| + |\psi||\eta_1|)$. In particular, Lemma 2.2 and (3.11) and (3.12) imply that the $L^2$ norms of these added terms are no greater than $c_0(1 + |E|)^{1/2}$. Meanwhile, the component in E of the second line in the $r \to s^2 r$ version of (1.7) consists of $\nabla_v \eta_0$ and no other derivatives of $\eta_0$. The remaining terms are combinations of derivatives of $\eta_1$ and terms that are, in any event, no greater than by $sr^{1/2}(|b_v| + |\phi|)|\psi|$. By virtue of (3.11) and (3.12), these latter terms have $L^2$ norms bounded by $c_0(1 + |E|)^{1/2}$.

Let j be as described in this last lemma. The next lemma gives an upper bound for the $L^2$ norm of j over certain balls in M when j vanishes at the ball's center.

**Lemma 3.4**: *There exists $\kappa \ge 1$ with the following significance: Make the same assumptions as in Lemma 3.3 and let j be as described in this same lemma. Let $U \subset M$ denote a ball of radius $\rho \le (1 + |E|)^{-1/2}$ and suppose that $r(1 - |\alpha|^2) \le \rho^{-2}$ on U. Suppose, in addition, that $j = 0$ at the center point of U. Let $\sigma \in (0, \frac{1}{4})$ and let $U_\sigma$ denote the ball with the same center as U and with radius $\sigma\rho$. Let $q = \mathcal{L}_s^2 j$. Then*

$$\int_{U_\sigma} |j|^2 \le \kappa\sigma^4 \left(\int_U |j|^2 + \rho^4 \int_U |q|^2\right)$$

*Proof of Lemma 3.4*: As before, the simplest case is that where $\mathcal{L}_s = D_{A_E + s(A-A_E)}$. To argue for this case, note first that the curvature, B, of the connection $A_s = A_E + s(A - A_E)$ is bounded by $c_0 \rho^{-2}$. As a consequence, there is a trivialization of $E|_U$ so that the connection $A_s$ appears as $A_I + \hat{a}$, where $d*\hat{a} = 0$ and $|\hat{a}| \le c_0 \rho^{-1}$ in the ball, U′, with radius $\frac{3}{4}\rho$ and center that of U. Use Gaussian coordinates to trivialize $TM|_U$. The trivialization of $E|_U$ and $TM|_U$ identifies $\mathbb{S}$ with $U \times \mathbb{C}^2$ and turns j into a $\mathbb{C}^2$-valued function on U. As such, j obeys an equation in U′ that has the form

$$d*dj = \Gamma_0 \cdot j + \Gamma_1 \cdot dj + q \,,$$

(3.13)

where $|\Gamma_0| + |d^*\Gamma_1| \leq c_0 \rho^{-2}$ and $|\Gamma_1| \leq c_0 \rho^{-1}$. For $x \in U'$, use $G_x(\cdot)$ to denote the Greens function of $d^*d$ with pole at x. Then j can be written on $U_\sigma$ as

$$j|_x = \int_U G_x \chi_* (\Gamma_0 \cdot j + \Gamma_1 \cdot dj + q) + \rho^{-2} \int_{U-U'} G_x \chi_{**} j \ ,$$

(3.14)

where $\chi_*$ and $\chi_{**}$ are smooth functions with sup-norms bounded by $c_0$ and with compact support in U. With x viewed as the variable, the Green's function $G_x$ is such that $|G_{(\cdot)}(y)| \leq c_0 |(\cdot) - y|^{-1}$ and $|dG_{(\cdot)}(y)| \leq c_0 |(\cdot) - y|^{-2}$. As is explained momentarily, these bounds and the fact that $j|_0 = 0$ can be used with (3.13) and (3.14) to prove that

$$|j|_x| \leq c_0 |x|^{1/2} (\rho^{-2} \|j\|_{2;U} + \|q\|_{2;U}) \ ,$$

(3.15)

at points $x \in U_{1/4}$. Here, $\|\cdot\|_{2;U}$ is used to denote the $L^2$ norm on U of the indicated function. This last inequality implies what is asserted by the lemma.

To see how (3.15) arises, let g denote an $L^2$ function on U. Break the integral

$$\int_U |G_x - G_0| g$$

(3.16)

into two parts, the first where the distance to the origin is greater than $4|x|$, and the second where the distance to the origin is no greater than $4|x|$. Where $y \in U$ obeys $|y| \geq 4|x|$, use the bound $|G_x(y) - G_0(y)| \leq c_0 |x|/|y|^2$. In the region where $|y| \leq 4|x|$, use instead the bound $|G_x(y) - G_0(y)| \leq c_0 (|x - y|^{-1} + |y|^{-1})$. These bounds with the standard inequality $|\int_U mg| \leq \|m\|_{2;U} \|g\|_{2;U}$ imply a bound for (3.16) of the form $c_0 |x|^{1/2} \|g\|_{2;U}$. Note that (3.13) is used to eliminate a term on the hand side of (3.15) that contains the $L^2$ norm of dj. Such a term can be replaced by the two terms present by using (3.13) to bound that the $L^2$ norm of dj on U' by $c_0 (\rho^{-1} \|j\|_{2;U} + \rho \|q\|_{2;U})$. To obtain this bound, contract both sides of (3.13) by j, multiply the result by a suitably chosen bump function with support on U, and then integrate over U. An integration by parts will lead to an expression that gives the asserted bound for the $L^2$ norm of dj. Equation (3.15) implies the lemma's assertion.

Now consider the case where $\mathcal{L}_s$ is the $r \to s^2 r$ version of (1.7). With $\mathbb{S}$ and TM trivialized over U as before, then $j = (b, \eta, \phi)$ can be viewed as an function on U with values in $i\mathbb{R}^3 \times \mathbb{C} \times i\mathbb{R}$. What with (3.10) giving $\mathcal{L}_s^2$, this function obeys an equation of the form

$$d^*dj + 2s^2 rj = \Gamma_0 \cdot j + \Gamma_1 \cdot dj + sr^{1/2} \Gamma_3(\nabla \psi) \cdot j + q,$$

(3.17)

where $|\Gamma_3(\nabla \psi)| \leq c_0 |\nabla \psi|$ at each point. Here again, $|\Gamma_0| \leq c_0 \rho^{-2}$ and $|\Gamma_1| \leq c_0 \rho^{-1}$.

To procede from here, note first that there is a bound for $|\nabla \psi|$ in U' of the form

$$|\nabla\psi| \le c_0 \rho^{-1} .$$

(3.18)

This bound is justified momentarily. Granted (3.18), the second point to note is that the Green's function to use with (3.17) is that for the operator $d^*d + 2s^2r$ Let $G_x$ now denote the Green's function for this operator with pole at x. This new version of $G_x$ obeys

$$|G_{(\cdot)}(y)| \le c_0 \frac{1}{|(\cdot) - y|} e^{-s(2r)^{1/2}|(\cdot) - y|} \quad \text{and} \quad |dG_{(\cdot)}(y)| \le c_0 \frac{(1 + sr^{1/2}|(\cdot) - y|)}{|(\cdot) - y|^2} e^{-s(2r)^{1/2}|(\cdot) - y|} .$$

(3.19)

But for some straightforward modifications, the arguments for (3.15) can be repeated using (3.18) and (3.19) to prove that

$$|j|_x| \le c_0 \min((s^2r)^{-1/4}, |x|^{1/2}) (\rho^{-2} \|j\|_{2;U} + \|q\|_{2;U}) .$$

(3.20)

when $x \in U_\sigma$. This last bound implies the inequality that is stated by the lemma.

To tie up a loose end, consider now (3.18). A bound for $|\nabla\psi|$ on U´ is obtained by differentiating the equation $D_A\psi = 0$ to obtain an equation for $\nabla\psi$ that has the schematic form $D_A\nabla\psi + \mathcal{R}\psi + B_A*\psi = 0$. Here, the components of $\mathcal{R}$ come from the Riemannian metric, and $B_A*$ is a homomorphism whose components come from $B_A$. To use this equation to bound $|\nabla\psi|$ at a given point $p \in U´$, it proves useful to change the trivialization of $E|_U$. To this end, trivialize E now by parallel transport via A along the radial geodesics from p. At the same time take a Gaussian coordinate chart centered at p and use these coordinates to trivialize $TM|_U$. These trivializations together induce a trivialization of $\mathbb{S}|_U$ and so identify its sections as $\mathbb{C}^2$ valued functions. Using this trivialization, the connection A appears as $A_I + a^p$, where $A_I$ is the product connection for the trivialization, and where $a^p$ obeys $|a^p| \le c_0\rho^{-2} \text{dist}(p, \cdot)$. The Christoffel symbols for the metric also vanish at p.

With the preceding understood, introduce next $G_{I,p}$ to denote the Green's function for the operator $D_{A_I}$ on $\mathbb{S} \otimes T^*M$ with pole at p. With the trivialization of $\mathbb{S}$ over U given, $G_{I,p}(x)$ at $x \in U$ appears as a matrix valued function on U. As such, it obeys $|G_{I,p}(x)| \le c_0 \text{dist}(p, x)^{-2}$ and $|dG_p|_x| \le c_0 \text{dist}(p, x)^{-3}$. Use $G_{I,p}$ to write $\nabla\psi|_p$ as

$$\nabla\psi|_p = -\int_U G_{I,p} \chi^U \text{cl}(a^p) \nabla\psi - \int_U G_{I,p} \chi^U (\mathcal{R}\psi + B_A * \psi) + \int_U G_{I,p} \text{cl}(d\chi^U) \nabla\psi ,$$

(3.21)

where $\chi^U$ is a smooth function with compact support on U that equals 1 where the distance to U´ is less than $\rho/8$. Moreover, $|d\chi^U| \le 32\rho^{-1}$ and $|\nabla d\chi^U| \le 1056\rho^{-2}$.

Granted these bounds and the bound $|\psi| \leq 1$, integration by parts finds the absolute value of the right most term in (3.21) bounded by $c_0 \rho^{-1}$. With $|B_A| \leq c_0 \rho^{-2}$, the bound on the norm of $G_{I,p}$ directly yield a bound by $c_0 \rho^{-1}$ for the absolute value of the middle integral on the right hand side of (3.21). As for the left most integral in (3.21), the bounds on $|a^p|$ and $|G_{I,p}|$ imply that their product has $L^2$ norm on U that is bounded by $\rho^{-3/2}$. As a consequence, the absolute value of the left most integral in (3.21) is no greater than

$$c_0 \, \rho^{-3/2} \, \| \chi^U \nabla \psi \|_2 \, .$$

(3.22)

A bound by $c_0 \rho^{-1}$ on the left most integral on the right side of (3.21) follows from (3.22) given that $\| \chi^U \nabla \psi \|_2 \leq c_0 \rho^{1/2}$. This last bound is obtained by using the Weitzenboch formula in (3.8) with $j = \psi$ and so $q = 0$. To elaborate, take the $L^2$ inner product of both sides of the $j = \psi$ and $q = 0$ version of (3.8) with $(\chi^U)^2 \psi$, and then integrate by parts. The result directly implies the desired norm bound given the bounds for $|B_A|$ and $|d\chi^U|$.

Lemma 3.4 is used to prove bounds for the $L^2$ norm of $j$ over cylinders.

**Lemma 3.5**: *There exists $\kappa \geq 1$ with the following significance: Make the same assumptions as in Lemma 3.3 and let $j$ be as described in that lemma. Fix R as described just prior to Lemma 3.2 and let $\phi: \Delta_R \times [-R, R] \to M$ denote an adapted coordinate chart map. Use $\phi$ to implicitly identify its domain and image. Fix $\rho \in (\frac{1}{2} r^{-1/2}, \frac{1}{4} R)$ to insure that $r(1 - |\alpha|^2) \leq \rho^{-2}$ on $\Delta_\rho \times [-\rho, \rho]$, and suppose that $j = 0$ at $(0, 0)$. Set $q = \mathcal{L}_s^2 j$. If $\sigma \in (0, \frac{1}{4})$, then*

$$\int_{\Delta_{\sigma\rho} \times [-R+\rho, R-\rho]} |j|^2 \leq \kappa \int_{\Delta_\rho \times [-R,R]} (R^2 | \nabla_v j |^2 \, + \, \sigma^3 |j|^2 \, + R\sigma^3 \rho^3 |q|^2)$$

*Proof of Lemma 3.5*: Let $z \in [-R + \rho, R - \rho]$, and with $\rho'$ either $\sigma\rho$ or $\rho$, let $\mathfrak{f}(z)$ denote the $L^2$ norm of $j$ over the ball in $\Delta_R \times [-R, R]$ with center at $(0, z)$ and radius $\rho'$. Note that

$$|\mathfrak{f}(z) - \mathfrak{f}(z')| \leq c_0 (R\rho' \int_{\Delta_\rho \times [-R,R]} | \nabla_z j |^2 \, )^{1/2} \, .$$

(3.23)

Let U denote the ball of radius $\rho$ in $\Delta_R \times [-R, R]$ with center $(0, 0)$, and use the $\rho' = \sigma\rho$ version of this last inequality with Lemma 3.4 to conclude that

$$\int_{\Delta_{\sigma\rho} \times [-R+\rho, R-\rho]} |j|^2 \leq c_0 \, R^2 \int_{\Delta_\rho \times [-R,R]} | \nabla_z j |^2 + c_0 (\tfrac{R}{\sigma\rho}) \sigma^4 ( \int_U (|j|^2 \, + \rho^4 |q|^2 )\, .$$

(3.24)

Next, use the ρ´ = ρ version of (3.23) to conclude that

$$\tfrac{R}{\rho} \int_U |j|^2 \leq c_0 \int_{\Delta_\rho \times [-R,R]} (R^2 |\nabla_z j|^2 + |j|^2).$$

(3.25)

These last two inequalities imply what is asserted by Lemma 3.5.

Lemmas 3.5 and 3.3 are key inputs to the proof of Proposition 3.1.

### c) Refined bounds for the norms of ψ and ∇ψ

Additional bounds for ψ and its covariant derivative are needed in order to exploit Lemma 3.5. The next lemma supplies the required bounds.

**Lemma 3.6**: *There exist constants $\kappa \geq 1$ and $\varepsilon_* \in (0, \tfrac{1}{4})$ with the following significance: Suppose that $r \geq \kappa$, $\mu \in \Omega$ has $C^3$ norm bounded by 1, and that $(A, \psi)$ is a solution to the r and μ version of (1.3). Let $U \subset M$ denote the subset where $(1 - |\alpha|^2) \geq \varepsilon_*$. Then*

$$((1 - |\alpha|^2) + r^{-1}|\nabla\alpha|^2 + |\nabla\beta|^2)|_{x \in M-U} \leq \kappa(e^{-\sqrt{r}\,\mathrm{dist}(x,\partial U)/\kappa} + r^{-1}).$$

*Proof of Lemma 3.6*: Let $w = 1 - |\alpha|^2$. This, the first step of the proof finds r and $(A, \psi)$ independent constants $c_1, c_2, c_3 > 0$ and $\varepsilon_* \in (0, 1)$ and $\kappa_* \in (0, \tfrac{1}{2})$ such that

$$g = w + c_1 r^{-1}|\nabla\alpha|^2 + c_2|\nabla\beta|^2 - c_3 r^{-1}$$

(3.26)

obeys the equation

$$d^\dagger dg + \kappa_* r g \leq 0$$

(3.27)

on the subset in M where $w = 1 - |\alpha|^2 < \varepsilon_*$. This is done by differentiating (6.1) so as to get an equation of the form $\nabla^\dagger \nabla(\nabla\psi) + \cdots = 0$, where the three dots indicate terms with either one or no derivatives of ψ. Take the inner product of the resulting equation first with $(\nabla\alpha, 0)$ and then with $(0, \nabla\beta)$. What with (6.5) in [T1], arguments that differ only cosmetically from those used in Section 2e of [T3] find constants that guarantee (3.27).

With (3.27) understood, the next step constructs a certain positive function that satisfies the analog of (3.27) with the reversed inequality. To this end, fix a smooth function χ: $[0, \infty) \to [0, 1]$ that is non-increasing, equals 1 on $[0, \tfrac{1}{4}]$, and is equal to zero on $[\tfrac{1}{2}, \infty)$. A constant $c \in (r^{-1/4}, 1)$ is needed next; it is fixed shortly so as to be independent of r and $(A, \psi)$. For now, fix any $c \in (r^{-1/4}, 1)$ and let $\chi_x$ denote the function

$$\chi_x(\cdot) = \chi(c\, r^{1/2} \text{dist}(\cdot, x)).$$

(3.28)

The assignment of a given pair $(x, y)$ to $\chi_x(y)$ defines a smooth function on $M \times M$ when $r \geq c_0$. Such a lower bound for $r$ is assumed in what follows. With $\chi_x$ understood, set $\upsilon_x$ to be the integral of $\chi_x$ and define $\rho \in C^\infty(M)$ by

$$\rho(x) = \upsilon_x^{-1} \int_M \chi_x(\cdot) \text{dist}(\cdot, U)$$

(3.29)

Note that $|d\rho| \leq c_0$ since the distance function is Lipshitz with norm 1. Note also that $|\nabla d\rho| \leq c_0 c\, r^{1/2}$.

Now set $z = e^{-c r^{1/2} \rho}$. The function $z$ obeys

$$d^\dagger dz + c_0 c^2\, r z \geq 0.$$

(3.30)

This understood, let $\kappa_*$ be as in (3.17) and set $c = (c_0^{-1} \kappa_*)^{1/2}$ so that the inequality in (3.30) reads $d^\dagger dz + \kappa_* r z \geq 0$. This understood, use Lemma 1.7 to obtain an $r$ and $(A, \psi)$ independent constant, $\kappa_1$, that bounds $g$ everywhere on $M$. Then

$$d^\dagger d(g - \kappa_1 z) + \kappa_* r(g - \kappa_1 z) \leq 0$$

(3.31)

on $M-U$ and $g \leq 0$ on $\partial U$. The maximum principle asserts that $g \leq \kappa_1 z$ on $M-U$. This implies what is stated by the lemma.

**d) Estimates for adapted coordinate chart maps**

The purpose of this step is to define a set $W_\phi \subset \Delta_R$ for an adapted coordinate chart map $\phi$. The following lemma is supplies some input to the definition. The constant $\varepsilon_*$ that appears in the lemma comes from Lemma 3.6.

**Lemma 3.7**: *There exist constants $\kappa \geq 1$, $\delta_1 \in (0, \delta)$ and $\varepsilon_1 \in (0, \varepsilon_*)$ with the following significance: Suppose that $r > \kappa$, $\mu \in \Omega$ has $C^3$ norm bounded by 1 and $(A, \psi)$ is a solution to the $r$ and $\mu$ version of (1.3). Let $\phi$: $C \times [-\delta, \delta] \to M$ denote an adapted coordinate chart map. Fix $R \in [100 r^{-1/2}, \delta_1]$. Then fix a minimal set of disks in $C$ with the following two properties. First, each disk in this set has center in $\Delta_R$ and radius $r^{-1/2}$. Second, $(1 - |\alpha|^2) \geq \varepsilon_1$ on $V \times \{0\}$ when $V$ is a disk from the set. Denote this set by $\Lambda$.*
- *If $(u, z) \in \Delta_R \times [-R, R]$ and $(1 - |\alpha|^2) \geq \varepsilon_*$, then $u \in \cup_{V \in \Lambda} V$.*
- *Let $N_\phi$ denote the number of elements in $\Lambda$. Then*

$$N_\phi \leq \kappa(R^{-1}E_\phi + R^2)$$

*where $E_\phi$ denotes the integral of $r|1 - |\alpha|^2|$ over $\Delta_R \times [-R, R]$*

*Proof of Lemma 3.7:*. The assertion in the first bullet follows using the uniform bound on $|\nabla'\beta|$ in Lemma 1.7 because the Dirac equation sets $\nabla_z\alpha$ equal to a linear combination of derivatives of $\beta$. To establish the assertion in the second bullet, note that each $V \in \Theta$ contains at least one point x such that $(1 - |\alpha|^2) \geq \varepsilon_1$ at $(0, x)$. Given the uniform bound on $|\nabla_z\alpha|$, this implies that there exists $\delta_2 \in (0, \delta)$ that is independent of r, $(A, \psi)$, the map $\phi$ and x; and is such that $(1 - |\alpha|^2) \geq \frac{1}{2}\varepsilon_1$ at any $(z, x)$ with $|z| \leq \delta_2$. It then follows using Lemma 1.7, now for its stated bound on $|\nabla\alpha|$, that

$$r\int_{V\times[-R,R]} (1 - |\alpha|^2) \geq c_1 R \varepsilon_1 - c_0 R^3.$$

(3.32)

for each $V \subset \Theta$. Here, $c_1 > 0$ and $c_0 > 1$ are independent of r, $(A, \psi)$, R and the map $\phi$. Since $\Theta$ is a minimal set of disks, there exists some $c_2 > 0$ that is independent of r, $(A, \psi)$, R, and the map $\phi$ such that at most $c_2$ disjoint disks from $\Theta$ contain any given point in C. This being the case, it follows from (3.32) and Lemma 1.6 that

$$E_\phi = r\int_{\Delta_R \times [-R,R]} |1 - |\alpha|^2| \geq c_2^{-1} c_1 \varepsilon_1 R N_\phi - c_0 R^3.$$

(3.33)

This last bound implies the assertion made in the second bullet of the lemma.

Let $\kappa$ denote the constant that appears in Lemma 3.6, and set

$$\rho_* = \kappa\, r^{-1/2} (\ln(1 + r/(1 + |E|)) + \ln(4\kappa)).$$

(3.34)

Fix an adapted, coordinate chart map $\phi$ and a set $\Lambda$ as described in Lemma 3.7. For each $V \in \Lambda$, let $W^V \subset C$ denote the disk with the same center as V with radius $4\rho_*$. Let $W^{V-}$ denote the disk with the same center as V with radius $\rho_*$. Granted this notation, set $W = W_\phi = \cup_{V \in \Lambda} W^V$ and $W^- = \cup_{V \in \Lambda} W^{V-}$.

**Lemma 3.8**: *There exists $\kappa > 1$ that is independent of r, $(A, \psi)$, R and $\phi$, and has the following significance: Define W as above. Then*

$$(1 - |\alpha|^2) \leq \tfrac{1}{2}\frac{1+|E|}{|E|+r} \quad \text{and} \quad |\nabla\alpha| \leq \kappa(1 + |E|)$$

on $(\Delta_R - W^-) \times [-R, R]$. *In addition, W is contained in a set of disks in C of radius $r^{-1/2}$ that has $N_{W\phi} = \kappa(R^{-1}E_\phi + R^2)(\ln(1 + r/(1 + |E|)))^2$ members.*

**Proof of Lemma 3.8**: The first assertion follows from Lemma 3.6. The second follows using the second bullet in Lemma 3.7 given that W is the union of $N_\phi$ disks, each of radius $\rho_*$.

The next task for this subsection is to specify a cover for $\Delta_R - W$ by disks with relatively large radii. A cover of the required sort is supplied by the next lemma.

**Lemma 3.9**: *There exists $\kappa > 1$ that is independent of r, $(A, \psi)$, R and $\phi$, and has the following significance: Define W as above. Then there is a cover of $\Delta_R - W$ by a set of disks in C with the following properties: First, each disk from the cover is disjoint from $W^-$ and each disk has radius at least $r^{-1/2}$. Second, the number of disks in the cover is bounded by $\kappa(R^{-1}E_\phi + R^2)\ln(1 + r/(1 + |E|))$. Third, the ratio of the respective radii of any two intersecting disks is bounded by $\kappa$. Finally, no more than $\kappa$ distinct disks from the cover contain any given point.*

**Proof of Lemma 3.9**: The set in question is the union of sets $K_0, \ldots, K_Z$ with a bound for Z given momentarily. The set $K_n$ is defined as follows: Let $\upsilon_0 = 0$ and for positive integer n, set $\upsilon_n = 2((\frac{5}{4})^n - 1)$. Let $L_n$ denote the set of points in C whose distance from W is equal to $\upsilon_n \rho_*$. Then $K_n$ is a maximal set of disks in C with the following properties: First, each disk has radius $(1 + \frac{1}{2}\upsilon_n)\rho_*$ and center on $L_n$. Second, the center of any given disk from $K_n$ has distance at least $\frac{1}{4}(1 + \frac{1}{2}\upsilon_n)\rho_*$ from the center of any other disk from $K_n$. The number of elements in $K_n$ is less than $16\pi N_\phi$ with $N_\phi$ as in Lemma 3.7. This is because the length of $L_n$ is at most $2\pi(2 + \upsilon_n)\rho_* N_\phi$. A straightforward induction argument proves that the union of the disks from the set $K_0 \cup \cdots \cup K_n$ contains the set of points in C with distance $\upsilon_{n+1}\rho_*$ or less from W. This construction guarantees that the number Z can be no greater than $(\ln\frac{5}{4})^{-1}\ln(1 + 8R/\rho_*)$. The existence of a uniform bound on the number of disks that contain any given point uses two observations: First, no point in a disk from any $n \geq 2$ version of $\vartheta_n$ has distance less than $(\frac{1}{2}\upsilon_n - 1)\rho_*$ from W. This implies that disks in $K_n$ and $K_{n-k}$ are disjoint when $k \geq 5$. Second, the center of any given disk in $K_n$ has distance at least $\frac{1}{4}(1 + \frac{1}{2}\upsilon_n)\rho_*$ from the center of any other disk in $K_n$. Granted these facts, apply the n-k version of this last observation for $k \in \{0, \ldots, 5\}$ with the fact that $(1 + \frac{1}{2}\upsilon_{n-k}) \geq (\frac{4}{5})^k(1 + \frac{1}{2}\upsilon_n)$ to obtain the desired bound.

**e) The definition of $\Theta_\phi$**

Let $\phi: C \times [-\delta, \delta]$ denote an adapted coordinate chart map. In order to define the set $\Theta_\phi$, introduce the subsets W and $W^-$ as defined in Subsection 3d.

Now, fix $\sigma \in (0, 1)$ for the moment. Its precise value is determined shortly. With $\sigma$ fixed, the following lemma is used to specify the additional points.

**Lemma 3.10**: *There exists $\kappa > 1$ which is independent of* r, $\mu$, (A, $\psi$), R *and* $\phi$ *and has the following significance: Fix $\sigma \in (0, 1)$, and there exists a set of disks in* C *with the following properties*:
- *The disks in this set cover $\Delta_R$.*
- *The set has at most $c_0 \sigma^{-2}(R^{-1}E_\phi + R^2)(\ln(1+r/(1 + |E|)))^2$ elements.*
- *No more than $\kappa$ of the disks from this set contain any given point in* C.
- *Let $\Delta$ denote any disk from this set. The radius of $\Delta$ can be written as $\sigma\rho$ where $\rho \in [\frac{1}{2} r^{-1/2}, \frac{1}{4} R]$. Moreover, $r(1 - |\alpha|^2) \leq \rho^{-2}$ on the radius $\rho$ ball with center at the center point of $\Delta$.*

*Proof of Lemma 3.10*: This lemma follows directly from Lemmas 3.8 and 3.9. Indeed, it follows from Lemma 3.8 that W has an open cover by a set of disks in C that has at most $c_0\sigma^{-2}(R^{-1}E_\phi + R^2)(\ln(1 + r/(1 + |E|)))^2$ members, and is such that each disk in the set has radius $\frac{1}{2}\sigma r^{-1/2}$. Moreover, the last two points in the lemma are also obeyed by the disks in this set. Indeed, the final point follows because Lemma 1.6 guarantees that $r(1 - |\alpha|^2) \leq 4r$. The remaining disks in the set cover $\Delta_R - W$ and are chosen using Lemma 3.10. In this case, Lemma 3.9 implies that each disk obeys the fourth point of Lemma 3.10.

Let $\mathcal{Q}_\phi$ denote the set of disks that are supplied by Lemma 3.10. The points in $\Theta_\phi$ are the centers of all disks of the form $\Delta \times \{0\}$ where $\Delta \in \mathcal{Q}_\phi$.

**f) The specification of $\sigma$ and $\upsilon$**

The constant $\sigma$ enters above in the definition of each $\Theta_\phi$, and the constant $\upsilon$ enters as $R = \max(\upsilon(1 + |E|)^{-1/2}, 100 r^{-1/2})$ in the case where $r > 2\delta^{-1}$. Lemmas 3.3 and 3.5 are the key inputs that determine $\sigma$ and $\upsilon$. As explained next there are values for $\upsilon$ and $\sigma$ that are independent of r, $\upsilon$ and (A, $\psi$) such that any $j \in \vartheta_E$ that vanishes at all points in $\Theta$ must be zero. To see that such is the case, suppose that $\upsilon$ and $\sigma$ have been chosen and that j is an element in $\vartheta_E$ with $L^2$ norm 1 that vanishes at all points in $\Theta$. Fix $\phi \in \Phi$. By virtue of Lemma 3.10 that each disk in $\mathcal{Q}$ can be written as $\Delta_{\sigma\rho}$ where $\rho \in [\frac{1}{2} r^{-1/2}, \frac{1}{4} R]$ and where the concentric disk, $\Delta_\rho$, of radius $\rho$ satifies the conditions in Lemma 3.5. Since the number of disks from $\mathcal{Q}_\phi$ that contain any given point of $\Delta_R$ is apriori bounded, it follows from Lemma 3.5 that

$$\int_{\Delta_R \times [-\frac{1}{4}R, \frac{1}{4}R]} |j|^2 \leq c_0 \int_{\Delta_R \times [-R, R]} (\sigma^{-2} R^2 |\nabla_v j|^2 + \sigma |j|^2 + \sigma R^4 |q|^2)$$

(3.35)

By virtue of Lemma 3.2, this then implies that

$$1 = \int_M |j|^2 \leq c_0 \sigma^{-2} R^2 \int_M |\nabla_v j|^2 + c_0 \sigma + c_0 \sigma R^4 \int_M |q|^2 .$$
(3.36)

Now, according to Lemma 3.3, the L² norm of $\nabla_v j$ is bounded by $c_0(1 + |E|)^{1/2}$, and that of q is bounded by $(1 + |E|)$. Thus, with $R = \upsilon(1 + |E|)^{-1/2}$, the inequality in (3.36) implies that

$$1 \leq c_0 (\sigma^{-2} \upsilon^2 + \sigma + \sigma \upsilon^4) .$$
(3.37)

This understood, choose $\sigma = \frac{1}{4}(1 + c_0)^{-1}$ and choose $\upsilon = \frac{1}{8}(1 + c_0)^{-3/2}$ to obtain from (3.37) the absurd conclusion that $1 < \frac{3}{4}$. Thus, no such $j \in \vartheta_E$ exists for these choices.

## 4. Overtwisting

This section supplies a proof of Theorem 1.4. The proof relies on a relative version of Corollary C from Mrowka and Rollin [MR]. The suggestion to use this corollary in the proof came from Tom Mrowka. As the reader will see, the proof also owes much to the work of Kronheimer and Mrowka in [KM2].

**a) The proof of Theorem 1.4**

To set the stage, choose the $\text{Spin}_\mathbb{C}$ structure on M where $\mathbb{S} = I_\mathbb{C} \oplus K^{-1}$. Let $r_I$ and $\delta$ be as described in Proposition 2.8. Fix $\mu \in \Omega$ with $C^3$ norm less than 1 as in Proposition 2.1. As noted in Proposition 2.8, for each $r > r_I$, there exists a unique gauge equivalence class of pairs $(A, \psi)$ that solves the r and $\mu$ version of (1.3) and is such that $|\psi| \geq 1 - \delta$. Use $\mathfrak{c}_I(r)$ in what follows to denote both the solution given by Proposition 2.8 to (1.3) and its corresponding gauge equivalence class. If follows from Lemma 5.4 in [T1] that $r_I$ can be chose so that the solution $\mathfrak{c}_I(r)$ has degree zero.

Now let $\{\rho_i\} \subset [1, \infty)$ be as in Proposition 2.1, and for $r \in (1, \infty) - \{\rho_i\}$, use the gauge equivalence classes of solutions to (1.3) to define the generators for the Seiberg-Witten Floer homology. Define the Seiberg-Witten Floer homology for such r in the manner that is described in Section 2, Propositions 2.3 and 2.4. Use $\mathcal{H}^0(r)$ to denote the degree zero Seiberg-Witten Floer *cohomology*. By definition, $\mathcal{H}^0(r)$ is a $\mathbb{Z}$-module whose elements are certain equivalence classes of linear functionals on the vector space of degree zero cycles for the Seiberg-Witten Floer homology. The equivalence relation identifies linear functionals t and t´ when $t´ = t + \upsilon(\delta(\cdot))$ where $\upsilon$ is a linear function on the vector space of degree -1 cycles. The module $\mathcal{H}^0$ consists of those equivalence classes of linear functionals that give zero to boundaries. Such a cocycle is said to be *closed*.

The constructions of Kronheimer and Mrowka in [KM2] can be modified using analysis from Chapter 24 of [KM1] to assign a reasonably canonical element to $\mathcal{H}^0(r)$. What follows outlines how this is done. To start, it is necessary to specify some auxiliary data. Here is the data:

- *A smooth, non-decreasing function* $\mathrm{T}\colon \mathbb{R} \to [0, \infty)$ *that has value* 0 *where* $s \leq 1$ *And equals* s *where* $s \geq 2$.
- *A smooth and suitably generic section,* $\hat{\mu}$, *of* $T^*M$ *over* $\mathbb{R} \times M$ *with compact support where* $2 \leq s \leq 4$ *and small* $C^6$ *norm. An upper bound for this norm is specified below, but require in any case that the norm is less than* 1.
- *A smooth, non-increasing function* $\chi\colon \mathbb{R} \to [0, 1]$ *that equals* 1 *where* $s \leq 0$ *and equals* 0 *where* $s \geq 1$.
- *The pair* $(A_I, \psi_I) \in \mathrm{Conn}(I_\mathbb{C}) \times C^\infty(M; \mathbb{S})$ *that solves (2.2) and has* $\mathrm{E}(A) = 0$.

(4.1)

Of interest in what follows are the smooth maps $s \to \mathfrak{d}(s) = (A(s), \psi(s))$ from $\mathbb{R}$ to $\mathrm{Conn}(I_\mathbb{C}) \times C^\infty(M; \mathbb{S})$ that obey

- $\frac{\partial}{\partial s} A = -B_A + r e^{2T}(\psi^\dagger \tau^k \psi - i a) + i\chi(*d\mu + \varpi_K) + i\hat{\mu}$,
- $\frac{\partial}{\partial s} \psi = -D_A \psi$.
- $\lim_{s \to -\infty} \mathfrak{d}(s)$ *exists and defines a non-degenerate solution to* (1.3) *with degree zero*.
- $\lim_{s \to \infty} \mathfrak{d}(s) = (A_I, \psi_I)$.

(4.2)

With regards to the nature of the limits in (4.2), it is sufficient to assume that these limits are defined with respect to the $C^\infty$ topology in $\mathrm{Conn}(I_\mathbb{C}) \times C^\infty(\mathbb{S}_I)$. Even so, the same moduli space arises if the limits are in $C^0$ or, for that matter, $L^2_1$.

As remarked on briefly below, arguments from Sections 2 and 3 of [KM2] for their Theorem 2.4 and the analysis given in Chapter 24 of [KM1] can be used to prove that there are only finitely many gauge equivalence classes of solutions to (4.2) for a suitably generic choice for $\hat{\mu}$. Moreover, each such equivalence classes has an associate sign, either +1 or -1. Granted that such is the case, let $\mathfrak{c}$ denote a solution to (1.3) with degree zero; and let $\sigma(\mathfrak{c})$ denote the sum of the signs that are associated to the solutions to (4.2) with $s \to -\infty$ limit equal to $u \cdot \mathfrak{c}$ with $u \in C^\infty(M; S^1)$. Set $\sigma(\mathfrak{c}) = 0$ if there are no gauge equivalence classes of solutions to (4.2) with $s \to -\infty$ limit in $\mathfrak{c}$'s gauge equivalence class. The assignment of the integer $\sigma(\mathfrak{c})$ to $\mathfrak{c}$ associates an integer to each generator of the Seiberg-Witten Floer complex with $\mathbb{Z}/p\mathbb{Z}$ degree 0. Use $\mathfrak{t}_r$ in what follows to denote the degree zero Seiberg-Witten Floer *cocycle* whose value on a generator $\mathfrak{c}$ is $\sigma(\mathfrak{c})$.

The following is a generalization of Theorem 2.4 in [KM2]:

**Theorem 4.1**: *Let* μ *be as described in Proposition 2.1, and suppose that* $r \notin \{\rho_i\}$. *Then the cocycle* $\mathfrak{t}_r$ *is closed, and so defines an element in* $\mathcal{H}^0(r)$.

The next result constitutes the fundamental input from [MR].

**Theorem 4.2**: *The cocycle* $\mathfrak{t}_r$ *defines the trivial class in* $\mathcal{H}^0(r)$ *when the kernel of* a *is overtwisted*.

*Proof of Theorems 4.1 and 4.2*: Consider first what is said leading up to Theorem 4.1. Kronheimer and Mrowka in Section 3 of [KM2] introduce the notion of an 'asymptotically flat almost Kahler' structure. The equations in (4.2) are the Seiberg-Witten equations on the 4-manifold $\mathbb{R} \times M$ where the geometry is such that the $s \to \infty$ end has such a structure. To elaborate, let $\mathfrak{m}$ denote the Riemannian metric on M. Then the asymptotically flat almost Kahler structure is defined using the metric $e^{2T}(ds^2 + \mathfrak{m})$ on $\mathbb{R} \times M$. The required symplectic form on the $s > 1$ part of $\mathbb{R} \times M$ is $\omega = e^{2T}(ds \wedge a + *a)$. Here, $*$ is the Hodge star on M that is defined by the metric $\mathfrak{m}$. Meanwhile, the end where $s \to -\infty$ has the sort of cylindrical end structure that is considered in Chapter 24 of [KM1].

      The techniques used in [KM2] can be applied to analyze the behavior of solutions to (4.2) on the $s \to \infty$ end of $\mathbb{R} \times M$, and those from Chapter 24 of [KM1] can be used to analyze the behavior of the solutions on the $s \to -\infty$ end. This limit business is discussed in Section 4b and it follows from Lemma 4.6 that the solutions to (4.2) have the same large s behavior as do elements in the modulis spaces that are defined in [KM2]. Meanwhile, the $s \to -\infty$ limit is analyzed for one particular case in Section 4b. In particular, Lemma 4b asserts for this one case that the relevant solutions to (4.2) have the same $s \to -\infty$ behavior as do elements in the moduli spaces that [KM1] define in their Chapter 24. The story for other cases is much the same. Given that the limiting behavior is such as to put the solution in the appropriate moduli spaces, then a straightforward application of the techniques from [KM2] and [KM1] prove that the set of gauge equivalence classes to (4.2) comprise a finite set when $\hat{\mu}$ is chosen in a suitably generic fashion, and that each such class has an associated ±1 assignment. Theorem 4.1 likewise follows in a direct fashion from what is said in these parts of [KM1] and [KM2]. The fact is that the discussion in Chapter 24 of [KM1] is written so as to readily accommodate the case where the 4-manifold has both cylindrical and asymptotically flat almost Kahler ends.

      To put all of this in a larger context, view $[0, \infty) \times M$ as a manifold with boundary, and Theorem 4.1 provides only another example of what is observed by Donaldson [D] in the context of SU(2) gauge theory: A construction that gives a numerical gauge theory

invariant to a 4-manifold with no boundary generalizes to a manifold with boundary so as to give an invariant in the boundary's Floer homology or cohomology.

As noted, Theorem 4.2 is the analog in this context of the vanishing result given by Corollary C of [MR]. The proof amounts to adapting the proof given by Mrowka and Rollins [MR] for their Corollary C to the context of $\mathbb{R} \times M$ with its asymptotically flat almost Kahler $s \to \infty$ end and its cylindrical $s \to -\infty$ end. This adaptation is completely straightforward with the analysis in Chapter 24 of [KM1] to control behavior on the $s \to -\infty$ end.

The proof of Theorem 1.4 exploits the dichotomy between what is said in Theorem 4.2 and what is said in the next proposition.

**Proposition 4.3**: *There exists $r_{I'} \geq r_I$ and $\varepsilon_0 > 0$ such that when $r \geq r_{I'}$ and $\|\hat{\mu}\|_{C^6} < \varepsilon_0$, then there is precisely one solution to (4.2), and its $s \to -\infty$ limit defines the class $c_I(r)$.*

This proposition is proved in Sections 4b and 4c modulo a proposition whose proof occupies Section 5.

Theorem 4.2 says that there is a cocycle, $v$, on the degree -1 Seiberg-Witten Floer complex such that $t_r(\mathfrak{z}) = v(\delta \mathfrak{z})$. Meanwhile, Proposition 4.3 guarantees that $t_r(c_I(r))$ is never zero. Thus, Proposition 4.3 and Theorem 4.2 demand that $\delta c_I(r) \neq 0$. This understood, let $\mathcal{B}(r)$ now denote the set of degree -1 Seiberg-Witten Floer cycles of the form $\delta(c_I(r) + \mathfrak{w})$ where $\mathfrak{w}$ is a linear combination of generators with $\mathfrak{a}^f < \mathfrak{a}^f(c_I(r))$. From what was just said, $\mathcal{B}(r) \neq 0$. For each $\mathfrak{n} \in \mathcal{B}(r)$, use $\mathfrak{a}^f(\mathfrak{n}, r)$ to denote the maximum of $\mathfrak{a}^f$ on the generators that appear in $\mathfrak{n}$ with non-zero weight. Now set $\mathfrak{a}^f_I(r)$ to denote the infimum of the set $\{\mathfrak{a}^f(\mathfrak{n}, r)\}_{\mathfrak{n} \in \mathcal{B}(r)}$. Define as before $\hat{E}(r)$, $\mathfrak{f}(r)$ and $v(r)$ using $\mathcal{B}(r)$.

The analog here of Proposition 2.9 holds; and as in the cases studied previously, there is a sequence $\{r_n, (A_n, \psi_n)\}$ with $\{r_n\}$ increasing and unbounded, with $(A_n, \psi_n)$ satisfying the $r = r_n$ and $\mu$ version of (1.3), with degree -1, and satisfying one or the other of the options in (2.3). As before, the second option is ruled out by Proposition 5.1 in [T1]. As a consequence, $\{E(A_n)\}$ is bounded. Since each $(A_n, \psi_n)$ has degree -1, Lemma 5.4 in [T1] guarantees the existence of some positive $\delta$ such that $\sup_M (1 - |\psi_n|) > \delta$ for all sufficiently large n. This being the case, the sequence $\{(r_n, (A_n, \psi_n)\}$ meets all of the conditions set forth by Theorem 1.5. Theorem 1.5 provides the set of closed integral curves of v for Theorem 1.4.

**b) Proof of Proposition 4.3: Existence**

This part of the proof establishes the existence of at least one solution to (4.2), this with $s \to -\infty$ limit equal to $\mathfrak{u} \cdot c_I(r)$ with $\mathfrak{u}$ a smooth map from M to $S^1$. As the reader will see, the arguments given in what follows are analogs for $\mathbb{R} \times M$ of arguments from Section 2d that prove Proposition 2.8.

To start the existence proof, introduce as before $(A_I, \psi_I)$. Here $A_i$ is a trivial connection on $I_\mathbb{C}$ and $\psi_I$ is $A_I$ covariantly constant. Moreover, $D_{A_I} \psi_I = 0$ and $\psi_I^\dagger \tau \psi_I = ia$. Thus, the constant map $s \to \mathfrak{d}_I(s) = (A_I, \psi_I)$ satisfies the version of (4.2) where the right hand side has $\mu$ and $\hat{\mu}$ equal to zero and $i\varpi_K$ absent. A solution to (4.2) is obtained from a section over $\mathbb{R} \times M$ of the bundle $iT^*M \oplus \mathbb{S} \oplus i\mathbb{R}$ over $\mathbb{R} \times M$. To say more, introduce $\mathfrak{b}_0 = (b_0, \eta_0, 0)$ to denote the fixed point to (2.5) with norm $\|\mathfrak{b}_0\| \leq c_0 r^{-1/2}$. The desired triple $(c, \varsigma, \phi) \in C^\infty(\mathbb{R} \times M; iT^*M \oplus \mathbb{S} \oplus i\mathbb{R})$ must have limit zero as $s \to \pm\infty$ and solve

- $\frac{\partial}{\partial s} c + \mathrm{T}'c + *dc - d\phi - 2^{-1/2}r^{1/2}e^\mathrm{T}(\psi_I^\dagger \tau \varsigma + \varsigma^\dagger \tau \psi_I) - r^{1/2}e^\mathrm{T}\varsigma^\dagger \tau \varsigma - \chi'b - \frac{i}{2}r^{-1/2}\hat{\mu}$
  $\qquad - 2^{-1/2}r^{1/2}e^\mathrm{T}\chi(\eta_0^\dagger \tau \varsigma + \varsigma^\dagger \tau \eta_0) - r^{1/2}\chi(1-\chi)\eta_0^\dagger \tau \eta_0 = 0$,
- $\frac{\partial}{\partial s}\varsigma + D_{A_I}\varsigma + 2^{1/2}r^{1/2}e^\mathrm{T}(\mathrm{cl}(c)\psi_I + \phi\psi_I) + 2^{1/2}r^{1/2}e^\mathrm{T}(\mathrm{cl}(c)\varsigma + \phi\varsigma) - \chi'\eta_0$
  $\qquad + 2^{1/2}r^{1/2}\chi(\mathrm{cl}(c)\eta_0 + \mathrm{cl}(b_0)\eta + \phi\eta_0) + 2^{1/2}r^{1/2}\chi(1-\chi)\mathrm{cl}(b_0)\eta_0 = 0$,
- $\frac{\partial}{\partial s}\phi + *d*c - 2^{-/2}r^{1/2}e^\mathrm{T}(\varsigma^\dagger \psi_I - \psi_I^\dagger \varsigma) = 0$.

(4.3)

In this last equation and in what follows, d denotes the exterior derivative along the M factor in $\mathbb{R} \times M$. Also, $\mathrm{T}'$ and $\chi'$ denote the respective derivatives of $\mathrm{T}$ and $\chi$. Here is one further constraint on $(c, \varsigma, \phi)$: Let

$$\sigma(s, \cdot) = 2^{1/2}r^{1/2}\int_0^s e^{\mathrm{T}(t)}\phi(t, \cdot)\, dt \ .$$

(4.4)

The function $\sigma$ must have a limit as $s \to \pm\infty$ as a smooth function on M. Granted that such is the case, let $\sigma_+$ denote the $s \to \infty$ limit of $\sigma$. Then the pair $(A, \psi)$ with the connection $A = A_I + 2^{-/2}r^{1/2}e^\mathrm{T}(\chi b_0 + c) - d(\sigma - \sigma_+)$ and the spinor $\psi = e^{\sigma - \sigma_+}(\psi_I + \chi\eta_0 + \varsigma)$ solves (4.2) and has the required limits as $s \to \pm\infty$.

To find $(c, \varsigma, \phi)$, introduce the Hilbert space completion of the space of compactly supported elements in $C^\infty(\mathbb{R} \times M; iT^*M \oplus \mathbb{S} \oplus i\mathbb{R})$ using the norm whose square gives the function

$$\mathfrak{h} \to \|\mathfrak{h}\|_\mathbb{X}^2 = \int_\mathbb{R} e^{2\mathrm{T}}\int_M (|\tfrac{\partial}{\partial s}\mathfrak{h}|^2 + |\nabla_I \mathfrak{h}|^2) + \tfrac{1}{4}re^{2\mathrm{T}}|\mathfrak{h}|^2) \ .$$

(4.5)

This Hilbert space is denoted by $\mathbb{X}$ in what follows. A solution to (4.3) is found in $\mathbb{X}$ using a contraction mapping argument much like that used to prove Proposition 2.8.

To set the stage for this contraction mapping business, introduce $\mathbb{L}$ to denote the completion of the space of compactly supported elements in $C^\infty(\mathbb{R} \times M; iT^*M \oplus \mathbb{S} \oplus i\mathbb{R})$ using the norm whose square is defined by

$$\|\mathfrak{h}\|_{\mathbb{L}}^2 = \int_{\mathbb{R}} e^{2T} \int_M |\mathfrak{h}|^2 \ .$$

(4.6)

Next, introduce the operator $\mathcal{D}: \mathbb{X} \to \mathbb{L}$ that sends $\mathfrak{h} = (c, \varsigma, \phi)$ to the element with respective $iT^*M$, $\mathbb{S}$ and $i\mathbb{R}$ components

- $\frac{\partial}{\partial s} c + T'c + *dc - d\phi - 2^{-1/2} r^{1/2} e^T (\psi_I^\dagger \tau \varsigma + \varsigma^\dagger \tau \psi_I)$,
- $\frac{\partial}{\partial s} \varsigma + D_{A_I} \varsigma + 2^{-1/2} r^{1/2} e^T (cl(c)\psi_I + \phi \psi_I)$,
- $\frac{\partial}{\partial s} \phi + *d*c - 2^{-1/2} r^{1/2} e^T (\varsigma^\dagger \psi_I - \psi_I^\dagger \varsigma)$.

(4.7)

Note that $\mathcal{D}$ can be written as $\mathcal{D} = \frac{\partial}{\partial s} + \mathfrak{L}_s + T' \Pi_{T^*M}$, where the notation uses $\mathfrak{L}_s$ to denote the version of (3.1) that has $(A, \psi) = (A_I, \psi_I)$ and $r$ replaced by $r e^{2T(s)}$. Meanwhile, $\Pi_{T^*M}$ denotes the fiberwise projection from the bundle $iT^*M \oplus \mathbb{S} \oplus i\mathbb{R}$ to $iT^*M$.

**Lemma 4.4**: *There exists $r_* > 1$ such that for $r \geq r_*$, the operator $\mathcal{D}$ is Fredholm with index zero. Moreover, $\|\mathcal{D}\mathfrak{h}\|_{\mathbb{L}} \geq \|\mathfrak{h}\|_{\mathbb{X}}$.*

*Proof of Lemma 4.4*: This follows from (5.23) in [T1] with some integration by parts.

Lemma 4.4 guarantees that $\mathcal{D}$ is invertible when $r \geq r_*$.
    One more lemma is required for the set up, this the following weighted version of a dimension 4 Sobolev inequality. Its proof is straightforward and so omitted.

**Lemma 4.5**: *Suppose that $r \geq 1$. Let $(\mathfrak{h}, \mathfrak{h}') \to \mathfrak{h} * \mathfrak{h}'$ denote a given homomorphism from the bundle $\otimes_2 (iT^*M \oplus \mathbb{S} \oplus i\mathbb{R})|_M$ to $(iT^*M \oplus \mathbb{S} \oplus i\mathbb{R})|_M$. The homomorphism $(\mathfrak{h}, \mathfrak{h}') \to e^T \mathfrak{h} * \mathfrak{h}'$ over $\mathbb{R} \times M$ induces a smooth, bilinear map from $\mathbb{X} \oplus \mathbb{X}$ to $\mathbb{L}$ whose norm obeys $\|e^T \mathfrak{h} * \mathfrak{h}'\|_{\mathbb{L}} \leq c_0 \|\mathfrak{h}\|_{\mathbb{X}}^2$.*

Granted these last two lemmas, it follows that what is written in (4.3) defines a smooth map from $\mathbb{X}$ to $\mathbb{L}$. Write this map as $\mathfrak{h} \to \mathcal{D}\mathfrak{h} + r^{1/2} e^T \mathfrak{h} * \mathfrak{h} + 2r^{1/2} \chi \mathfrak{h}_0 * \mathfrak{h} - r^{-1/2} \mathfrak{v}_\mathbb{X}$, and it then follows that the solutions to (4.3) in $\mathbb{X}$ constitute the fixed point set of the self map

$$\mathfrak{h} \to T_\mathbb{X}(\mathfrak{h}) = \mathcal{D}^{-1}(r^{-1/2} \mathfrak{v}_\mathbb{X} - r^{1/2} \chi \mathfrak{h}_0 * \mathfrak{h} - r^{1/2} e^T \mathfrak{h} * \mathfrak{h}) \ .$$

(4.8)

To see that $T_\mathbb{X}$ has a fixed point for large $r$, reintroduce the norm $\|\cdot\|_{\mathbb{H}}$ that is defined by (2.6). By design, the fixed point $\mathfrak{h}_0 = (b_0, \eta_0, 0)$ to (2.5) obeys $\|\mathfrak{h}_0\|_{\mathbb{H}} \leq c_0 r^{-1/2}$.

It then follows from (2.8) that the $L^4$ norm on M of $\mathfrak{h}_0$ is bounded by $c_0 r^{-5/8}$. Meanwhile, that of $\mathfrak{h}$ on any given slice $\{s\} \times M \subset \mathbb{R} \times M$ is bounded by $r^{-1/8} \|\mathfrak{h}|_s\|_{\mathbb{H}}$. As a consequence, $\nu_{\mathbb{X}}$ and $\mathfrak{h}_0 * \mathfrak{h}$ obey

$$\|\nu_{\mathbb{X}}\|_{\mathbb{L}} \leq c_0 (\|\hat{\mu}\|_{\mathbb{L}} + r^{-1/4}) \text{ and } \|\chi\, \mathfrak{h}_0 * \mathfrak{h}\|_{\mathbb{L}} \leq c_0 r^{-3/4} \|\mathfrak{h}\|_{\mathbb{X}}$$

(4.9)

when r is large. It then follows from Lemmas 4.4 and 4.5 that

$$\|T_{\mathbb{X}}(\mathfrak{h})\|_{\mathbb{X}} \leq c_0 (r^{-1/2} \|\hat{\mu}\|_{\mathbb{L}} + r^{-1/4} \|\mathfrak{h}\|_{\mathbb{X}} + r^{1/2} \|e^{\tau} \mathfrak{h} * \mathfrak{h}\|_{\mathbb{L}})$$

(4.10)

As a consequence, the map $T_{\mathbb{X}}$ on the ball in X of radius $R = 4 c_0^{-1} \|\hat{\mu}\|_{\mathbb{L}} r^{-1/2}$ obeys

$$\|T_{\mathbb{X}}(\mathfrak{h})\|_{\mathbb{X}} \leq \tfrac{1}{4} r^{-1/2} R (1 + 4 c_0 (r^{-1/4} + R))$$

(4.11)

This understood, there exists $r_0 > 1$, $\varepsilon_0 > 0$ and $R_0 > 0$ such that $T_{\mathbb{X}}$ maps the ball of radius $R_0 r^{-1/2}$ to itself when $r > r_0$ and when $\|\hat{\mu}\|_{\mathbb{L}} \leq \varepsilon_0$. A very similar calculation proves that the constants $r_0$ and $\varepsilon_0$ can be chosen so that $T_{\mathbb{X}}$ maps this ball to itself as a contraction mapping when $r > r_0$ and $\|\hat{\mu}\|_{\mathbb{L}} < \varepsilon_0$. The contraction mapping theorem finds a unique fixed point of $T_{\mathbb{X}}$ in the ball of radius $R_0 r^{-1/2}$ for such r and $\hat{\mu}$. Very much standard elliptic regularity arguments using the equation in (4.3) prove that this fixed point is a smooth section of $iT^*M \oplus \mathbb{S} \oplus i\mathbb{R}$ over $\mathbb{R} \times M$. Note for later reference that the elliptic regularity assertions are of the following sort:

*There exists an r-independent constant $\rho > 0$ and, given $r \geq 1$ and integers $k, n \geq 0$, there exists $c_{k,n,r}$, with the following significance: Let $B \subset \mathbb{R} \times M$ denote a ball of radius $\rho$ and let $B' \subset B$ denote the concentric ball of radius $\tfrac{1}{2} \rho$ as defined with the product metric. Let $\mathfrak{h}$ denote a solution to (4.3) on B. Then*

$$\sup_{x \in B'} |\tfrac{d^n}{ds^n} (\nabla_1)^k \mathfrak{h}|^2 \leq c_{k,n,r}\, e^{(n+k)\tau(p)} \int_B e^{2\tau} (|\tfrac{\partial}{\partial s} \mathfrak{h}|^2 + |\nabla_1 \mathfrak{h}|^2 + |\mathfrak{h}|^2)$$

(4.12)

As is explained momentarily, the existence of the desired $s \to \pm\infty$ limits of (4.4)'s function $\sigma$ follows from.

**Lemma 4.6**: *There exists $\kappa > 0$ and $\varepsilon_0 > 0$, and given $\varepsilon \in (0, \varepsilon_0)$, there exists $r_\varepsilon \geq 1$, and these have the following significance: Suppose that $r \geq r_\varepsilon$ and $s_0 > 2$  Let $\mathfrak{h} = (c, \varsigma, \phi)$ denote a smooth solution to (4.3) on either the half cylinder $(-\infty, -s_0] \times M$, or on the half cylinder $[s_0, \infty) \times M$. Suppose in addition that $|\mathfrak{h}| < \varepsilon$. Then*

- $\lim_{s \to -\infty} \exp(\kappa r^{1/2}|s|)|\mathfrak{h}|$ *is finite in the case when $\mathfrak{h}$ is defined on* $(-\infty, -s_0] \times M$.
- $\lim_{s \to \infty} \exp(\kappa r^{1/2} e^s)|\mathfrak{h}|$ *is finite in the case when $\mathfrak{h}$ is defined on* $[s_0, \infty) \times M$.

What with (4.12) and this exponential decay, it follows that derivatives of $\mathfrak{h}$ to any given order decay at an exponential rate as $s \to -\infty$ with rate constant proportional to $r^{1/2}$ and it decays as $s \to \infty$ as $\exp(-\kappa r^{1/2} e^s)$. For example, such is the case for the solution found in the previous subsection. Thus, the function $\sigma$ that is defined in (4.4) limits as $s \to \pm\infty$ in the $C^\infty$ topology to a smooth function on M. In the case when $\mathfrak{h}$ is the solution to (4.3) found in the previous subsection, the second item of the lemma plus (4.12) proves that the function $\sigma$ in (4.4) also has a limit as $s \to \infty$ in the $C^\infty$ topology on the space of smooth functions on M.

***Proof of Lemma 4.6***: To prove the first item, let $\mathcal{L}$ denote the version of (1.7) that is defined by $(A_I + 2^{1/2} r^{1/2} b_0, \psi_I + \eta_0)$. It is a consequence of Lemma 5.4 in [T1] that $\mathcal{L}$ has no zero eigenvalue on $L^2(M; iT^*M \oplus \mathbb{S} \oplus i\mathbb{R})$. Indeed it follows from this lemma that $\|\mathcal{L}\mathfrak{b}\|_2 \geq (1 - c_0 r^{-3/2})\|\mathfrak{b}\|_\mathbb{H}$ for any $\mathfrak{b} \in \mathbb{H}$. Thus, all eigenvalues of $\mathcal{L}$ have absolute value greater than $\frac{1}{4} r^{1/2}$ when r is large. Let $\Pi_+$ denote the $L^2$ orthogonal projection on M onto the span of those eigenvectors of $\mathcal{L}$ that have positive eigenvalue. Let $\Pi_- = (1 - \Pi_+)$. Let $h_\pm(s) = \|\Pi_\pm \mathfrak{h}|_s\|_2$.

Fix $\varepsilon > 0$ and suppose that $s_0 > 0$ and that $\mathfrak{h}$ is a solution to (4.3) such that $|\mathfrak{h}| \leq \varepsilon$ where $s < -s_0$. It then follows as a consequence of (4.3) that the functions $h_\pm$ on $(-\infty, -s_0)$ obey

$$\tfrac{d}{ds} h_+ + \tfrac{1}{4} r^{1/2} h_+ \leq c_0 \varepsilon r^{1/2} (h_+ + h_-) \quad\text{and}\quad \tfrac{d}{ds} h_- - \tfrac{1}{4} r^{1/2} h_- \geq -c_0 \varepsilon r^{1/2} (h_+ + h_-).$$

(4.13)

Taking linear combinations of these two equations results in the following conclusion: There exist r-independent constants $c_1 > 0$ and $c_2$ such that if $\varepsilon < c_1$, then $h_* = h_+ - c_2 \varepsilon h_-$ obeys

$$\tfrac{d}{ds} h_* + \tfrac{1}{8} r^{1/2} h_* \leq 0.$$

(4.14)

Integrating this last equation finds that

$$h_*(s) \geq e^{r^{1/2}(s'-s)/8} h_*(s') \text{ for any } s \leq s' < s_0.$$

(4.15)

Since $h_\pm(s)$ have bounded absolute value as $s \to -\infty$, this last equation requires that $h_*(s) \leq 0$ for all $s \leq s_*$. This says that $h_+ \leq c_2 \varepsilon h_-$ for $s \leq s_\varepsilon$. Granted the latter, it follows from the right most equation in (4.15) that there exists $c_1' > 0$ such that when $\varepsilon < c_1'$, then

$$\tfrac{d}{ds} h_- - \tfrac{1}{8} r^{-1/2} h_- \geq 0 \quad \textit{where } s \leq s_0.$$

(4.16)

Integrating (4.16) finds that $h_-(s) \leq e^{r^{1/2}(s-s')/8} h_-(s')$ when $s \leq s'$. Thus, $h_-(s)$ limits to zero at an exponential rate as $s \to -\infty$, and since $h_+ \leq c_2 \varepsilon h_-$, this is true for $h_+$ also.

Consider next the story for the second item. Thus, assume now that $\mathfrak{h}$ is a solution to (4.3) on $[s_0, \infty) \times M$ with $|\mathfrak{h}| < \varepsilon$. To study this case, let $\mathcal{L}_s$ denote the operator that is defined by the version of (1.7) that has $re^{2s}$ replacing $r$ and $(A_I, \psi_I)$ for the pair of connection and section of $\mathbb{S}$. In addition, take $\mathfrak{t}$ and $\mathfrak{s}$ equal to zero. Let $\prod_{s+}$ and $\prod_{s-}$ denote the respective projections to the span of the eigenvectors of $\mathcal{L}_s$ with positive and negative eigenvalues.

**Lemma 4.7**: *The operators $\prod_{s\pm}$ on $L^2(M; iT^*M \oplus \mathbb{S} \oplus \mathbb{R})$ vary in a real analytic fashion as a function of s; and the norm of $\tfrac{d}{ds} \prod_{s\pm}$ has an* r *and* s *independent bound.*

This lemma is proved momentarily.

With Lemma 4.7 in hand, let $h_\pm(s) = \| \prod_{s\pm} \mathfrak{h}|_s \|_2$. These functions are defined on the half line $[s_0, \infty)$ where they obey the following analogs of (4.14):

$$\tfrac{d}{ds} h_+ + \tfrac{1}{4} r^{1/2} e^s h_+ \leq c_0 (\varepsilon r^{1/2} e^s + 1)(h_+ + h_-) \quad \textit{and} \quad \tfrac{d}{ds} h_- - \tfrac{1}{4} r^{1/2} e^s h_- \geq -c_0 (\varepsilon r^{1/2} e^s + 1)(h_+ + h_-).$$

(4.17)

Note that the extra factor of $c_0(h_+ + h_-)$ that appears on the right hand side of each inequality comes from s-derivatives of $\prod_{s\pm}$ courtesy of Lemma 4.7. With the preceding understood, change variables from s to $t = e^s$. Doing so makes (4.17) look like a version of (4.13) when $r^{1/2} e^s > 1/\varepsilon$. This understood, the arguments that are used for (4.13) can be used here with the roles of $h_+$ and $h_-$ reversed to prove the second item of the lemma.

**Proof of Lemma 4.7**: The fact that $\prod_{s+}$ varies in a real analytic fashion follows from the results in Chapter 7 of [Ka] given that the arguments in Section 5e of [T] prove that $\mathcal{L}_s$ has no eigenvalues with absolute value less than $\tfrac{1}{4} r^{1/2} e^s$. The bound on the norm of the derivative follows by writing $\prod_{s+}$ as a contour integral using its resolvent:

$$\prod_{s+} = \tfrac{1}{2\pi i} \int_C \frac{1}{z - \mathcal{L}_s} dz$$

(4.18)

where the contour $C \subset \mathbb{C}$ is the path that starts at $(\infty, ir^{1/2} e^s)$ and moves with constant imaginary part to $(0, ir^{1/2} e^s)$, then crosses the real axis to $(0, -ir^{1/2} e^s)$ and then moves with constant imaginary part to $(\infty, -ir^{1/2} e^s)$. The derivative of $\prod_{s+}$ is thus given by the contour integral

$$\tfrac{d}{ds} \Pi_{s+} = \tfrac{1}{2\pi i} \int_C \frac{1}{z - \mathfrak{L}_s} \tfrac{d}{ds} \mathfrak{L}_s \frac{1}{z - \mathfrak{L}_s} dz \ .$$

(4.19)

Now, the derivative of $\mathfrak{L}_s$ is an endomorphism of $iT^*M \oplus \mathbb{S} \oplus i\mathbb{R}$ with norm bounded by $c_0 r^{1/2} e^s$. Meanwhile, the operator norm of $(z - \mathfrak{L}_s)^{-1}$ at a point $(x, \pm i r^{1/2} e^s)$ on the coutour is bounded by $(x + r^{1/2} e^s)^{-1/2}$. Thus, the integral along the horizontal parts of the contour in (4.19) define an operator with norm bounded by

$$c_0 \int_0^\infty \frac{r^{1/2} e^s}{(x + r^{1/2} e^s)^2} dx \ ,$$

(4.20)

which has an r and s independent bound. The operator norm of $(z - \mathfrak{L}_s)^{-1}$ along the vertical part of the contour in (4.19) is bounded by $c_0 r^{-1/2} e^{-s}$ since there are no eigenvalues of $\mathfrak{L}_s$ whose eigenvalue has absolute value less than $\tfrac{1}{4} r^{1/2} e^s$. Given this bound, it then follows that the vertical part of the contour integral in (4.19) defines an operator whose norm has an r and s independent upper bound.

The r-independence of the sign that is associated to this solution follows from with a proof that the linearization of (4.3) at the solution $\mathfrak{h}$ has just trivial kernel and cokernel when r is large. The latter fact follows from the contraction mapping construction of $\mathfrak{h}$. Indeed, an element in the kernel of the linearization of (4.3) at $\mathfrak{h}$ is in the kernel of the differential of the map $T_\mathbb{X}$ at $\mathfrak{h}$, and this differntial has trivial kernel because $T_\mathbb{X}$ is a contraction at $\mathfrak{h}$.

### c) Proof of Proposition 4.3: Uniqueness

The uniqueness assertion in Proposition 4.3 is proved with the help of the following:

**Proposition 4.8**: *Given $\varepsilon > 0$, there exists $r_\varepsilon > 1$ with the following significance: Suppose that $r \geq r_\varepsilon$ and that $s \to \mathfrak{d}(s) = (A, \psi = (\alpha, \beta))$ is a solution to (4.2). Then*

$$(1 - |\alpha|^2)^{1/2} + |\beta| + r^{-1/2}(|\tfrac{d}{ds} \alpha| + |\nabla \alpha|) + r^{-1}(|\tfrac{\partial^2}{\partial s^2} \alpha| + |\tfrac{\partial}{\partial s} \nabla \alpha| + |\nabla \tfrac{\partial}{\partial s} \alpha| + |\nabla^2 \alpha| < \varepsilon \ .$$

(4.21)

*at all points in $\mathbb{R} \times M$.*

A proof of this proposition is given shortly, so assume it to be true so as to complete the proof of Proposition 4.3.

The proof of the uniqueness assertion in Proposition 4.3 uses Proposition 4.8 to obtain a solution to (4.3) from a solution to (4.2). This first step in the proof invokes some conclusions from the next lemma. To set the stage for this lemma, suppose that $I \subset \mathbb{R}$ is either all of $\mathbb{R}$, or $(-\infty, s_1]$, or $[s_1, \infty)$, with $s_1 \in \mathbb{R}$.

**Lemma 4.9**: *There exists $\kappa > 4$, and given $\varepsilon \in (0, \kappa^{-1})$, there exists $r_\varepsilon > 1$ with the following significance: Suppose that $r \geq r_\varepsilon$ and that $s \to \mathfrak{d}(s) = (A, \psi)$ satisfies the first two equations in (4.2) on $I \times M$. Also assume that (4.21) holds on $I \times M$. Then there exists*
- *A subset $I' \subset I$ that is equal to $\mathbb{R}$ when $I = \mathbb{R}$ and is a half-line otherwise.*
- *A solution, $\mathfrak{h} = (c, \varsigma, \phi) \in C^\infty(I' \times M; iT^*M \oplus \mathbb{S} \oplus i\mathbb{R})$ to (4.3) on $I' \times M$ with $|\mathfrak{h}| \leq \kappa\varepsilon$ and such that the pairs $(A_I + 2^{1/2}r^{1/2}e^T(\chi b_0 + c - \phi ds), \psi_I + \chi\eta_0 + \varsigma)$ and $(A, \psi)$ of connection on the trivial bundle over $I' \times M$ and section over $I' \times M$ of $\mathbb{S}_I$ are gauge equivalent.*

**Proof of Lemma 4.9**: Since $\alpha \neq 0$ on $I \times M$, there exists a map $\hat{u}: I \times M \to S^1$ with the property that $\hat{u}\alpha = |\alpha| I_\mathbb{C}$. It then follows that $A - \hat{u}^{-1}d\hat{u} - (\hat{u}^{-1}\frac{\partial}{\partial s}\hat{u}) ds = A_I + \hat{a} - \hat{a}_0 ds$, where $\hat{a}$ is the imaginary part of the $\mathbb{C}$-valued 1-form $\alpha^{-1}\nabla\alpha$ and $\hat{a}_0$ is the imaginary part of the $\mathbb{C}$-valued function $\alpha^{-1}\frac{\partial}{\partial s}\alpha$. Write $\hat{a} = 2^{1/2}r^{1/2}e^T(\chi b_0 + b)$, write $\psi = \psi_I + \chi\eta_0 + \eta$, and write $\hat{a}_0 = 2^{1/2}r^{1/2}e^T\sigma$. Then the triple $\mathfrak{b} = (b, \eta, \sigma)$ obeys $|\mathfrak{b}| \leq c_0\varepsilon$. In addition, the top two equations in (4.3) are obeyed with $\mathfrak{b}$ replacing $\mathfrak{h}$. To obtain $\mathfrak{h}$ as described by the lemma, consider an additional gauge transformation of the form $e^x$ where $x: I \times M \to i\mathbb{R}$. Write $u = e^x \hat{u}$ and use $u$ and $(A, \psi)$ to define the triple $\mathfrak{h} = (c, \varsigma, \phi)$ as in the lemma. Thus,

$$c = b - 2^{-1/2}r^{-1/2}e^{-T}dx, \quad \varsigma = (e^x - 1)(\psi_I + \chi\eta_0) + e^x\eta \quad and \quad \phi = \sigma - 2^{-1/2}r^{-1/2}e^{-T}\tfrac{\partial}{\partial s}x \,.$$
(4.22)

Note that the top two equations in (4.3) are still obeyed by $\mathfrak{h}$. The bottom equation in (4.3) is obeyed if $x$ obeys an equation with the form

$$-\tfrac{\partial^2}{\partial s^2}x + T'\tfrac{\partial}{\partial s}x + d^*dx + 2r\,e^{2T}x + \tfrac{\partial}{\partial s}\sigma + d^*b + r\,e^{2T}\wp(x) = 0$$
(4.23)

where $\wp(\cdot)$ at each point in $I \times M$ is a smooth function on $(-c_0, c_0) \subset \mathbb{R}$ with the property that $|\wp(t)| \leq c_0(t^2 + |\eta|)$ and $|(\tfrac{d}{dt}\wp)(t)| \leq c_0 t$ at any given point in $I \times M$. The desired bound $|\mathfrak{h}| \leq \kappa\varepsilon$ is obeyed if $x$ solves (4.23) and $|x| + r^{-1/2}e^{-T}(|\tfrac{\partial}{\partial s}x| + |dx|) \leq c_0\varepsilon$. The existence of the desired solution can be derived using 4-dimensional analogs of the arguments that proved Lemma 2.11. There is nothing tricky about this, and so the details are left to the reader.

Given Proposition 4.8, use Lemma 4.9 to deduce the following: Give $\varepsilon > 0$, there exists $r_\varepsilon > 0$ such that any solution to an $r > r_\varepsilon$ version of (4.2) supplies a solution, $\mathfrak{h}$, to (4.3) with $|\mathfrak{h}| < \varepsilon$. Use Lemma 4.6 to deduce that $\mathfrak{h}$ is in the Hilbert space $\mathbb{X}$. Granted that $\mathfrak{h}$ is in $\mathbb{X}$, if it has small norm, it will have to be the solution to (4.3) from Subsection 4b. Thus, it is a fixed point of the map $T_\mathbb{X}$ that is depicted in (4.8). This being the case, (4.10) implies that

$$\|\mathfrak{h}\|_\mathbb{X} \leq c_0 (r^{-1/2} \|\hat{\mu}\|_\mathbb{L} + r^{-1/4} \|\mathfrak{h}\|_\mathbb{X} + r^{1/2} \|e^T \mathfrak{h} * \mathfrak{h}\|_\mathbb{L}) .$$

(4.24)

Since $|\mathfrak{h}| < \varepsilon$, this last inequality implies that

$$\|\mathfrak{h}\|_\mathbb{X} \leq 2c_0 (r^{-1/2} \|\hat{\mu}\|_\mathbb{L} + r^{1/2} \varepsilon \|e^T \mathfrak{h}\|_\mathbb{L})$$

(4.25)

when $r > c_0$. Meanwhile, it follows from (4.5) and (4.6) that $\|e^T \mathfrak{h}\|_\mathbb{L} \leq 2r^{-1/2} \|\mathfrak{h}\|_\mathbb{X}$ and so (4.25) implies that

$$\|\mathfrak{h}\|_\mathbb{X} \leq 2c_0 r^{-1/2} \|\hat{\mu}\|_\mathbb{L} + 4c_0 \varepsilon \|\mathfrak{h}\|_\mathbb{X} .$$

(4.26)

Thus, if $\varepsilon \ll (4c_0)^{-1}$, then $\mathfrak{h}$ will lie inside the ball in $\mathbb{X}$ where $T_\mathbb{X}$ acts as a contraction mapping. In this event, $\mathfrak{h}$ is the solution from Subsection 4b.

## 5. The proofs of Propositions 4.8 and 2.9

The proof of Proposition 4.8 is given first, as many of the estimates and techniques used in the latter are then applied with only minor modifications to prove Proposition 2.9. The proof of Proposition 2.9 appears in the final subsection.

The proof of Proposition 4.8 invokes a set of apriori bounds for any solution $s \to (A(s), \psi(s))$ to (4.2). The first such bound is the analog for $\mathbb{R} \times M$ of what is stated in Lemma 1.6.

**Lemma 5.1**: *There exists a constant, $\kappa > 1$ such that if $r \geq \kappa$ and if $(A, \psi)$ is a solution to (4.2), then*
- $|\alpha| \leq 1 + \kappa r^{-1} e^{-2T}$.
- $|\beta|^2 \leq \kappa r^{-1} e^{-2T} (1 - |\alpha|^2) + \kappa r^{-2} e^{-4T}$.

The next apriori bound concerns the curvature of the connection A. This proposition refers to the functional $\mathfrak{a}$ that is depicted in (1.4).

**Proposition 5.2**: *There exists a constant $\kappa_* > 1$ with the following significance: Suppose that $r \geq \kappa_*$ and that $s \to \mathfrak{d}(s) = (A, \psi)$ is a solution to (4.2). Then $|\frac{\partial}{\partial s} A| + |B_A| \leq \kappa_* e^{2T} r$.*

The last of the required results asserts a part of Proposition 4.8:

**Proposition 5.3**: *Let $(A, \psi)$ and $\kappa_*$ be as described in Proposition 5.2. Given $\delta > 0$, there exists $r_\delta > \kappa_*$ such that if $r > r_\delta$ and $(A, \psi)$ is a solution to (4.2), then $1 - |\alpha|^2 \leq \delta$.*

Lemma 5.1 is proved in Subsection 5b. The proofs of Propositions 5.2 and 5.3 occupy the remainder of Section 5. First up is the proof of Proposition 5.2 for points where $s \geq -2$, and then comes the proof of Proposition 5.3 for points where $s \geq -1$. Next is the proof of Propositions 5.2 for points where $s \leq -1$; the latter starts in Subsection 5g. The final subsection proves Proposition 5.3 for points where $s \leq -1$.

Note that Proposition 5.3 has the following implication: If $c_* > 0$ has been fixed, if r is sufficiently large, and if $s \to \mathfrak{d}(s) = (A, \psi)$ is a solution to (4.2) with $\lim_{s \to -\infty} \mathfrak{a}(\mathfrak{d}(s)) \leq c_* r^2$, then $\lim_{s \to \infty} \mathfrak{d}(s) = \mathfrak{c}_I(r)$. Note in this regard that (4.11) implies that $\mathfrak{a}(\mathfrak{c}_I(r)) \leq c_0 r^{-1/2}$.

**a) Proof of Proposition 4.8 given Lemma 5.1 and Propositions 5.2 and 5.3**

Needed are upper bounds for the derivatives of $\alpha$ up through second order. These are obtained using a scaling argument of the sort used for the proof of Proposition 4.2 in [T3]. To elaborate, the bound is proved by rescaling Gaussian coordinates centered at a given point $p \in \mathbb{R} \times M$ so that the ball of radius $r^{-1/2} e^{-T(p)}$ has radius 1 with respect to the Euclidean metric for the new coordinates. When the pair $(A, \psi)$ are pulled back to $\mathbb{R}^4$ by this scaling map, the equations in the first two lines of (4.2) appear with coefficients that are either independent of r, or are bounded by $c_0 r^{-1/2} e^{-T(p)}$. By virtue of Lemma 5.1, the pull-back of $\psi$ is uniformly bounded. Meanwhile, Proposition 5.2 guarantees that the curvature of the pull-back of A is also uniformly bounded, and thus so is A after an appropriate gauge transformation. This being the case, standard elliptic regularity arguments guarantee the following: Fix $R \geq 1$, $\varepsilon > 0$ and $k \in \{0, 1, 2, \ldots\}$ and there exists $r(R, \varepsilon, k) > 1$ that is independent of p, r and $(A, \psi)$ and has the following significance: If $r > r(R, \varepsilon, k)$, the pull back of $(A, \psi)$ via the composition of the Gaussian coordinate chart map and the rescaling map has $C^k$ distance $\varepsilon$ or less in the radius R ball centered at the origin in $\mathbb{R}^4$ to a pair, $(A_0, \psi_0 = (\alpha_0, 0))$, defined on all of $\mathbb{R}^4$, and obeying:

- $\bar{\partial}_0 \alpha_0 = 0$.
- $|\alpha_0| \leq 1$.
- $P^+ F_0 = \frac{1}{2}(1 - |\alpha_0|^2) \omega_0$.
- $|P^- F_0| \leq c$.

(5.2)

Here, $\mathbb{R}^4$ is written as $\mathbb{C}^2$, the connection $A_0$ is a type 1-1 connection, and $\bar{\partial}_0$ is the d-bar operator that is defined by $A_0$. Meanwhile, $F_0$ is the curvature 2-form of $A_0$, $P^\pm$ are the self-dual and anti-self dual projections for $\wedge^2 T^*\mathbb{R}^4$ with respect to the standard metric, and $\omega_0 = dx^1 \wedge dx^2 + dx^3 \wedge dx^4$. The constant c that appears in the final item of (5.2) is twice the constant $\kappa_*$ from Proposition 5.2.

Note that (5.2) does not assume the conclusions of Proposition 5.3. When Proposition 5.3 is assumed, then it is also the case that

$$1 - \delta \leq |\alpha_0| .$$
(5.3)

Granted that the pull-back of $(A, \psi)$ by the rescaling map is close to $(A_0, \psi_0)$, then Proposition 4.8 follows with a proof that $|\alpha_0| = 1$ and that $\alpha_0$ is $A_0$-covariantly constant when both (5.2) and (5.3) hold. Indeed, these conclusions are valid, then the covariant derivatives of the rescaled $\alpha$ at the origin are small. This last statement gives what is required by Proposition 4.8 by rescaling to the original Gaussian coordinates.

Given what was just said, the next lemma supplies what is needed to finish the proof of Proposition 4.8.

**Lemma 5.4**: *Suppose that $\delta < 1$ and $(A_0, \alpha_0)$ obeys (5.2) and (5.3). Then $|\alpha_0| = 1$ and $\alpha_0$ is $A_0$-covariantly constant.*

***Proof of Lemma 5.4***: Note that as $\alpha_0$ is nowhere zero, it is gauge equivalent on $\mathbb{R}^4$ to a map to $\mathbb{C}$ that is everywhere real and positive. Thus, $e^{-u}$ for some function u. Note that $u \geq 0$ since $|\alpha_0| \leq 1$. By virtue of the third equation in (5.2), the function u obeys

$$\Delta u = 1 - e^{-2u} ,$$
(5.4)

where $\Delta = \sum_{j=1,\ldots,4} \frac{\partial^2}{\partial x_j^2}$ is the Laplacian. By assumption, u is non-negative and bounded. Note that the maximum principle implies that $u > 0$ or u is everywhere zero. By virtue (5.3), the function u is bounded, and so there exists $\varepsilon > 0$ such that $1 - e^{-2u} - \varepsilon u > 0$. Thus,

$$\Delta u - \varepsilon u \geq 0 .$$
(5.5)

The only bounded solution to this last equation is $u = 0$. To see that such is the case, introduce $G_x$ to denote the Green's function on $\mathbb{R}^4$ with pole at x for the operator $-\Delta + \varepsilon$. Then $G_x$ is positive, and it decays to 0 on $\mathbb{R}^4$ as fast as $\exp(-\frac{1}{2} \varepsilon^{1/2}|\cdot|)$. Multiply both sides of (5.5) by $G_x$ to see that $-u(x) \geq 0$.

### b) Proof of Lemma 5.1

The starting point for the analysis is the second order equation for $\psi$ that is obtained from the second equation in (4.2):

$$-\frac{\partial^2}{\partial s^2}\psi + \nabla^\dagger\nabla\psi - cl(\tfrac{\partial}{\partial s}A + B_A)\psi + R\psi = 0.$$

(5.6)

Here, the covariant derivative $\nabla$ differentiates only along directions tangent to M and it is defined at a given $s \in \mathbb{R}$ by the connection $A|_s$. In (5.6), $R$ denotes an endomorphism of $\mathbb{S}$ that is independent of $r$ and $(A, \psi)$ that is uniformly bounded on $\mathbb{R} \times M$.

What with the first equation in (4.2), the equation in (5.6) implies that

$$-\tfrac{1}{2}\tfrac{\partial^2}{\partial s^2}|\psi|^2 + \tfrac{1}{2}d^\dagger d|\psi|^2 + |\tfrac{\partial}{\partial s}\psi|^2 + |\nabla\psi|^2 + re^{2T}(|\psi|^4 + i\psi^\dagger cl(a)\psi) + \psi^\dagger\mathfrak{r}\psi = 0.$$

(5.7)

Here, $\mathfrak{r}$ also denotes an endomorphism of $\mathbb{S}$ that is independent of $r$ and $(A, \psi)$, and is uniformly bounded on $\mathbb{R} \times M$. Set $c_* = \sup_{\mathbb{R} \times M}|\mathfrak{r}|$, set $z = |\psi|^2$ and set $x = r^{-1}e^{-2T}$. Then (5.6) implies that

$$-\tfrac{1}{2}\tfrac{\partial^2}{\partial s^2}z + \tfrac{1}{2}d^\dagger dz + re^{2T}z(z - 1 - c_* x) \leq 0.$$

(5.8)

As $\tfrac{d^2}{ds^2}x = -2T''x + 4(T')^2 x$, so (4.20) implies the following: For any $\sigma > c_R$,

$$-\tfrac{1}{2}\tfrac{\partial^2}{\partial s^2}(z-1-\sigma x) + \tfrac{1}{2}d^\dagger d(z-1-\sigma x) + re^{2T}z(z-1-\sigma x) \leq -re^{2T}z(\sigma - c_*)x + (4(T')^2 - 2T'')x.$$

(5.9)

Granted the preceding, suppose that $\sigma \geq 1 + c_*$ and that $r \geq \sup(4(T')^2 - 2T'')$. Then (5.9) precludes a positive, local maximum on $\mathbb{R} \times M$ for the function $(z-1-\sigma x)$. By virtue of Lemma 1.6 and the third item in (4.2), there exists $c_0$ such that if $\sigma$ is also greater than $c_0$, then $(z-1-\sigma x)$ is negative as $s \to -\infty$. Granted the fourth item in (4.2), this function also limits to zero as $s \to +\infty$. Because it can not have a positive local maximum, it follows that $z < 1 + \sigma x$ everywhere on $\mathbb{R} \times M$. This last bound gives the first item in Lemma 5.1.

To obtain the second item in Lemma 5.1, project (5.6) onto the respective $I_\mathbb{C}$ and $K^{-1}$ summands in $\mathbb{S}$ to obtain second order equations for $|\alpha|^2$ and $|\beta|^2$. These look much like what is written in (6.4) of [T1] (see (2.3) and (2.4) of [T3]) as they have the form

- $-\tfrac{1}{2}\tfrac{\partial^2}{\partial s^2}|\alpha|^2 + \tfrac{1}{2}d^\dagger d|\alpha|^2 + |\tfrac{\partial}{\partial s}\alpha|^2 + |\nabla\alpha|^2 - re^{2T}(1 - |\alpha|^2 - |\beta|^2)|\alpha|^2$
  $+ \mathfrak{r}_0(\alpha,\beta) + \mathfrak{r}_1(\alpha,\nabla'\beta) + \mathfrak{r}_2|\alpha|^2 = 0.$

- $-\frac{1}{2} \frac{\partial^2}{\partial s^2} |\beta|^2 + \frac{1}{2} d^\dagger d |\beta|^2 + |\frac{\partial}{\partial s} \beta|^2 + |\nabla' \beta|^2 + r e^{2T}(1 + |\alpha|^2 + |\beta|^2)|\beta|^2$
$$+ \mathfrak{r}'_0 |\beta|^2 + \mathfrak{r}'_1 (\beta, \nabla \alpha) + \mathfrak{r}'_2 (\alpha, \beta) = 0.$$
(5.10)

Here, $\mathfrak{r}_2$ and $\mathfrak{r}'_2$ have support where $s \leq 2$. Introduce $w = (1 - |\alpha|^2)$, and note that the top equation in (5.10) implies that $w$ obeys the inequality

$$-\frac{1}{2} \frac{\partial^2}{\partial s^2} w + \frac{1}{2} d^\dagger d w + r e^{2T} w - (|\frac{\partial}{\partial s} \alpha|^2 + |\nabla \alpha|^2 + r e^{2T}(w^2 + |\beta|^2 |\alpha|^2)) + (|\alpha|^2 + |\nabla' \beta|^2) \geq 0.$$
(5.11)

Arguing as in the derivation of (5.9), this last equation plus the lower equation in (5.8) implies the following: There are constants $c_0$ and $c_1$ that are independent of $r$ and $(A, \psi)$ and are such that the function $u = |\beta|^2 - c_0 r^{-1} e^{-2T} w - c_1 r^{-2} e^{-4T}$ obeys

$$-\frac{1}{2} \frac{\partial^2}{\partial s^2} u + \frac{1}{2} d^\dagger d u + \frac{1}{2} r e^{2T} u \leq 0.$$
(5.12)

The argument is essentially that given to prove Proposition 2.3 in [T3].

To apply this last equation, fix a smooth, non-negative function, $\sigma$, on $\mathbb{R}$ that is equal to 1 on $[-1, 1]$ and vanishes where $|s| \geq 2$. For $R \geq 1$, define $\sigma_R$ in $\mathbb{R} \times M$ by the rule $\sigma_R(s) = \sigma(R^{-1} s)$. Let $u_+ \colon \mathbb{R} \times M \to [0, \infty)$ denote the maximum of $u$ and 0. Note that the apriori bound on $|\psi|$ given previously implies that $u_+ \leq 2$ if $r$ is large. Multiply both sides of (5.12) by $\sigma_R$ and integrate over the support of $u_+$. Integration by parts finds that

$$-\frac{1}{2} R^{-2} \int_{\mathbb{R} \times M} (\frac{d^2}{ds^2} \sigma)_R u_+ + \frac{1}{2} r \int_{\mathbb{R} \times M} \sigma_R e^{2T} u_+ \leq 0.$$
(5.13)

Now, the integrand of the left most integral in (5.13) is supported where $R \leq |s| \leq 2R$, and its integrand is bounded by $2 \sup_{\mathbb{R}} |\frac{d^2}{ds^2} \sigma|$. As a consequence, the left most integral in (5.13) is bounded in absolute value by $c_0 R^{-1}$. Taking $R$ ever larger proves that $u_+$ must vanish. This then proves the second assertion of the lemma.

### c) Norm bounds on $[-6, \infty) \times M$

The first step to a proof of Proposition 5.2 gives weighted $L^2$ bounds on $\mathbb{R} \times M$ for $(1 - |\alpha|^2)$, $\beta$, and covariant derivatives of $\alpha$ and $\beta$.

**Lemma 5.5**: *There exists a constant $\kappa \geq 1$ with the following significance: Suppose that $r \geq \kappa$ and that $(A, \psi)$ is a solution to (4.2). Then*

- $\int_{\mathbb{R} \times M} e^{2s} (|\frac{\partial}{\partial s} \alpha|^2 + |\nabla \alpha|^2 + r e^{2T} (1 - |\alpha|^2)^2) \leq \kappa r^{-1}$,

- $\int_{\mathbb{R}\times M} e^{2s}(|\frac{\partial}{\partial s}\beta|^2 + |\nabla'\beta|^2) + \frac{1}{4}re^{2T}|\beta|^2) \le \kappa r^{-1}$

***Proof of Lemma 5.5***: The proof given here modifies an argument from [KM2]. To start, multipy both sides of the lower equation in (5.10) by $e^{2s}$ and integrate over $\mathbb{R}\times M$. Lemma 4.6 guarantees that this integral is absolutely convergent, and that it is permissable to integrate by parts so as to obtain an identity of the form

$$\int_{\mathbb{R}\times M} e^{2s}(|\tfrac{\partial}{\partial s}\beta|^2 + |\nabla'\beta|^2) + \tfrac{1}{4}re^{2T}|\beta|^2) \le r^{-1}(1 + \int_{\mathbb{R}\times M} e^{2s-2T}|\nabla\alpha|^2)$$

(5.14)

when r is larger than some $(A, \psi)$ independent constant. Save this equation.

Introduce the differential form $\omega = e^{2s}(ds \wedge a + *a)$. This is a closed, exact self-dual 2-form on $\mathbb{R}\times M$; it is the exterior derivative of the 1-form $\tfrac{1}{2}e^{2s}a$. Also introduce $F_A$ to denote the 2-form $ds \wedge \tfrac{\partial}{\partial s}A + *B_A$ on $\mathbb{R}\times M$. This is also a closed 2-form on $\mathbb{R}\times M$; it is the curvature 2-form of the connection A when the latter is viewed as a connection on the trivial line bundle over $\mathbb{R}\times M$. Let $\partial_A$ denote the covariant exterior derivative on this trivial line bundle as defined by the connection A on $\mathbb{R}\times M$.

Now define the 2-form

$$\wp = \partial_A\alpha \wedge \partial_A\bar{\alpha} + (1 - |\alpha|^2)F_A.$$

(5.15)

This 2-form is closed. Moreover, it follows from Lemma 5.1 and from Lemma 4.6 that $e^{kT}|\wp|$ is bounded on $\mathbb{R}\times M$ for any $k > 0$. This being the case, the 4-form $\omega \wedge \wp$ is integrable over $\mathbb{R}\times M$. It also follows that the integral is zero since $\wp$ is closed and $\omega$ is the exterior derivative of $\tfrac{1}{2}e^{2s}a$.

To see what to make of this, introduce $\partial_A = \tfrac{1}{2}(\partial_A + ie^{-2s}*(\omega \wedge \partial_A))$; this operator sends a given section over $\mathbb{R}\times M$ of the relevant line bundle to one of the bundle's tensor product with $T_\mathbb{C}^*(\mathbb{R}\times M)$. The conjugate, $\bar{\partial}_A = \tfrac{1}{2}(\partial_A - ie^{-2s}*(\omega \wedge \partial_A))$, is the d-bar operator for the almost complex structure that is defined using the product metric and the self-dual 2-form $e^{-2s}\omega$. Granted this notation, then the vanishing of the integral of $\omega \wedge \wp$ is equivalent to the assertion that

$$\int_{\mathbb{R}\times M} e^{2s}(|\partial_A\alpha|^2 - |\bar{\partial}_A\alpha|^2 + re^{2T}((1 - |\alpha|^2)^2 + (1 - |\alpha|^2)|\beta|^2)) =$$
$$\int_{\mathbb{R}\times M} e^{2s}(1 - |\alpha|^2)\langle a, 2\chi\mu + \hat{\mu}\rangle.$$

(5.16)

Here, the notation $\langle , \rangle$ denotes the inner product that is defined by the metric on M. Since the second equation in (4.2) writes $\bar{\partial}_A \alpha$ in terms of $\partial_A \beta$, and since $\hat{\mu}$ has compact support, this last equation with (5.14) implies the first of Lemma 5.5's inequalities. The first inequality in Lemma 5.5 with (5.14) implies the second.

### d) $L^2_1$ bounds for $B_A$ over $[-6, 6] \times M$

The next step towards a proof of Propositions 5.2 establishes a a preliminary bound for the $L^2$ norm of $B_A$ over the portion of $\mathbb{R} \times M$ where $-6 \leq s \leq 6$. Here is the required $L^2$ bound:

**Lemma 5.6**: *There exists a constant $\kappa > 1$ such that if $r \geq \kappa$, and if $(A, \psi)$ is a solution to (4.2), then $\int_{[-6,6] \times M} |B_A|^2 \leq \kappa r$.*

***Proof of Lemma 5.6:*** Let $\{e_\nu\} \subset C^\infty(M; T^*M)$ denote a complete, $L^2$-orthonormal basis of eigenvectors of $*d$ on the subspace of co-closed 1-forms. Thus, $*de_\nu = \lambda_\nu e_\nu$ where $\lambda_\nu \in \mathbb{R}$. For each index $\nu$, introduce the function $s \to b_\nu(s)$ on $\mathbb{R}$ by the rule

$$b_\nu = \int_{\{s\} \times M} e_\nu \wedge *B_A .$$

(5.17)

The top equation in (4.2) finds this function obeying the equation

$$\tfrac{d}{ds} b_\nu + \lambda_\nu b_\nu = \mathfrak{w}_\nu ,$$

(5.18)

where $\mathfrak{w}_\nu|_s$ is defined to be the inner product on $\{s\} \times M$ between $e_\nu$ and the section of $iT^*M$ over $\mathbb{R} \times M$ that is defined by $\mathfrak{m} = *d[re^{2T}(\psi^\dagger \tau^k \psi - i a) + i\chi(*d\mu + \varpi_K) + i\hat{\mu}]$. The equation in (5.18) can be integrated so as to give

- $b_\nu(s) = -e^{-\lambda_\nu s} \int_s^\infty e^{\lambda_\nu t} \mathfrak{w}_\nu(t)\, dt$ *when $\lambda_\nu < -1$.*
- $b_\nu(s) = e^{-\lambda_\nu s} \int_{-\infty}^s e^{\lambda_\nu t} \mathfrak{w}_\nu(t)\, dt$ *when $\lambda_\nu > 1$.*

(5.19)

It follows from Lemma 5.5 that these integrals are well defined; this is because

$$|\mathfrak{m} + 2ir e^{2T} a| \leq c_0 (1 + r e^{2T}(|\alpha||\beta| + |\beta|^2 + |\nabla \alpha| + |\nabla \beta|)) .$$

(5.20)

In particular, (5.19) implies that

- $|b_v(s)|^2 \leq |\lambda_v + 1|^{-1} \int_s^\infty e^{-2t} |\mathfrak{w}_v(t)|^2 \, dt$ when $\lambda_v < -1$,
- $|b_v(s)|^2 \leq |\lambda_v - 1|^{-1} \int_{-\infty}^s e^{2t} |\mathfrak{w}_v(t)|^2 \, dt$ when $\lambda_v > 1$;

(5.21)

and so (5.20), (5.21) with Lemma 5.1 give the bound

$$\sum_v{}' |b_v(s)|^2 \leq c_0 r \quad \text{when } s \in [-6, 6].$$

(5.22)

Here, the prime on the summation sign is meant to indicate that the sum is restricted to the eigenvectors whose eigenvalue has absolute value greater than 1.

To obtain bounds for the remainder of $B_A$, introduce now $d_A$ to denote the covariant, exterior derivative along the constant s slices of $\mathbb{R} \times M$ as defined by the connection A. Thus, $*B_A |\alpha|^2 = \bar{\alpha} d_A \nabla_A \alpha$. Granted this, write

$$b_v(s) = \int_{\{s\} \times M} e_v \wedge *B_A = \int_{\{s\} \times M} e_v \wedge *B_A(1 - |\alpha|^2) + \int_{\{s\} \times M} e_v \wedge \bar{\alpha} d_A \nabla_A \alpha \,.$$

(5.23)

Now integrate by parts in the right most integral to obtain the bound

$$|b_v(s)| \leq c_v (\|B_A|_s\|_2 \|(1 - |\alpha|^2)|_s\|_2 + \|\nabla\alpha|_s\|_2 + \|\nabla\alpha|_s\|_2^2)$$

(5.24)

Here, $c_v$ depends on the eigenvector $e_v$. However, there is an upper bound for all such $c_v$ for which $|\lambda_v| \leq 1$.

Meanwhile, by virtue of Lemma 5.5, there exist $s \in [-6, 6]$ where $\|(1 - |\alpha|^2)|_s\|_2 \leq c_0 r^{-1}$ and $\|\nabla\alpha|_s\|_2 \leq c_0 r^{-1/2}$. For such s,

$$|b_v(s)| \leq c_0 r^{-1} \|B_A|_s\|_2 + c_0 r^{-1/2} \quad \text{when } |\lambda_v| \leq 1.$$

(5.25)

This last estimate plus (5.22) finds $s \in [-6, 6]$ where $\|B_A|_s\|_2 \leq c_0 r^{1/2}$ and

$$|b_v(s)| \leq c_0 r^{-1/2} \quad \text{when } |\lambda_v| \leq 1.$$

(5.26)

To finish the story, fix some point in [-6, 6] where (5.26) holds, and integrate (5.17) from that point using (5.26) and Lemma 5.5 to conclude that

$$|b_v(s)| \leq c_0 (1 + r^{-1}(\int_{[-6,6] \times M} |B_A|^2)^{1/2}) \,,$$

(5.27)

at all values of $s \in [-6, 6]$. Equations (5.27) and (5.22) imply the lemma's claim.

The $L^2$ bound for $B_A$ over the cylinder $[-6, 6] \times M$ can be coupled with what is said in Lemma 5.6 to obtain an $L^2_1$ bound for $B_A$ over $[-5, 5] \times M$.

**Lemma 5.7**: *There exists $\kappa \geq 1$ such that if $r \geq \kappa$ and if $(A, \psi)$ is a solution to (4.2), then*

$$\int_{[-5,5] \times M} (|\tfrac{\partial}{\partial s} B_A|^2 + |\nabla B_A|^2) \leq \kappa r.$$

*Proof of Lemma 5.7*: Let $\sigma$ here denote a smooth, non-negative function with compact support on $[-6, 6]$ with value 1 on $[-5, 5]$. As $d_A * B_A = 0$, so (5.17) implies that

$$\int_{[-6,6] \times M} (\sigma^2 |\tfrac{\partial}{\partial s} B_A|^2 + \sigma^2 |dB_A|^2 + \sigma^2 |d*B_A|^2 - 2\sigma'\sigma \langle B_A, *dB_A \rangle) = \int_{[-6,6] \times M} \sigma^2 |\mathfrak{m}|^2.$$
(5.28)

This plus the bounds in Lemma 5.5 for the $L^2$ norm of $\mathfrak{m}$ and the bound in Lemma 5.6 for the $L^2$ norm of $B_A$ imply what is claimed by the lemma.

### e) Proof of Proposition 5.2's assertion for points in $[-2, \infty) \times M$

The assertions of the proposition on $[-2, \infty) \times M$ follow from

**Lemma 5.8**: *There exists $\kappa > 1$ such that if $r > \kappa$ and if $(A, \psi)$ is a solution to (4.2), then*

$$|\tfrac{\partial}{\partial s} A + B_A| \leq r e^{2T}(1 - |\alpha|^2) + \kappa \quad \text{and} \quad |\tfrac{\partial}{\partial s} A - B_A| \leq r e^{2T}(1 - |\alpha|^2) + \kappa$$

*at all points where $s \geq -2$.*

*Proof of Lemma 5.8*: The bound for $|\tfrac{\partial}{\partial s} A + B_A|$ follows directly from what is written on the top line of (4.2). The proof of the asserted bound for $|\tfrac{\partial}{\partial s} A - B_A|$ is a modification of an argument from [KM2]. The proof for this bound is given here in four steps.

<u>Step 1</u>: It proves useful to introduce as notation $F_- = \tfrac{\partial}{\partial s} A - B_A$ and $F_+ = \tfrac{\partial}{\partial s} A + B_A$. Both $F_-$ and $F_+$ are sections over $\mathbb{R} \times M$ of $iT^*M$; they are obtained from the 2-form $F_A = ds \wedge \tfrac{\partial}{\partial s} A + *B_A$ by contracting the respective anti-self dual and self-dual parts of the latter form with the vector field $\tfrac{\partial}{\partial s}$. Note in this regard that $F_A$ is the curvature 2-form for the connection A when the latter is viewed as a connection for the trivial bundle over $\mathbb{R} \times M$. The fact that $F_A$ is a closed 2-form is equivalent to the pair of equations

$$\tfrac{\partial}{\partial s} F_- + *dF_- = \tfrac{\partial}{\partial s} F_+ - *dF_+ \quad \text{and} \quad d*F_- = d*F_+.$$
(5.29)

Here, d again denotes the exterior derivative along the constant s slices of $\mathbb{R} \times M$. Note that the right hand sides of both of the equations in (5.29) can be written in terms of $\psi$ using the first two lines in (4.2). Doing so, and then differentiating the preceding equations using appropriate combinations of partial dervatives gives a second order equation for $e^{-2T} F_-$. The latter equation implies a differential inequality for the function $q_- = e^{-2T}|\frac{\partial}{\partial s} A - B_A|$ that has the form

$$(-\tfrac{\partial^2}{\partial s^2} + d^\dagger d) q_- + 2r e^{2T}|\alpha|^2 q_- \leq c_0 q_- + r \mathfrak{z} + 2r(|\tfrac{\partial}{\partial s}\alpha|^2 + |\nabla\alpha|^2 + |\tfrac{\partial}{\partial s}\beta|^2 + |\nabla'\beta|^2).$$

(5.30)

Here $\mathfrak{z}$ has support where $s \leq c_0$; and $\mathfrak{z} \leq c_0$ in any event. This equation implies the following:

**Lemma 5.9**: *There exist $\kappa > 1$ and positive constants $\kappa_1$, $\kappa_2$ and $\kappa_3$ such that if $r \geq \kappa$ and if $(A, \psi)$ is a solution to (4.2), then $q_0 = \max(q_- + \kappa_1 r|\beta|^2 - (r + \kappa_2) w - \kappa_3 e^{-2T}, 0)$ obeys*

$$(-\tfrac{\partial^2}{\partial s^2} + d^\dagger d) q_0 + 2r e^{2T}|\alpha|^2 q_0 \leq c_0 (q_0 + r w^2)$$

(5.31)

*Proof of Lemma 5.9*: By virtue of (5.11), $x_1 = q_- - (r + c) w$ obeys

$$(-\tfrac{\partial^2}{\partial s^2} + d^\dagger d) x_1 + 2r e^{2T}|\alpha|^2 x_1 \leq c_0 x_1 + c_0(r + c) w_+ + r c_0 + 3r(|\tfrac{\partial}{\partial s}\beta|^2 + |\nabla'\beta|^2) - c|\nabla\alpha|^2.$$

(5.32)

Here, $w_+ = \max(w, 0)$. Set $x_2 = x_1 + 4r|\beta|^2$; given the bottom equation in (5.10), $x_2$ obeys

$$(-\tfrac{\partial^2}{\partial s^2} + d^\dagger d) x_2 + 2r e^{2T}|\alpha|^2 x_2 \leq c_0 x_2 + c_0(r + c) w_+ + r c_0$$

(5.33)

when $c \geq c_0$ and $r \geq c_0$. Now let $x_3 = x_2 - e^{-2T} c'$. Then this function obeys

$$(-\tfrac{\partial^2}{\partial s^2} + d^\dagger d) x_3 + 2r e^{2T}|\alpha|^2 x_3 \leq c_0 (x_3 + r w_+) + r c_0 (1 + c'/r) - 2r|\alpha|^2 c'.$$

(5.34)

when $r \geq c_0$. Now, take $c' > c_0 (1 + c'/r)$ so that the right hand side of (5.34) is no greater than

$$c_0 (x_3 + r w_+) + r c_0 (1 + c'/r)(1 - |\alpha|^2) - r |\alpha|^2 c' \leq c_0 x_3 + r c_0 (2 + c'/r) w_+ - r|\alpha|^2 c'.$$

(5.35)

If it is also the case that $c' \geq c_0(2 + c'/r)$, then the right most two terms in (5.35) are no greater than $r c_0(2 + c'/r)(w_+ - |\alpha|^2)$. Since $w_+ - |\alpha|^2 \leq w_+^2 \leq w^2$, the claim follows.

Step 2: This step bounds $q_0$ on the cylinder $[-4, 4] \times M$.

**Lemma 5.10**: *There exists a constant, $\kappa > 1$, such that when $r \geq \kappa$ and $(A, \psi)$ is a solution to (4.2), then $q_0 \leq \kappa\, r^{1/2}$ at points in $[-4, 4] \times M$.*

***Proof of Lemma 5.10***: Fix a smooth function, $\chi\colon [0, \infty) \to [0, 1]$ that equals 1 on $[0, \frac{1}{4}]$ and vanishes on $[\frac{1}{2}, \infty)$. For $p \in \mathbb{R} \times M$, set $\chi_p = \chi(\mathrm{dist}(\cdot, p))$. Meanwhile let $G_p$ denote the Green's function for the operator $(-\frac{\partial^2}{\partial s^2} + d^\dagger d)$ with pole at $p$. Note in particular that $|G_p| \leq c_0\, \mathrm{dist}(\cdot, p)^{-2}$. Take $p \in [-4, 4] \times M$, multiply both sides of (5.31) by $\chi_p G_p$ and integrate. A judicious integration by parts finds that

$$q_0|_p \leq c_0\, r^{1/2} + c_0 \int_{[-5,5] \times M} \mathrm{dist}(p, \cdot)^{-2}(q_0 + r w^2)\ .$$

(5.36)

Here, Lemma 5.6 has been invoked to bound the $L^2$ norm of $q_0$ on the support of $d\chi_p$. Note next that it follows from Lemmas 5.6 and 5.7 that the $L^2$ norm over $[-5, 5] \times M$ of the function $\mathrm{dist}(p, \cdot)^{-1} q_-$ is finite, and bounded by $c_0\, r^{1/2}$. As a consequence of this and Lemma 5.5, the size of the integral on the right hand side of (5.36) is less than $c_0 r^{1/2}$.

Step 3: Fix a smooth function, $\sigma\colon \mathbb{R} \to [0, 1]$ that vanishes where $s \leq -2.75$ and is equal to 1 where $s \geq -2.5$. The function $\sigma q_0$ has support where $s \geq -3$ and it obeys the inequality

$$(-\tfrac{\partial^2}{\partial s^2} + d*d)(\sigma q_0) + r e^{2T} |\alpha|^2 (\sigma q_0) \leq c_0 (\sigma q_0 + r \sigma w^2) - \sigma'' q_0 - 2\sigma' \tfrac{\partial}{\partial s} q_0$$

(5.37)

when $r \geq c_0$. Equation (5.37) is used to bound the $L^1$ norm of $e^{2T}\sigma q_0$; this done as follows: First, introduce $z$ to denote the supremum of $r^{-1}\sigma q_0$. Now, divide both sides of (5.37) by $r$ and integrate to find that

$$\int_{\mathbb{R}\times M} e^{2T}\sigma q_0 \leq c_0 \int_{\mathbb{R}\times M} e^{2T} \sigma q_0 w_+ + c_0 \left( \int_{[-3,\infty)\times M} w^2 + r^{-1/2} \right)\ .$$

(5.38)

Here, Lemma 5.10 is used to obtain a bound of $c_0 r^{1/2}$ for the absolute value of the integral of $q_0$ over the region where $\sigma''$ is non-zero. Equation (5.38) with Lemma 5.5 implies that

$$\int_{\mathbb{R}\times M} e^{2T}\sigma q_0 \leq c_0 \sup(|\sigma q_0|^{1/2})\, \Big( \int_{\mathbb{R}\times M} e^{2T}\sigma q_0 \Big)^{1/2} r^{-1} + c_0\, r^{-1/2}\ .$$

(5.39)

Since $\sup(\sigma q_0) = r z$, this then implies that

$$\int_{\mathbb{R}\times M} e^{2T}\sigma q_0 \leq c_0(r^{-1}z + r^{-1/2}) \;.$$

(5.40)

Step 4: Let $p \in [-2, \infty) \times M$. Multiply both sides of (5.37) by $\chi_p G_p$ and integrate. Since $\sigma = 1$ on the domain of $\chi_p$, what with Lemma 5.5 and (5.40), the resulting equation implies that

$$\sigma q_0|_p \leq c_0 (zr \int_{\text{dist}(\cdot,p)\leq 1} \text{dist}(\cdot,p)^{-2} w_+ + r \int_{\text{dist}(\cdot,p)\leq 1} \text{dist}(\cdot,p)^{-2} w^2) \leq c_0 e^{-2T}(z(1+\ln r) + 1)$$

(5.41)

when $r \geq c_0$. A similar argument for the case when $p \in [-3, -1] \times M$ using also Lemma 5.10 finds that

$$\sigma q_0|_p \leq c_0(z(1 + \ln(r)) + r^{1/2}).$$

(5.42)

Together, (5.41) and (5.42) imply that $z \leq c_0 r^{1/2}$. Granted that such is the case, then (5.41) implies the assertion made in Proposition 5.2 where $s \geq -2$.

### f) Proof of Proposition 5.3 at points in $[-1, \infty) \times M$

Suppose that $p \in \mathbb{R} \times M$, that $c \geq 0$ is a given constant, and that $|\frac{\partial}{\partial s} A| + |B_A| \leq cr$ on the ball of radius $\frac{1}{2}$ centered at p. Rescale a Gaussian coordinate system centered at p so that the ball of radius $r^{-1/2} e^{-T}$ is mapped to the ball of radius 1, centered at the origin in $\mathbb{R}^4$. As argued just prior to (5.2), given $R \geq 1$, $\varepsilon > 0$ and $k \in \{0, 1, 2, \ldots\}$, there exists $r(R, \varepsilon, k, c) > 1$ that is independent of p, r and $(A, \psi)$ and has the following significance: If $r > r(R, \varepsilon, k, c)$, then the pull back of $(A, \psi)$ via the composition of the Gaussian coordinate chart map and the rescaling map will have $C^k$ distance $\varepsilon$ or less in the radius R ball centered at the origin in $\mathbb{R}^4$ to a pair, $(A_0, \psi_0 = (\alpha_0, 0))$, that is defined on all of $\mathbb{R}^4$ and obeys (5.2). As the norm of the curvature of $A_0$ is uniformly bounded, the absolute value of the $A_0$-covariant derivative of $\alpha_0$ is bounded by some c-dependent constant, $c_0(c)$. Rescaling back to the original Gaussian coordinates gives $|\nabla \alpha| \leq c_0(c) r^{1/2} e^{T(p)}$ at p.

Now, suppose that $\delta > 0$, that $p \in [-1, \infty) \times M$, and that $1 - |\alpha|^2 > \delta$ at p. As Proposition 5.2 holds on $[-2, \infty) \times M$, the preceding bound for $|\nabla \alpha|$ holds near p, and so it follows that there exits $c_0 > 1$ that is independent of p, r, and $(A, \psi)$ such that when $r \geq c_0$, then $1 - |\alpha|^2 > \frac{1}{2}\delta$ in the ball of radius $c_0 r^{-1/2} e^{-T(p)}$ centered at p. This implies that the ball of radius $r^{-1/2} e^{-T(p)}$ centered at p contributes $c_0^{-1} \delta^2 r^{-1}$ to the first integral in Lemma 5.5. As a consequence, there is a set of at most $N_\delta = c_0 \delta^{-2}$ disjoint balls in $[-1.5, \infty) \times M$ with the following two properties: Let U denote any given ball from this set and p its center point. Then U has radius $r^{-1/2} e^{-T(p)}$. Second, $(1 - |\alpha|^2)$ is greater than $\delta$ at p.

To see that there are no points in $[-1.5, \infty) \times M$ where $1 - |\alpha|^2 > \delta$ when r is large, suppose to the contrary $p \in [-1.5, \infty) \times M$ is such a point. Fix $R \geq 4N_\delta$ and $\varepsilon \in (0, \frac{1}{100}\delta)$. Suppose that $r \geq r(R, \varepsilon, 0)$. Pull back $(A, \psi)$ as described in the opening paragraph using a Gaussian coordinate system centered at p. Then the pull-back of $\alpha$ is withing $\varepsilon$ of $\alpha_0$ where $(A_0, \alpha_0)$ obeys (5.2) on $\mathbb{R}^4$. Moreover, $1 - |\alpha_0|^2 > \frac{3}{4}\delta$ at the origin but $1 - |\alpha_0|^2 < \frac{1}{4}\delta$ on a sphere of radius 1 or more centered at the origin. As explained next, this is not possible.

To see why no such $(A_0, \alpha_0)$ exists, suppose first that $|\alpha_0| > 0$ in the interior of the sphere in question. Then there is a gauge transformation that renders $\alpha_0$ as $e^{-u}$ with $u \geq 0$ a smooth solution to (5.4) on the interior of the sphere. Note that u is smaller at all points on the boundary than it is at the origin. This being the case, u must have a local maximum in the interior of the ball. However, the maximum principle with (5.4) rules out such a point. Meanwhile, if $\alpha_0 = 0$ in the ball, then its zero locus is non-empty, and so a complex analytic curve. In particular, the zero locus can't be compact and so must intersect the boundary of the ball.

**g) The behavior of $\mathfrak{a}(A, \psi)$ on $\{s\} \times M$ for $s \in [-1, 0]$.**

The proof of Proposition 5.2 for points in $(-\infty, -2] \times M$ starts here by saying somewhat more about the $(A, \psi)$ on $[-2, 0] \times M$. The discussion that follows refers to the version of $\mathfrak{a}$ that appears in (1.4). With $s \in (-\infty, 0)$ given, set $\mathfrak{a}|_s$ to equal $\mathfrak{a}(A(s), \psi(s))$ when $(A, \psi)$ is a solution to (4.2).

**Lemma 5.11**: *There exists a constant, $\kappa \geq 1$, with the following significance: Suppose that $r \geq \kappa$ and that $(A, \psi)$ is a solution to (4.2). Then $\int_{-1}^{0} |\mathfrak{a}|_s| ds \leq \kappa r^{1/2}$.*

*Proof of Lemma 5.11*: The conclusions of Proposition 5.3 where $s \geq -1$ assert that the section $\alpha$ of $I_\mathbb{C}$ has norm very near 1 on the whole of $[-1, \infty) \times M$ when $r > c_0$. As a consequence, there is a smooth map, u: $[-1, \infty) \times M \to S^1$ such that $u\alpha = |\alpha| 1_\mathbb{C}$. Meanwhile, $A - u^{-1}du$ can be written as $A_I + \hat{a}$ on $[-1, \infty) \times M$ with

$$|\hat{a}| \leq |\alpha|^{-1} |\nabla \alpha|.$$
(5.43)

It follows from the fourth point in (4.2) that u is homotopically trivial. Since the map u is null-homotopic, the function $\mathfrak{a}$ has the same value on $(A(s), \psi(s))$ as it does on the pair $(A_I + \hat{a}, (|\alpha|1_\mathbb{C}, u\beta))$. This understood, it follows from (5.43), Lemma 5.5 and Lemma 5.8 that $\int_{-1}^{0} |\mathfrak{cs}|_s| \leq c_0 r^{-1/2}$. It also follows from these results that $\int_{-1}^{0} |\mathfrak{e}_\mu|_s| \leq c_0$. Meanwhile, integrating by parts finds

$$\int_{[-1,0]\times M} |E|_s| \le 2\int_{[-1,0]\times M} |\hat{a}| \le c_0 r^{-1/2}.$$

(5.44)

Here, the right most inequality follows using (5.43) with Lemma 5.5. Thus, the term with $rE$ in (1.4) contributes at most $c_0 r^{1/2}$ to the integral of $|\mathfrak{a}|_s|$. Finally, the integrand for the right most term in (1.4) is bounded by the integral of $c_0 (|\beta| + |\nabla\alpha| + |\nabla\beta|)$. As a consequence, Lemmas 5.1 and 5.5 imply that the integral from -1 to 0 of the absolute value of right most term in (1.4) is also bounded by $c_0 r^{1/2}$.

### h) $L^2$ bounds for $\frac{\partial}{\partial s} A$ and $\frac{\partial}{\partial s}\psi$

The statement of these bounds uses

$$\mathfrak{B}_{(A,\psi)} = B_A - r(\psi^\dagger \tau^k \psi - i\,a) - i(*d\mu + \varpi_K).$$

(5.45)

The next lemma says something about the $L^2$ norms of $\frac{\partial}{\partial s} A$, $\mathfrak{B}_{(A,\psi)}$, $\frac{\partial}{\partial s}\psi$ and $D_A\psi$.

**Lemma 5.12**: *There exists $\kappa \ge 1$ such that if $r \ge \kappa$ and if $(A, \psi)$ is a solution to (4.2) with $s \to -\infty$ limit equal to $\mathfrak{c}$, then*

$$\int_{(-\infty,-1/2]\times M} \left(|\tfrac{\partial}{\partial s}A|^2 + |\mathfrak{B}_{(A,\psi)}|^2 + 2r(|\tfrac{\partial}{\partial s}\psi|^2 + |D_A\psi|^2)\right) \le \mathfrak{a}(\mathfrak{c}) + \kappa r^{1/2}.$$

*Proof of Lemma 5.12*: Fix $R \gg 1$ and a point $t \in (-R, 0)$. The top two equations in (4.2) imply the integral identity

$$\int_{[-R,t]\times M} \left(|\tfrac{\partial}{\partial s}A|^2 + |\mathfrak{B}_{(A,\psi)}|^2 + 2r(|\tfrac{\partial}{\partial s}\psi|^2 + |D_A\psi|^2)\right) = -\int_{-R}^{t} \tfrac{d}{ds}(\mathfrak{a}_s) = \mathfrak{a}_{-R} - \mathfrak{a}_t.$$

(5.46)

By virtue of the third point in (4.2), the $R \to \infty$ limit in (5.46) can be taken and doing so finds

$$\int_{(-\infty,t]\times M} \left(|\tfrac{\partial}{\partial s}A|^2 + |\mathfrak{B}_{(A,\psi)}|^2 + 2r(|\tfrac{\partial}{\partial s}\psi|^2 + |D_A\psi|^2)\right) \le \mathfrak{a}(\mathfrak{c}) - \mathfrak{a}|_t.$$

(5.47)

Meanwhile, Lemma 5.11 guarantees points in $[-\tfrac{1}{2}, 0]$ where $|\mathfrak{a}|_t| \le c_0 r^{1/2}$. Take $t$ in (5.47) to be such a point and the claim in Lemma 5.12 follows directly.

### i) Proof of Proposition 5.2 for points in $(-\infty, -1] \times M$

The first step here is to obtain a local $L^2$ bound for $B_A$. This is the part of the proof that assumes something about the $s \to -\infty$ limit of $\mathfrak{a}((A, \psi))$. To set things up, let $I \subset (-\infty, -1/2]$ denote a compact interval. The bound given in Lemma 5.12 is going to be used in what follows to derive $L^2$ bounds for $B_A$ and $\nabla \psi$ on $I \times M$. To do so, use the Bochner-Weitzenboch formula for the Dirac operator to write

$$\int_{\{s\} \times M} (|\mathfrak{B}_{(A,\psi)}|^2 + r|D_A\psi|^2) = \int_{\{s\} \times M} (|B_A|^2 + r^2|\psi^\dagger \tau \psi - i a|^2 + r(|\nabla \psi|^2 + \psi^\dagger R \psi))$$
$$+ 2i \int_{\{s\} \times M} (B_A - r(\psi^\dagger \tau \psi - ia)) \wedge d\mu - 2ir \int_{\{s\} \times M} (\psi^\dagger \tau \psi - ia) \wedge *\varpi_K - 2r\,\mathrm{E}(A(s)).$$

(5.48)

Here, $R$ is an endomorphism of $\mathbb{S}_I$ that is independent of $r$ and $(A, \psi)$. As can be seen from this last equation, bounds on the integral over $I \times M$ of $|B_A|^2$, $|\nabla \psi|^2$ and $|\psi^\dagger \tau \psi - ia|^2$ will follow from the bound given in Lemma 5.12 together with a suitable bound for $\mathrm{E}$. In any event, it follows from Lemma 5.12 and (5.48) that

$$\int_{[s,s+1] \times M} (|B_A|^2 + r^2|\psi^\dagger \tau \psi - ia|^2 + r(|\nabla \psi|^2 + \psi^\dagger R \psi)) \leq c_0(c_* + 1)r^2$$

(5.49)

when $\lim_{s \to -\infty} \mathfrak{a}((A, \psi)) \leq c_* r^2$. This last equation gives the desired $L^2$ bound for $B_A$. Note in this regard that it is always the case that such a bound exists for the $s \to -\infty$ limit of $\mathfrak{a}$. Indeed, this follows from (1.6) with Proposition 1.10.

To continue with the proof of Proposition 5.2, repeat the argument that led to Lemma 5.9 to deduce the existence of $r$ and $(A, \psi)$ independent constants $\kappa > 1$ and $\kappa_{1-3}$ such that when $r \geq \kappa$ and $(A, \psi)$ is a solution to (4.2), then

$$q_0 = \max(|\tfrac{\partial}{\partial s} A - B_A| + \kappa_1 r|\beta|^2 - (r + \kappa_2)w - \kappa_3, 0)$$

(5.50)

obeys the $T = 0$ version of (5.31) where $s \leq 0$. For the sake of argument, assume that $q_0 > 10r$ on $(-\infty, 0] \times M$. It then follows from Lemma 5.8 and the $s \geq -1$ version of Proposition 5.3 that $q_0$ takes on its maximum on $(-\infty, -1] \times M$ if $r \geq c_0$. Let $p$ a point where this is the case. Let $s_0$ be such that the support of $\chi_p$ lies in $[s_0 - 1, s_0 + 1] \times M$. Multiply both sides of (5.31) by $\chi_p G_p$, integrate the result, and integration by parts finds

$$q_0|_p \leq c_0 \int_{\mathbb{R} \times M} q_0 w_+ \chi_p \mathrm{dist}(p,\cdot)^{-2} + c_0 r \int_{\mathbb{R} \times M} w^2 \chi_p \mathrm{dist}(p,\cdot)^{-2} + c_0 \left( \int_{[s_0-1,s_0+1] \times M} q_0^2 \right)^{1/2}.$$

(5.51)

By virtue of (5.49), the right most integral in (5.51) is bounded by $c_0(1 + c)r$. The middle integral on the right hand side of (5.51) is also bounded by $c_0 r$. To make something of the left most integral, fix $\rho > 0$ and divide the integration domain into the part where the distance from p is greater than ρ, and the complementary part. The contribution from the former is no greater than an r-independent multple of $|\ln\rho|^{1/2}$ times the $L^2$ norm of $q_0$ over $[s_0 - 1, s_0 + 1] \times M$. The part where $\text{dist}(p, \cdot) \leq \rho$ contributes at most $c_0 \rho^2 q_0|_p$. This is because $q_0|_p$ is the maximum value for $q_0$. Take ρ so that this last contribution is equal to $\frac{1}{2} q_0|_p$ to obtain the bound asserted by the lemma.

### j) Proof of Proposition 5.3 for points in $(-\infty, -1] \times M$

To start, reintroduce on $\mathbb{R} \times M$ the closed 2-forms $\omega = e^{2s}(ds \wedge a + *a)$ and $F_A = ds \wedge \frac{\partial}{\partial s} A + *B_A$; and reintroduce $\partial_A$ to denote the covariant, exterior derivative on the trivial bundle over $\mathbb{R} \times M$ that is defined by the connection A.

Fix an increasing function $\sigma: [0, \infty) \to [0, 1]$ that has value 0 on $[0, \frac{1}{2}]$ and equals 1 on $[1, \infty)$. Write $\sigma'$ for the derivative of σ. Use $\sigma_\delta$ to denote the function on $\mathbb{R} \times M$ given by $\sigma(\delta^{-1}(1 - |\alpha|^2))$. Likewise, use $\sigma'_\delta$ to denote $\sigma'(\delta^{-1}(1 - |\alpha|^2))$. With all of this notation set, let

$$\wp_\delta = \delta^{-1} \sigma'_\delta \, \partial_A \alpha \wedge \partial_A \bar{\alpha} + \sigma_\delta F_A .$$

(5.52)

This is a closed 2-form on $\mathbb{R} \times M$. The $s \geq -1$ version of Proposition 5.3 implies that this form has support where $s < -1$ when r is large. It is also the case that $|\wp_\delta|$ is bounded on $\mathbb{R} \times M$ by $c_0 \delta^{-1} r$; this a consequence of Lemmas 1.6 and 1.7.

Granted what was just said, it follows that $|\omega \wedge \wp_\delta|$ is integrable on $\mathbb{R} \times M$, and thus so is the 4-form $\omega \wedge \wp_\delta$. Moreover, because $|\wp_\delta|$ is bounded and s is bounded from above on its support, an integration by parts can be employed to verify that the integral of $\omega \wedge \wp_\delta$ over $\mathbb{R} \times M$ is zero. The vanishing of the integral of $\omega \wedge \wp$ asserts that

$$\int_{\mathbb{R} \times M} e^{2s}(\delta^{-1} \sigma'_\delta (|\partial_A \alpha|^2 - |\bar\partial_A \alpha|^2) + re^{2T} \sigma_\delta (1 - |\alpha|^2 + |\beta|^2 - r^{-1} e^{-2T} \langle a, 2\chi\mu + \hat{\mu}\rangle)) = 0 .$$

(5.53)

Two auxillary results are needed in order to make something from (5.53). These are stated in the next lemma.

**Lemma 5.13**: *There exists $\kappa > 0$, and given $c \geq 0$ and $\delta > 0$, there exist constants $r_\delta \geq 1$ and $\kappa_\delta > 0$; and these have the following significance: Suppose that $r \geq r_\delta$ and that $(A, \psi)$ is a solution to (4.2) such that $|\frac{\partial}{\partial s} A| + |B_A| \leq c r e^{2T}$. Then*
- $|\bar\partial_A \alpha| \leq \kappa$.
- *Let $\Omega_\delta$ denote the integral over $\mathbb{R} \times M$ of $e^{2s} e^{2T} \sigma_\delta$ and let $\Omega'_\delta$ denote the integral over $\mathbb{R} \times M$ of $e^{2s} \sigma_{\delta/2}$. Then $\Omega'_\delta < \kappa_\delta \Omega_\delta$.*

This lemma is proved momentarily.

Granted that it is true, note first that its assumptions are implied by Lemma 5.1 and Proposition 5.2. This understood, then (5.15) implies that

$$(-\kappa_\delta \kappa^2 \delta^{-1} + r\delta(1 - r^{-1}\delta^{-1}c_0))\,\Omega_\delta < 0.$$

(5.54)

when $r > r_\delta$. It follows from this last inequality that there exists $r'_\delta \geq r_\delta$ such that $\Omega_\delta = 0$ when $r > r'_\delta$. Of course, this implies what is asserted by Proposition 5.3.

*Proof of Lemma 5.13*: The second equation in (4.2) equates $\bar{\partial}_A \alpha$ to a linear combination of covariant derivatives of $\beta$. This understood, the first item follows with an appropriate bound for $|\frac{\partial}{\partial s}\beta| + |\nabla'\beta|$. To obtain the latter and to prove the second item, fix a point $p \in \mathbb{R} \times M$ and a Gaussian coordinate system centered at p. Rescale the Gaussian coordinate as done in Subsection 5a so that the ball of radius $r^{-1/2}e^{-T(p)}$ in the original coordinates has radius 1 in the rescaled coordinates. Fix $R \geq 1$, $\varepsilon > 0$, and $k \in \{0, 1, 2, \ldots\}$ and reintroduce $r(R, \varepsilon, k)$ from this same subsection. As in Subsection 5a, if $r > r(R, \varepsilon, k)$, there exists $(A_0, (\alpha_0, 0))$ with the following two properties: First, $A_0$ and $\alpha_0$ are defined on all of $\mathbb{R}^4$ where they obey (5.2). Second, $(A_0, (\alpha_0, 0))$ has distance less than $\varepsilon$ in the $C^k$ topology on the ball of radius R about the origin in $\mathbb{R}^4$ to the pull-back of $(A, (\alpha, \beta))$ via the composition of the Gaussian coordinate chart map and the rescaling map.

Let $\beta_r$ denote the pull-back of $\beta$, now viewed as a complex number valued function on the ball of radius R in $\mathbb{R}^4$ by using the Gaussian coordinates to trivialize the bundle K. This understood, the pull-back of (5.6) implies an equation for $\beta_r$ in the ball of radius 1 about the origin in $\mathbb{R}^4$ that has the form

$$-\Delta\beta_r + \sum_{j=1,\ldots,4} \mathcal{R}_j \tfrac{\partial}{\partial x_j} \beta_r + \mathcal{R}(\beta_r) + r^{-1}\mathcal{R}' = 0,$$

(5.55)

where $\{|\mathcal{R}_j|\}_{j=1,\ldots,4}$, $|\mathcal{R}_0|$ and $|\mathcal{R}'|$ are all bounded by $c_0 c$. Granted that $\beta \leq c_0 r^{-1/2}$, it is a straightforward matter to use the Green's function for the Laplacian $\Delta$ to see that $|\frac{\partial}{\partial x_j}\beta_r|$ $\leq c_0 c\, r^{-1/2}$. The rescaled version of this last bound implies what is asserted in the first item of the lemma.

Now consider the second item in the proposition. Its proof requires a related assertion about solutions to (5.2).

**Lemma 5.14**: *Given $c > 0$ and $\delta \in (0, 1)$, there exits $R_\delta \geq 2$ and $\kappa_\delta \geq 1$ with the following significance: Suppose that $(A_0, \alpha_0)$ satisfies (5.2) on $\mathbb{R}^4$. Let V denote the volume of the subset of the ball about the origin in $\mathbb{R}^4$ with radius $R_\delta$ where $(1 - |\alpha|^2) > \delta$; and let V'* *denote the volume of the subset of the concentric ball of radius $\frac{1}{2} R_\delta$ where $\delta > (1 - |\alpha|^2) \geq \frac{1}{2}\delta$. Then $V' \leq \kappa_\delta V$.*

**Proof of Lemma 5.14**: Suppose that $\delta$ is fixed, as is $R \geq 2$. Let $(A_0, \alpha_0)$ obey the stated conditions. Suppose, for the sake of the argument, that $V' > 0$. There are two cases to consider. In the first case, there exists a point $x \in \mathbb{R}^4$ with $|x| \leq \frac{3}{4} R$ where $1 - |\alpha|^2 > 2\delta$. Under the stated conditions, there is an apriori upper bound for the norm of the covariant derivative of $\alpha$ at each point in $\mathbb{R}^4$. This bound depend on the constant c in (5.2); in any event, let $c_*$ denote a bound with $c_* \geq 1$. Then $1 - |\alpha|^2 > \delta$ in the ball of radius $\frac{1}{2}\delta c_*^{-1}$ centered at x, and so V is no less than $c_0 c_*^{-4}\delta^4$. Meanwhile, V´ is no greater than the volume of a radius $\frac{1}{2} R$ ball in $\mathbb{R}^4$, thus no greater than $c_0 R^4$. This implies that the conclusions of the lemma hold with $\kappa_\delta = c_0 c_*^{-4}\delta^{-4} R^4$.

The second possibility is that $(1 - |\alpha|^2) < 2\delta$ at all points in the ball of radius $\frac{3}{4} R$ ball centered on the origin in $\mathbb{R}^4$. To see what this implies, recall that there is a gauge transformation in this ball such that takes $\alpha$ to $e^{-u}$, with u a positive function. Moreover, as $u \leq \frac{1}{2} |\ln(1 - 2\delta)|$, there exists $\varepsilon_\delta$ such that $1 - e^{-2u} - \varepsilon_\delta u \geq 0$ in this radius $\frac{3}{4} R$ ball. The function u also obeys the equation $\Delta u = 1 - e^{-2u}$, and so this means that $\Delta u - \varepsilon_\delta u \geq 0$ in the radius $\frac{3}{4} R$ ball centered on the origin. It then follows using the comparison principle that

$$u|_x \leq \tfrac{1}{2} |\ln(1 - 2\delta)| \exp(-\varepsilon_\delta^{1/2}(\tfrac{3}{4} R - |x|))$$

(5.56)

at all points with $|x| \leq \frac{3}{4} R$. Note that this is not compatible with the assumption $V' \neq 0$ if $R \geq R_\delta = 8\varepsilon_\delta^{-1/2} \ln(1 - 2\delta)/\ln(1 - \tfrac{1}{2}\delta)$.

To continue the proof of Lemma 5.13, fix $R = R_\delta$ as given in Lemma 5.14, and fix a set, $\Lambda$, of balls in $\mathbb{R} \times M$ with the following properties: First, each ball in $\Lambda$ has radius $\tfrac{1}{4} Rr^{-1/2}$ and center where $\delta \geq (1 - |\alpha|^2) \geq \tfrac{1}{2}\delta$. Second, distinct balls from $\Lambda$ are disjoint. Third, $\Lambda$ is maximal amongst all sets that satify the first two criteria. For each ball $U_* \subset \Lambda$, let U´ denote the concentric ball with radius $\tfrac{1}{2} Rr^{-1/2}$, and let U denote the concentric ball with radius $Rr^{-1/2}$. When r is large (as determined solely by $\delta$), the set $\cup_\Lambda U'$ covers the set of points in $\mathbb{R} \times M$ where $\delta \geq (1 - |\alpha|^2) \geq \tfrac{1}{2}\delta$. This follows by virtue of the fact that $\Lambda$ is maximal. Note that there is an r, (A, $\psi$) and $\delta$ independent bound to the number of distinct versions of U that can simultaneously intersect. Granted this, it follows from Lemma 5.14 that there exists $r_\delta$ that is independent of (A, $\psi$) and is such that if $r \geq r_\delta$, then the volume of the set of points in $\mathbb{R} \times M$ where $\delta \geq 1 - |\alpha|^2 \geq \tfrac{1}{2}\delta$ is bounded by an r and (A,$\psi$) independent multiple of the volume of the subset in $\mathbb{R} \times M$ where $1 - |\alpha|^2 \geq \delta$.

**k) Proof of Proposition 2.9**

The key to the proof that $\mathcal{B}(r)$ is non-empty when r is large is the proof of the following assertion:

**Proposition 5.15**: *There exists $r_{I*} > r_I$ with the following significance: Suppose that $r \geq r_I$. If $\mathfrak{p} \in \mathcal{P}$ is sufficiently small, then the $(r, \mathfrak{g} = e_\mu + \mathfrak{p})$ version of the moduli space $\mathcal{M}_{\iota=1}(\mathfrak{c}, \mathfrak{c}_I)$ is empty when $\mathfrak{c}$ is a non-degenerate solution to (1.3).*

This proposition implies that $\delta(\mathfrak{c}_I + \mathfrak{w})$ can not vanish if $\mathfrak{w}$ lacks the generator $\mathfrak{c}_I$. Indeed, were this class to vanish, then there would exist a class $\mathfrak{z}$ of degree 1 such that $\delta\mathfrak{z} = \mathfrak{c}_I + \mathfrak{w}$ because the Seiberg-Witten Floer homology is assumed to vanish. This is impossible if there are no instantons with $s \to \infty$ limit equal to $\mathfrak{c}_I$.

*Proof of Proposition 5.15*: It is enough to prove that $\mathcal{M}_{\iota=1}(\cdot, \mathfrak{c}_I)$ is empty for the $\mathfrak{p} = 0$ case. Indeed, if this is true, then some very minor modifications to the compactness arguments from Chapter 16 in [KM1] forbid a non-empty version of $\mathcal{M}_{\iota=1}(\cdot, \mathfrak{c}_I)$ when $\mathfrak{p}$ is non-zero but has small norm. The proof of the $\mathfrak{p} = 0$ version of Proposition 5.15 has four steps.

Step 1: Suppose that $\mathfrak{c}$ and $\mathfrak{c}'$ are non-degenerate solutions to (1.3) and that the moduli space $\mathcal{M}_\iota(\mathfrak{c}, \mathfrak{c}') \neq \emptyset$. Let $\mathfrak{d} = (A, \psi) \in \mathcal{M}_\iota(\mathfrak{c}, \mathfrak{c}')$.

The arguments that lead to the conclusions of Lemma 5.1 work here to prove that $|\alpha|^2 \leq 1 + \kappa r^{-1}$ and $|\beta|^2 \leq \kappa r^{-1}(1 - |\alpha|^2) + \kappa r^{-2}$. Indeed, simply set $T = 1$ in the proof of Lemma 5.1 given in Section 5a and use Lemma 1.6 to control the large $|s|$ behavior of $|\alpha|$ and $|\beta|$.

Step 2: The first point to make is that $|\mathfrak{a}(\mathfrak{c}) - \mathfrak{a}(\mathfrak{c}')| \leq c_0 r^2$. Indeed, this follows from (1.6) with Proposition 1.10 since the degree of $\mathfrak{c}$ is $\iota$ less than the degree of $\mathfrak{c}'$. The next point is that (5.47) holds in this case, but now where $t \in \mathbb{R}$ is allowed. In particular, taking $t \to \infty$ finds that

$$\int_{\mathbb{R} \times M} (|\tfrac{\partial}{\partial s} A|^2 + |\mathfrak{B}_{(A,\psi)}|^2 + 2r(|\tfrac{\partial}{\partial s}\psi|^2 + |D_A \psi|^2)) \leq c_0 r^2$$

(5.57)

since $|\mathfrak{a}(\mathfrak{c}) - \mathfrak{a}(\mathfrak{c}')| \leq c_0 r^2$. It is also the case that (5.48) holds, where now s can have any value in $\mathbb{R}$. This understood, it follows from (5.57) that (5.49) holds with $c_* < c_0$ for any $s \in \mathbb{R}$.

The arguments used to prove Lemma 5.9 can be repeated once again to find r and $(A, \psi)$ independent constants $\kappa > 1$ and $\kappa_{1-3}$ such that when $r \geq \kappa$, then $q_0$ as defined in (5.50) obeys (5.51). Since $\lim_{|s|\to\infty} q_0 \leq c_0 r$, the argument used at the end of Section 5i to finish the proof of Proposition 5.2 work just as well in the present context to prove that $q_0 \leq c_0 r$ on the whole of $\mathbb{R} \times M$. Here is an immediate consequence:

**Lemma 5.16**: *Given $\iota \in \mathbb{Z}$, there is a constant $\kappa > 1$ with the following significance: Suppose that $\mathfrak{c}$ and $\mathfrak{c}'$ are non-degenerate solutions to (1.3) and that $(A, \psi) \in \mathcal{M}_\iota(\mathfrak{c}, \mathfrak{c}')$ is a solution to the $\mathfrak{r}$ and $\mathfrak{g} = \mathfrak{e}_\mu$ version of (1.10). Then $|\frac{\partial}{\partial s} A| + |B_A| \leq \kappa_* \mathfrak{r}$.*

Step 3: Now specialize to the case where $\mathfrak{c}' = \mathfrak{c}_I(\mathfrak{r})$. Fix $\delta > 0$ and reintroduce the function $\sigma_\delta$ from Section 5j and the 2-form $\wp_\delta$ from (5.52). There exists $\mathfrak{r}_\delta > \mathfrak{r}_I$ such that when $\mathfrak{r} > \mathfrak{r}_\delta$, then the function s has an upper bound on the support of $\wp_\delta$. This is the case for any choice of $\mathfrak{c}$ and for any $(A, \psi) \in \mathcal{M}(\mathfrak{c}, \mathfrak{c}_I(\mathfrak{r}))$. Let $\omega = e^{2s}(ds \wedge a + *a)$ as in Section 5j. Then $\omega \wedge \wp_\delta$ is integrable when $\mathfrak{r} > \mathfrak{r}_\delta$, and the integration by parts can be performed to see that its integral is zero. This implies that (5.53) holds with T and $\hat{\mu}$ set to zero and with $\chi$ identically equal to 1.

The following analog of what is claimed by Lemma 5.13 is needed to exploit this new version of (5.53).

**Lemma 5.17**: *Given $\iota \in \mathbb{Z}$ and $\delta > 0$, there exists $\mathfrak{r}_\delta > \mathfrak{r}_I$ and $\kappa_\delta > 0$ with the following significance: Suppose that $\mathfrak{r} \geq \mathfrak{r}_\delta$, that $\mathfrak{c}$ is a non-degenerate solution to (1.3), and that the pair $(A, \psi)$ is in the $\mathfrak{r}$ and $\mathfrak{g} = \mathfrak{e}_\mu$ version of $\mathcal{M}_\iota(\mathfrak{c}, \mathfrak{c}_I(\mathfrak{r}))$. Then*
- $|\bar{\partial}_A \alpha| \leq \kappa$.
- *Let $\Omega_\delta$ denote the integral over $\mathbb{R} \times M$ of $e^{2s}\sigma_\delta$ and let $\Omega'_\delta$ denote the integral over $\mathbb{R} \times M$ of $e^{2s}\sigma_{\delta/2}$. Then $\Omega'_\delta < \kappa_\delta \Omega_\delta$.*

*Proof of Lemma 5.17*: Given what is said by Lemma 5.14, the proof of this lemma is essentially identical to that of Lemma 5.13.

With Lemma 5.17 in hand, then (5.54) is also true in this case, as its conclusion, that $\Omega_\delta = 0$ when $\mathfrak{r} > \mathfrak{r}'_\delta$, where $\mathfrak{r}'_\delta$ is independent of $\mathfrak{c}$ and of $(A, \psi) \in \mathcal{M}_\iota(\mathfrak{c}, \mathfrak{c}_I)$.

Step 4: To complete the argument, fix some small, but positive $\delta$ and take $\iota = 1$ to define $\mathfrak{r}'_\delta$ as in the previous step. Suppose that $\mathfrak{r} > \mathfrak{r}_\delta$, that $\mathfrak{c} = (A, \psi = (\alpha, \beta))$ is a solution to (1.3) and that $\mathcal{M}_{\iota=1}(\mathfrak{c}, \mathfrak{c}_I(\mathfrak{r})) \neq \emptyset$. It then must be the case that $1 - |\alpha|^2 < \delta$ at all points in M. If $\delta$ is sufficiently small, this and Lemma 1.6 imply via the uniqueness assertion of Proposition 2.8 that $\mathfrak{c} = u \cdot \mathfrak{c}_I(\mathfrak{r})$ with u a smooth map from M to $S^1$. But this isn't possible as the spectral flow between $\mathfrak{c}$ and $\mathfrak{c}_I(\mathfrak{r})$ is assumed equal to 1 while the spectral flow between any gauge transform of $\mathfrak{c}_I(\mathfrak{r})$ and $\mathfrak{c}_I(\mathfrak{r})$ is an even number.

**References**

[ACH] C. Abbas, K. Cielebak and H. Hofer, *The Weinstein conjecture for planar contact structures in dimension three*, preprint arXiv:math.SG/0409355v2 March 2005.

[D] S. K. Donaldson. *Floer homology groups in Yang-Mills theory*, Volume 147 of Cambridge Tracts in Mathematics, Cambridge University Press, Cambridge, 2002. With the assistance of M. Furuta and D. Kotschick.

[E] Y. Eliashberg, *Contact 3-manifolds twenty years since J. Martinet's work*, Ann. Inst. Fourier, vol. **42** (1992), 165-192.

[H] H. Hofer, *Pseudoholomorphic curves in symplectizations with applications to the Weinstein conjecture in dimension three*, Invent. Math. **114** (1993), 515-563.

[JT] A. Jaffe and C. H. Taubes, Vortices and Monopoles, Birkhäuser, Boston 1980.

[Ka] T. Kato, Perturbation theory for linear operators, Springer-Verlag, Berlin-Heidelberg-New York, 1966.

[KM1] P. Kronheimer and T. Mrowka, Monopoles and Three-Manifolds, Cambridge University Press, to appear.

[KM2] P. Kronheimer and T. Mrowka, *Monopoles and contact structures*, Invent. Math. **130** (1997) 209-255.

[MR] T. Mrowka and Y. Rollin, *Legendrian knots and monopoles*, Algebr. Geom. Topol. **6** (2006), 1-69.

[MT] G. Meng and C. H. Taubes, *SW and Milnor torsion*, Math. Res. Lett. **3** (1996), 661-674.

[T1] C. H. Taubes, *The Seiberg-Witten equations and the Weinstein conjecture*, preprint arXiv:mathSG/0611007.

[T2] C. H. Taubes, *Asymptotic spectral flow for Dirac operators*, preprint arXiv:math.DG/0612126.

[T3] C. H. Taubes, *SW =>Gr: From Seiberg-Witten equations to pseudo-holomorphic curves*, in Seiberg-Witten and Gromov Invariants for Symplectic 4-manifolds by C. H. Taubes, International Press, Somerville MA, 2005.
78